\documentclass[reqno,11pt]{amsart}
\usepackage{geometry}
\usepackage[numbers,sort&compress]{natbib}
\usepackage{mathtools,amssymb,amsthm,mathrsfs,color,lineno,paralist,graphicx,float}
\usepackage[colorlinks,
linkcolor=red,
anchorcolor=green,
citecolor=blue,
]{hyperref}

\usepackage[T1]{fontenc}
\usepackage[utf8]{inputenc}

\usepackage{tikz}
\usetikzlibrary{positioning}
\usepackage{calc}
\definecolor{bleu1}{RGB}{0,57,128}
\def\bleu1{\color{bleu1}}

\usepackage{etoolbox}
\patchcmd{\section}{\normalfont}{\normalfont \bleu1}{}{}
\patchcmd{\subsection}{\normalfont}{\normalfont \bleu1}{}{}
\patchcmd{\subsubsection}{\normalfont}{\normalfont \bleu1}{}{}
\renewcommand{\proofname}{\it \bleu1 Proof}

\def\a{\alpha}

\def\e{\varepsilon}

\let\newpf\proof \let\proof\relax 
\newenvironment{pf}{\newpf[\proofname]}{\qed\endtrivlist}

\newcommand{\ba}{\overline{A}}

\def\be{\begin{equation}}
\def\ee{\end{equation}}

\def\ba{{\begin{align}}}
\def\ea{{\end{align}}}

\def\bm{\begin{matrix}}
\def\em{\end{matrix}}

\def\a{{\alpha}}

\def\0{{\mathbf 0}}

\newtheorem{Theorem}{Theorem}[section]
\newtheorem{Lemma}{Lemma}[section]
\newtheorem{Proposition}{Proposition}[section]
\newtheorem{Corollary}{Corollary}[section]
\newtheorem{Remark}{Remark}[section]
\newtheorem{Example}{Example}[section]
\newtheorem{Definition}{Definition}[section]

\numberwithin{equation}{section}

\theoremstyle{definition}

\renewcommand{\mod}{\operatorname{mod}}

\newcommand{\C}{{\mathbb C}}

\newcommand{\N}{{\mathbb N}}
\newcommand{\Q}{{\mathbb Q}}
\newcommand{\R}{{\mathbb R}}
\newcommand{\T}{{\mathbb T}}

\newcommand{\Z}{{\mathbb Z}}

\def\B0{{\bold{0}}}

\catcode`\@=12

\def\Empty{}
\newcommand\oplabel[1]{
  \def\OpArg{#1} \ifx \OpArg\Empty {} \else
    \label{#1}
  \fi}

\newcommand{\comm}[1]{}
\newcommand{\comment}[1]{}

\begin{document}
\title[]{Kotani theory, Puig's argument, and stability of The Ten Martini Problem}

\author {Lingrui Ge}
\address{Beijing International Center for Mathematical Research, Peking University, Beijing, China
	} \email{gelingrui@bicmr.pku.edu.cn}
\author{Svetlana Jitomirskaya}
\address{
	Department of Mathematics, University of California, Irvine CA, 92717} \email{zhitomi@uci.edu}
\author{Jiangong You}
\address{
	Chern Institute of Mathematics and LPMC, Nankai University, Tianjin 300071, China} \email{jyou@nankai.edu.cn}

      \begin{abstract}
 We solve the ten martini problem (Cantor spectrum with no condition on
irrational frequencies, previously only established for the almost
Mathieu) for a large class of one-frequency
quasiperiodic operators, including nonperturbative analytic neighborhoods of
several popular explicit families. 
The proof is based on the structural analysis of dual cocycles as
introduced in \cite{gjyz}. 
As a part of the proof, we develop several general ingredients of  independent interest:
Kotani theory, for a class of finite-range
operators over general minimal underlying dynamics, making the first
step towards and providing a partial solution of the Kotani-Simon problem,  simplicity of
point spectrum for the same class, and the
all-frequency version of Puig's argument.

\end{abstract}

\maketitle
\tableofcontents
\section{Introduction}

The Hofstadter butterfly \cite{hof}, a plot of the  band spectra of almost
Mathieu operators
\begin{align}\label{amo}
(H_{\lambda,\alpha,x}u)_n=u_{n+1}+u_{n-1}+2\lambda\cos2\pi(x+n\alpha)u_n,
\end{align}
at rational frequencies $\alpha$, has become a pictorial symbol
of the field of quasiperiodic operators. It is visually clear
from this plot that for {\it all irrational frequencies} the spectrum must be
a Cantor set, a statement that has been dubbed the ten martini problem
by Barry Simon \cite{barry} after an 1981 offer of Mark Kac \cite{kac}. The problem itself is
considered iconic in the field of quasiperiodic operators, its final
solution \cite{aj} requiring a combination of many ideas and techniques and
significant ingenuity. The proof in \cite{aj} used the specific nature of the
almost Mathieu operator in several key ways and was based on different
approaches at the Diophantine and Liouville sides that miraculously
met at the middle \cite{solving}, enough so that in a field driven by
bold conjectures (e.g. \cite{sim15,simXXI}) a conjecture that the same
statement could hold for other operators \eqref{sch} has never even
been explicitly made.

Indeed, while there were several Cantor spectrum results for
operators \eqref{sch} with analytic $v$, other than for the almost
Mathieu family, all required various (often implicit) conditions on frequencies,
among other unnatural restrictions. At the same time, the physics
nature and relevance of the almost Mathieu family strongly
suggest that the ten martini problem has to be robust and hold at least in the
entire analytic neighborhood of \eqref{amo}.  Here we prove for the first time the {\it robust ten
martini problem}, for an open set of analytic one-frequency operators \eqref{sch},
by developing a method that does not rely neither on the arithmetics
nor on the almost Mathieu specifics.

One-frequency analytic quasiperiodic Schr\"odinger operators on
$\ell^2(\Z)$ are given by \eqref{amo} with $2\lambda\cos $ replaced by
a 1-periodic non-constant real analytic function $v,$ that is
\begin{align}\label{sch}
(H_{v,\alpha,x}u)_n=u_{n+1}+u_{n-1}+ v(x+n\alpha)u_n,\ \ n\in\Z,
\end{align}
where $\alpha\in\R\backslash\Q$ and $x\in\R$ are parameters (called
the {\it frequency} and the {\it phase} respectively). Their theory has been developed extensively (see \cite{bbook,damanik,youcongr,jitcongr} for more recent surveys).  


The almost Mathieu operator \eqref{amo}  (AMO) is the
central/prototypical model, lying both at the physics origin and the
center of current physics interest of the
field \cite{Peierls, harper,R,aos,oaag,h,ntw,nichen}, as well as driving many
of the mathematical developments. The latter is, at least historically
due to Barry Simon's problems \cite{sim15,simXXI} prominently featuring several
almost Mathieu questions. This has remarkably lead to all of them being solved,
and then many new ones appearing.

The study of general operators \eqref{sch} with analytic
$v$ has long been developed in the perturbative regime, with the key
highlights in \cite{ds,fsw,sin} and especially Eliasson \cite{Eli92,Eli97}. 
The nonperturbative analysis has taken off after the work of Bourgain and collaborators
(see \cite{bbook}), who significantly developed theory of operators
\eqref{sch}, especially in
the regime of positive Lyapunov exponents, an important catalyst to
these developments being again the almost Mathieu result \cite{j}. The
development of nonperturbative/Liouvillean KAM \cite{hy,afk} and quantitative
reducibility (see \cite{youcongr} and references therein) has led to
many strong results in the (almost) reducibility  regime.

Avila's global theory \cite{avila0} of operators \eqref{sch},  based on the analysis of
complexified Lyapunov exponents, has brought new vision and
understanding, in particular, introducing a simple yet fundamental concept of acceleration,
as an important feature that allows to divide the spectrum into more
manageable subsets. With a rough division of the spectrum  into subcritical, critical,
and supercritical energies, Avila showed that critical ones are very
rare in a strong sense \cite{avila0}, while almost reducibility  becomes a
corollary of subcriticality \cite{arc1,arc2}, see also \cite{ge} for a
different proof for the Diophantine case.

Despite all these remarkable advances, many major results that do not
require unnecessary and/or non-explicit parameter exclusion
e.g. \cite{jliu1,jliu2,liuresonant,aj,ak, alsz, last, dryAYZ}, still exist only for the almost Mathieu operators and have heavily used
several different almost Mathieu specifics. Of those, the ten martini
problem particularly stands out \footnote{Along with its dry version.}. First, while the statement is, by
design, about  {\it all} irrational $\alpha,$ historically, the proofs 
developed very different arguments depending on the arithmetic
properties of $\alpha,$
all using the specific features of \eqref{amo}.   Moreover, it is the only frequency-universal
almost-Mathieu statement that has not been extended to {\it any} other
operator \eqref{sch}. Indeed, absolutely continuous spectrum and
absolute continuity of the density of states
are now known for all $\alpha$ and all subcritical operators
\eqref{sch}, an open set, by the proof of the almost reducibility
conjecture \cite{arc1,arc2}, while the ten-martini proofs have
remained almost-Mathieu only.


Here we present a method of proof of Cantor spectrum that does not depend on
the almost Mathieu symmetry, self-duality or low degree of the
potential, treats all
irrational frequencies simultaneously, and works for a large open set
of analytic potentials.

Our method is based on the  recently developed quantitative global
theory \cite{gjyz}: the Aubry duality based approach 
to Avila's global theory \cite{avila0}, that, in particular,  linked
the dynamics of dual cocycles to properties of 
Lyapunov exponents of complexified  Schr\"odinger cocycles and spectral features of \eqref{sch}.

  Lyapunov exponent of complexified  Schr\"odinger cocycles are defined as
\begin{align}\label{multiergodicsch}
L_\e(E)=\lim\limits_{n\rightarrow\infty}\frac{1}{n}\int_\T\ln \|S_E^v(x+i\e+(n-1)\alpha)\cdots S_E^v (x+i\e)\|dx
\end{align}
where
\begin{equation}\label{S}
S_E^v (x)=\begin{pmatrix}E-v(x)&-1\\ 1&0\end{pmatrix}.
\end{equation}
See Section \ref{pre} for details. Avila showed \cite{avila0} that $L_\e(E)$ is an even convex piecewise
affine function with integer slopes. Moreover, he argued in
\cite{avila0} that the exact value of the first slope is the most
important quantity, leading to the definition of the {\it acceleration},
$$
\omega(E)=\lim\limits_{\e\rightarrow 0^+}\frac{L_\e(E)-L_0(E)}{2\pi\e}.
$$
and showed, in particular, that $L$ is stratified analytic on the
spectrum, with strata defined by the values of acceleration. At the
same time, acceleration is not locally constant with
respect to analytic perturbations, as is easily illustrated by the
almost Mathieu family.

\begin{figure}[htbp]\label{fig}
\centering
\begin{minipage}[t]{0.3\linewidth}
\centering
\includegraphics[width=1.2\linewidth]{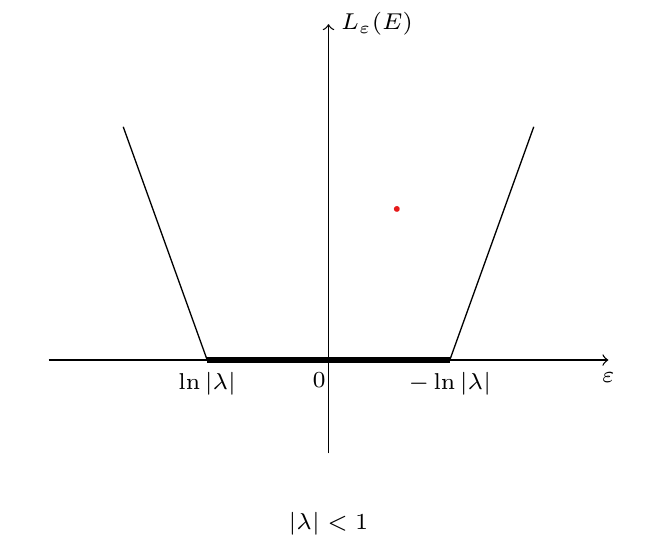}
\end{minipage}
\begin{minipage}[t]{0.3\linewidth}
\centering
\includegraphics[width=1.2\linewidth]{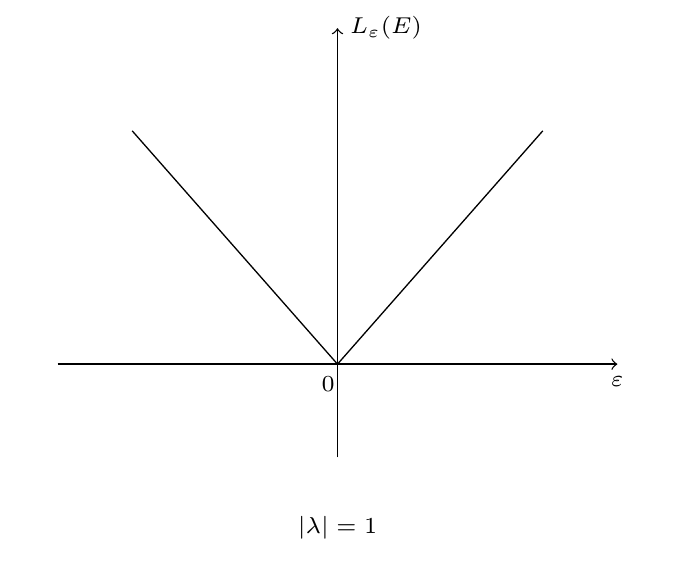}
\end{minipage}
\begin{minipage}[t]{0.3\linewidth}
\centering
\includegraphics[width=1.2\linewidth]{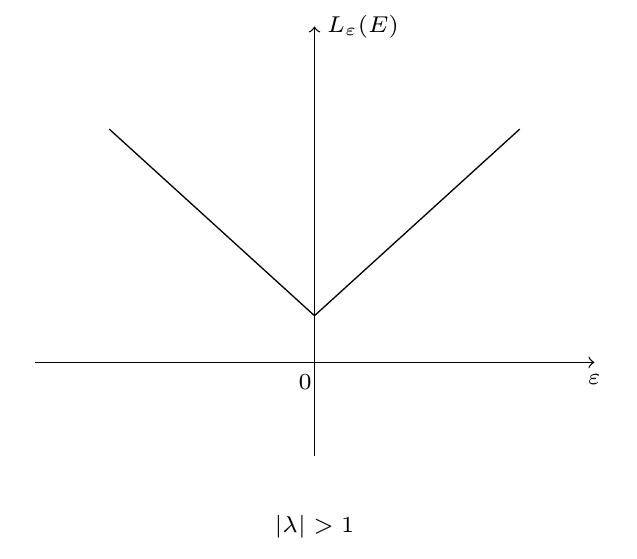}
\end{minipage}
\end{figure}

Indeed, for $E$ in the spectrum of the almost Mathieu operators \eqref{amo},
$L_\e(E),$ as a function of $\e,$ is particularly simple, allowing only
three fundamental possibilities, as in Figure \ref{fig}, with, in
particular, acceleration on the spectrum changing from $0$ to $1$ at
$\lambda=1.$ However, these
pictures also illustrate that the almost Mathieu acceleration is always bounded by $1,$ and it is this
feature that turns out to be both very robust and important for our proof of Cantor spectrum.

This paper is second in the series started with
\cite{gjz} where we first embarked on a project to extend the almost Mathieu facts to its
analytic neighborhood. The fact that some {\it  supercritical} almost Mathieu proofs can be extended
to the case of  acceleration $1$ has been clear ever since the
introduction of acceleration in \cite{avila0} where Avila, in particular,
showed that, just like for the almsot Mathieu, the Lyapunov exponent
restricted to the spectrum is an analytic function in this case. Moreover, Avila's proof
essentially showed that, for energies with acceleration $1$, traces of transfer-matrices 
(i.e., determinants of block-restrictions with
periodic boundary conditions) of size $q_n$ effectively behave like
trigonometric polynomials of
degree $q_n$, which they are for the almost Mathieu, a crucial feature
for localization proofs in \cite{j,aj,liu1,
  jliu1,jliu2,liuresonant} that does lead to a sharp spectral
transition result for supercritical type I operators
 \cite{gj}. At the same time, the AMO
ten martini proof is designed for the subctitical
case \footnote{utilizing self-duality of the family \eqref{amo} to
  obtain the result also for the
supercritical regime.}. While, in some sense, all subcritical energies
are alike \cite{arc1,arc2,ge}, in many other ways they are not, and the AMO ten martini proof did require
several specific features of almost Mathieu operators. However,  the concept of acceleration did not allow to distinguish what
makes the {\it subcritical} almost Mathieu special.
Here  we introduce a new concept that achieves that, allowing to
divide {\it both} sub and (super)critical parts of the spectrum into more
manageable sets, and at the same time leading to a stable property
encompassing both the sub and (super)critical regimes.

\begin{Definition}[T-acceleration]\label{gege1}
{\rm The {\it T-acceleration} is defined by
$$
\bar{\omega}(E)=\lim\limits_{\e\rightarrow \e_1^+}\frac{L_\e(E)-L_{\e_1}(E)}{\e-\e_1}
$$
where $0\le\e_1<\infty$ is the first turning point \footnote{It is an
  easy corollary of  the results of \cite{gjyz} that for
non-constant trigonometric polynomial $v$ we always have $\e_1<\infty$.} of the piecewise affine
function $L_\e(E)$.  If there is no
turning point, we set $\bar{\omega}(E)=1$.}
\end{Definition}
\begin{Remark}
{\rm Obviously, $\omega(E)\leq \bar{\omega}(E)$ for any $E\in\R$ and the equality holds if and only if $\omega(E)>0$. In particular, for the almost Mathieu operator, $\bar{\omega}(E)=1$ for all $E$ in the spectrum.}
\end{Remark}

\begin{Definition}[Type I]
  {\rm We say $E$ is a {\it type I energy} for operator
    $H_{v,\alpha,x}$ if
  $\bar{\omega}(E)=1$. We say $H_{v,\alpha,x}$ is a {\it type I operator,}
  if every $E$ in the spectrum of
  $H_{v,\alpha,x}$ is type I.
  }
\end{Definition}

Just like the acceleration, by continuity of
the Lyapunov exponent \cite{bj} and convexity, T-acceleration is upper-semincontinuous in
$\R\backslash \Q\times C^{\omega}(\R/\Z, SL(2,\C))$. Thus, since we always
have $\bar{\omega}(E)\ge 1,$  unlike the property of having
acceleration $1$, the property of T-acceleration being equal to $1$ is stable, and the set of type I operators
includes, in particular, the following sets, for any  1-periodic real analytic function $f:$ 
\begin{Example}\label{exa1}{\rm The almost Mathieu operator and its
    analytic perturbations, i.e., $v=2\lambda(\cos2\pi(x)+\delta
    f(x))$ where
    $|\delta|<\delta(\lambda, \|f\|_0)$ 
  }
\end{Example}
\begin{Example}\label{exa2}{\rm The GPS model in \cite{gps} and its
    analytic perturbations, i.e.,
    $v(x)=\frac{2\cos2\pi(\theta)}{1-b\cos2\pi(\theta)}+\delta f(x)$
    with $b\in(-1,1)$, with $|\delta|<\delta(a,b, \|f\|_0)$ 
  }
\end{Example}
\begin{Example}\label{exa3}{\rm The supercritical generalized Harper's
    model of \cite{hk,se} and its analytic perturbations, i.e.,
    $v(x)=2a\cos2\pi(x)+2b\cos4\pi(x)+\delta f(x)$ with $b\in(-1,1)$
    restricted to the positive Lyapunov exponent regime where
    $|\delta|<\delta(a,b, \|f\|_0)$
  }
\end{Example}

While type I operators generally have none of the other nice almost
Mathieu features such as symmetry, self-duality, or low-degree, we prove

\begin{Theorem}\label{main11}
Type I operators \eqref{sch} with  non-constant even trigonometric potentials $v$ have  Cantor spectrum for all $\alpha\in \R\backslash\Q.$ 
\end{Theorem}

Moreover, the result can actually be localized to the set of energies with
$\bar{\omega}(E)=1.$ Let $\Sigma^1_{v,\alpha}=\{E\in
\Sigma_{v,\alpha}: \bar{\omega}(E)=1\}$ where $\Sigma_{v,\alpha}$ is
the spectrum of $H_{v,\alpha,x}$ \footnote{It  does not depend on $x$
  when $\alpha$ is irrational \cite{as}.}. Theorem \ref{main11} is a
direct corollary of 
\begin{Theorem}\label{main12}
For any $\alpha\in\R\backslash\Q$ and $v$ a non-constant even trigonometric polynomial, $\Sigma_{v,\alpha}^1$ is a Cantor set.
\end{Theorem}

\begin{Remark}
{\rm Neither the assumption that  $v$ is even nor that it is a
  trigonometric polynomial are essential. Both theorems in fact hold
  for general real analytic $v.$ However,  rather technical
  quantitative arguments required both for  removing the evenness
  assumption and for extending the proof to the general analytic case
are unrelated to the core novel ideas
  of this paper, and would lengthen it unnecessarily. To streamline
  the presentation we chose to present here our key ideas for this
  technically simplified case that contains all the important
  features, 
 while the extensions to the general analytic case and the removal of
 the evenness requirement will be included instead in the forthcoming paper \cite{gjy2}.}
\end{Remark}
\begin{Remark}
{\rm To the best of our knowledge, Theorem \ref{main11} is not only the first Cantor spectrum
  result for non-almost Mathieu operators \eqref{sch} without an assumption on
  the frequency, but also the first such result without an
  assumption on the Lyapunov exponent.}
 \end{Remark}
 \begin{Remark}
{\rm Cantor spectrum results for the almost Mathieu operators have had
  a long history \cite{bs,sin,hs,cey,last,ak,Puig} prior to the proof of the ten martini
  problem  in \cite{aj}. For non-almost Mathieu operators \eqref{sch}
  with analytic $v,$
  the only existing results were for either (unspecified) typical  $v$
  in the zero Lyapunov exponents regime \cite{avilajams} \footnote{or 
  \cite{aj1} that provides the dry
  version but, again, only for typical $v$ and with a further
  Diophantine/smallness restriction.} or (unspecified) typical $\alpha$ in
  the regime of positive $L(E)$, proved in a combination of very
  technically complicated 
  \cite{gs1,gs2}, or in the perturbative regime \footnote{So, again,
    not for all irrational $\alpha$ for any given $v$.} \cite{sin,wz} for
  $\cos$-type $v.$}
\end{Remark}
\begin{Remark}
{\rm Contrary to the almost Mathieu operator, $L_\e(E)$ of general  type I
  operators may have many turning points, see Fig \ref{fig2}
  .}
\begin{figure}[htbp]\label{fig2}
\centering
\begin{minipage}[t]{0.3\linewidth}
\centering
\includegraphics[width=0.944\linewidth]{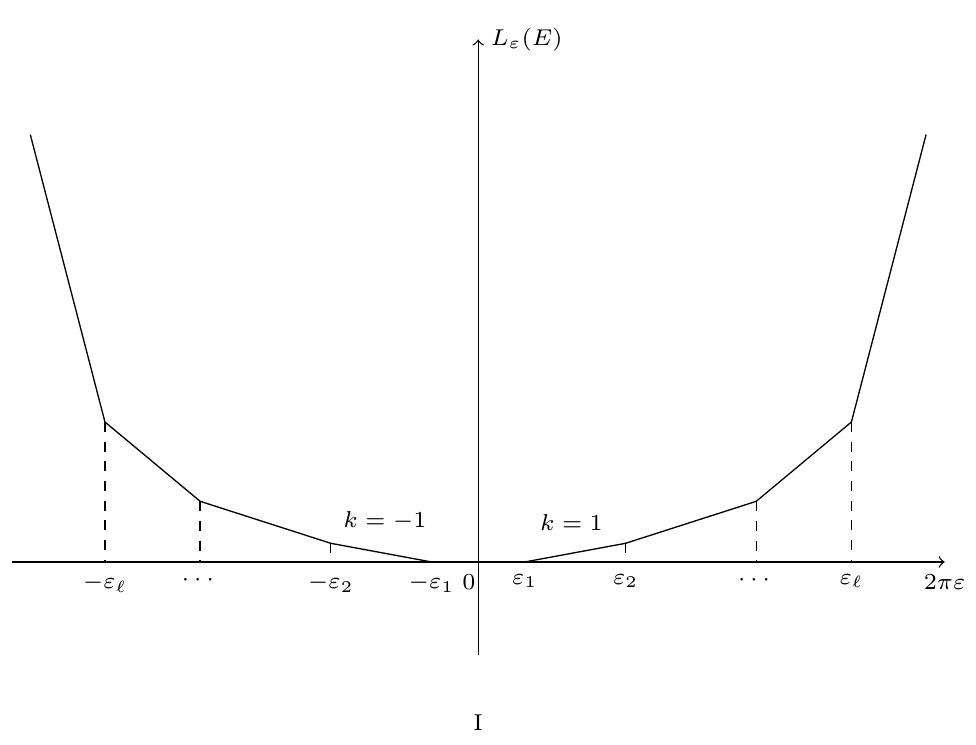}
\end{minipage}
\begin{minipage}[t]{0.3\linewidth}
\centering
\includegraphics[width=0.944\linewidth]{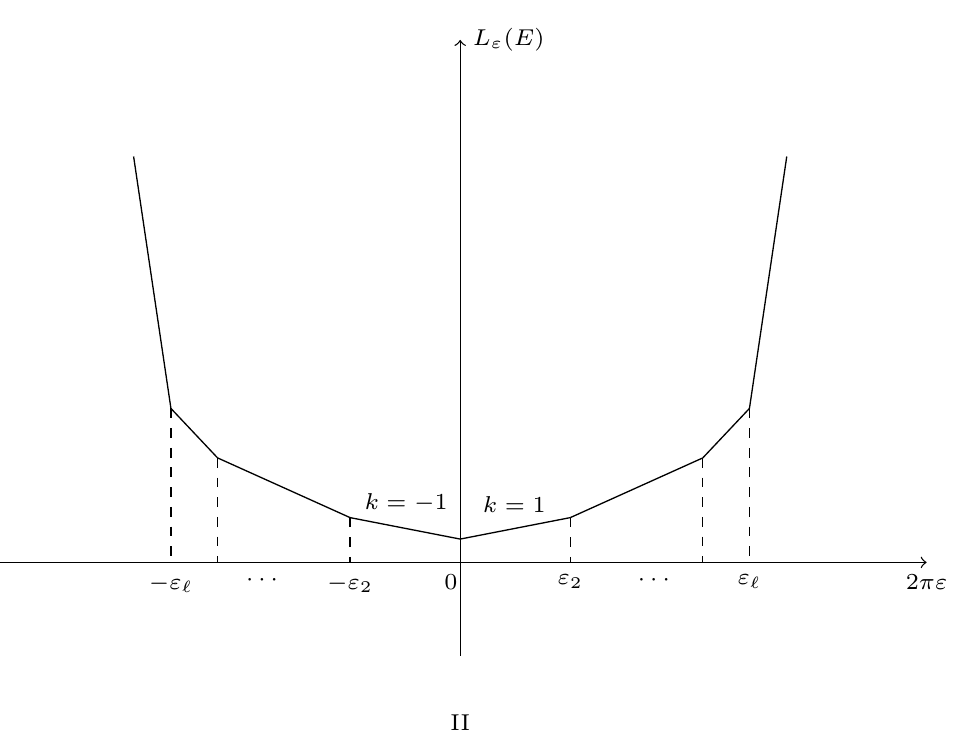}
\end{minipage}
\begin{minipage}[t]{0.3\linewidth}
\centering
\includegraphics[width=0.944\linewidth]{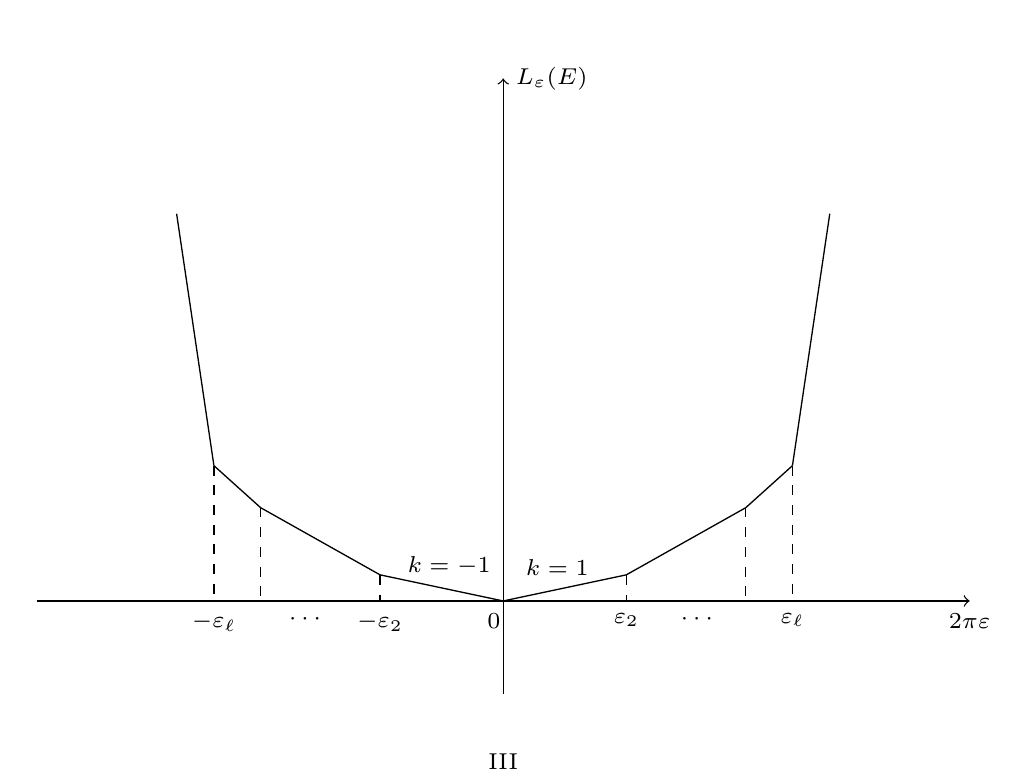}
\end{minipage}
\end{figure}
\end{Remark}

Besides settling the ten martini problem for all type I operators/
parts of the spectrum, our proof contains several
ingredients that we believe are of independent interest.

The first one is a partial solution of the Kotani-Simon problem
\cite{ks}.

Let $(\Omega,T)$ be an ergodic dynamical system, and $f:\Omega\to\R.$
A corresponding family of 
ergodic Schr\"odinger operators $H_{f,\omega}$,  $\omega \in\Omega,$ on
$\ell^2(\Z)$ is given by
\begin{align}\label{ergs}
(H_{f,\omega}u)_n=u_{n+1}+u_{n-1}+ f(T^n\omega)u_n,\ \ n\in\Z,
\end{align}
so that \eqref{sch} is a particular case, with $\Omega=\T$ and
$T\omega=\omega+\alpha.$ The transfer-matrix cocycles and Lyapunov
exponents for $H_{f,\omega}$ are defined in a similar way, see Section
\ref{3.2}.

Celebrated Kotani theory \cite{kot,sim83} proves that
ergodic Schr\"odinger operators with Lyapunov exponents vanishing on a
set of positive measure are deterministic, allowing, in particular,
the holomorphic extension of the
$m$-function through the interval with zero Lyapunov exponents. The
latter has been crucially used in the ten martini proof of
\cite{aj}. For general type I operators, our duality based argument
requires dealing with $Sp(2d,\R)$ cocycles for $d>1,$ for which this  result is
not available.

Indeed, Kotani theory has been extended by Kotani-Simon \cite{ks} to
Jacobi matrices on the strip, or matrix-valued Schr\"odinger operators, with corresponding  $Sp(2d,\R)$
transfer-matrix cocycles. However, the key result on the
reflectionness of $M$-matrix  \footnote{See Section \ref{6.2} for the
  definition.}, enabling, in particular, the holomorphic extension
described above, was proved in \cite{ks} only under the condition that {\it all
 } Lyapunov exponents are zero. At the same time, it was conjectured
 in \cite{ks} that under certain additional assumptions at least a
 partial version of the result
 should only require vanishing of {\it some
 } Lyapunov exponent \footnote{Additional assumptions/partial form are needed as
   demonstrated by an example of a model with decoupled potentials in
   \cite{ks}.}. Over the years, there have been some extensions of
 \cite{ks}, e.g.\cite{damanik1,oliviera}, and, most notably, \cite{xu},
 but the above problem stubbornly defied progress.

 Here we, for the first time, find and prove the desired partial reflectionless (therefore making the corresponding Kotani theory complete) in a
 situation where some
 Lyapunov exponents are positive.
 
 Let $\Omega$ be a compact metric space, and $T:\Omega\rightarrow\Omega$
 be a minimal homeomorphism \footnote{That is
   $\overline{\{T^n\omega\}_{n\in\Z}}=\Omega$ for any
   $\omega\in\Omega$.}.  For continuous $f:\Omega\rightarrow\R,$
 {minimal} finite-range operators on $\ell^2(\Z)$ are defined by
\begin{equation}\label{finite operator11}
(L_{f,\omega}u)(n)=\sum\limits_{k=-d}^d a_ku(n+k)+f(T^n\omega)u(n),\ \ n\in\Z,
\end{equation}
where $a_k=a_{-k}$ is a real sequence. As usual, the eigenvalue
equation $L_{f,\omega}u=Eu$ defines a  complex symplectic cocycle
$(T,L_E^f)$ \footnote{with respect to $S$ defined in \eqref{sdef}.} (see \eqref{1111}) and we denote its non-negative
Lyapunov exponents by
$\{L^i_f(E)\}_{i=1}^d$ (see Section \ref{3.2} for the definitions). Let
$\Sigma_f$ be the spectrum of $L_{f,\omega}$.

\begin{Definition}\label{defph2}
{\rm We  say that a complex symplectic {\it cocycle} is {\it $PH2$} if it is
{\it partially hyperbolic with two-dimensional center}. We say that an {\it energy $E$} is {\it $PH2$} for operator $L_{f,\omega}$ if the cocycle
$(T,L_E^f)$ is $PH2$. We
say that {\it operator $L_{f,\omega}$}  is {\it $PH2$} if every $E\in
\Sigma_f$ is $PH2$.}
\end{Definition}

In other words, $L_{f,\omega}$ is $PH2$ if

\begin{enumerate}
\item{$L^1_f(E)\geq \cdots\geq L^{d-1}_f(E)>L_f^d(E)$ for all $E\in\Sigma_f$;}
\item $(T,L_{E}^f)$ is $(d-1)$ and $(d+1)$-dominated for all
  $E\in\Sigma_f$ \footnote{See Section \ref{3.3} for the definition of domination.}.
\end{enumerate}

The results of \cite{ks} imply that the $M$ matrix (defined in Section
\ref{6.2}) is reflectionless on $\Sigma_f^0=\{E: L_f^1(E)=0\}$.


Here we solve the Kotani-Simon problem for all $PH2$ operators, thus a
class allowing positive Lyapunov exponents.


Let
$$
Sp_{2d\times 2}(\R)=\{F\in M_{2d\times 2}(\R): F^*S F=J\}, \ \ J=\begin{pmatrix}0&1\\ -1&0\end{pmatrix}.
$$
A crucial part of our proof of the Cantor spectrum is the following
general result on
{\it partial reflectionless.}
\begin{Theorem}\label{L2 reducibility}
For $PH2$ cocycles $(T,L_{E}^f)$ with minimal $T,$ for almost every $E$ in $\{E: L_f^d(E)=0\}$, there exist $H_E\in L^2(\Omega,Sp_{2d\times2}(\R))$ and $R_E(\omega)\in SO(2,\R)$ such that
$$
L_E^f(\omega)H_E(\omega)=H_E(T\omega)R_E(\omega).
$$
\end{Theorem}
\begin{Remark}
{\rm The theorem as stated is fully sufficient for our current
  purposes, however the simplicity of $L_f^d(E)$ (or, equivalently, the
  two-dimensionality of the center), while used substantially in the proof, is not essential, and this
  condition will be removed in the upcoming work of the first author
  \cite{g}.  What we currently see as crucially important is that the
  cocycle $(T,L_E^f)$ is partially hyperbolic. It is an interesting
  question whether this assumption is necessary for the result of Theorem \ref{L2 reducibility}.}
\end{Remark}
We will actually use Theorem \ref{L2 reducibility} through the
following corollary, also of independent interest
\begin{Corollary}\label{C0reducibility}
For $PH2$ cocycles $(T,L_{E}^f)$ with minimal $T,$  If $L^d_f(E)=0$ on an interval $I\subset \R$, there exist   $H_E\in C^0(\Omega,Sp_{2d\times2}(\R))$, and $R_E\in C^0(\Omega,SO(2,\R))$, depending analytically on $E\in I$ such that
$$
L_E^f(\omega)H_E(\omega)=H_E(T\omega)R_E(\omega).
$$
\end{Corollary}
\begin{Remark}
{\rm Corollary \ref{C0reducibility} implies in a standard way that if
  $\Sigma_f^0$ contains an open interval $I$, then $L_{f,\omega}$ has
   purely absolutely continuous spectrum on $I$ for any $\omega\in\Omega$.}
\end{Remark}

Another general  ingredient is a criterion of simplicity of point
spectrum for long-range operators. It is well known, by an easy
Wronskian argument, that point spectra of second-difference operators
are always simple. Certainly, the argument breaks down for
higher-difference operators, only implying absence of point spectra of
correspondingly high multiplicity. It turns out, however, that for
minimal $PH2$
operators
$L_{f,\omega},$ 
point spectrum is always simple, for any $d$, thus $PH2$ property implies certain {\it essential
  second-differenceness} of these higher-difference operators.
 \begin{Theorem}\label{tsimp}
 
 Minimal $PH2$ operators $L_{f,\omega}$ with continuous $f:\Omega\rightarrow \R$ have simple point spectrum for any $\omega\in\Omega$.
 \end{Theorem}

This allows an extension of Puig's argument originally designed for
the almost Mathieu operators. Puig showed \cite{puig} that Schr\"odinger cocycles
associated with almost Mathieu operators cannot be reduced to the
identity. This was a key element in his proof of the Cantor spectrum
for Diophantine $\alpha$. The argument itself is almost Mathieu
specific precisely because it is based on simplicity of point spectrum of
the dual operator, which, for the almost Mathieu, is again in the
almost Mathieu family. For general operators \eqref{sch}, a dual
operator is
defined by 
\begin{align}\label{fi}
(L_{v,\alpha,\theta}u)_n=\sum\limits_{k=-\infty}^{\infty} \hat{v}_k u_{n+k}+2\cos2\pi(\theta+n\alpha)u_n,\ \ n\in\Z,
\end{align}
where $\hat{v}_k$ is the $k$-th Fourier coefficient of $v$. In
particular, if $v$ is a
trigonometric polynomial potential of degree $d,$ operator
$L_{v,\alpha,\theta}$ is finite-range and of the form \eqref{finite
  operator11}, so for every $E$ we can define a corresponding
complex-symplectic cocycle. Slightly abusing the language, we will
call such cocycle the {\it dual cocycle} of $(\alpha,S_E^v)$ and/or of
the operator \eqref{sch}.

Theorem \ref{tsimp} immediately leads to
the following generalized Puig's argument

\begin{Theorem}\label{gjy-puig}

A Schr\"odinger cocycle $(\alpha,S_E^v)$ whose dual is $PH2$,
cannot be (analytically) reduced to the identity, i.e. there does not exist $B\in C^{\omega}(\T,SL(2,\R))$ such that
$$
B^{-1}(x+\alpha)S_E^v(x)B(x)=Id.
$$

\end{Theorem}

The reducibility-based arguments however require Diophantine
conditions, thus cannot work within an {\it all $\alpha$} argument. Here we
develop a uniform in $\alpha$ scheme, by replacing localized
eigenfunctions and reducibility to the identity in
a Puig-type argument with {\it almost localized eigenfunctions} and
{\it rotations reducibility }. We prove

\begin{Theorem}\label{rot1}

A Schr\"odinger cocycle $(\alpha,S_E^v)$ whose dual is $PH2,$
cannot be (analytically) reduced to  a rotation with zero
rotation number, i.e. there does not exist $B\in C^{\omega}(\T,SL(2,\R))$ and $\psi\in C^\omega(\T,\R)$ with $\int_\T \psi(x)dx=0$ such that
$$
B^{-1}(x+\alpha)S_E^v(x)B(x)=R_{\psi(x)}.
$$
\end{Theorem}

As should be clear from the above, the $PH2$ property of the dual
cocycles is key to our argument. Indeed, our Cantor spectrum result
for type I operators is largely  enabled by the fact that dual cocycles
of subcritical such operators are always partially hyperbolic with two dimensional
center, thus are $PH2.$ This is an immediate
corollary of the {\it recent duality approach to Avila's global theory}
\cite{gjyz}, which underlies this analysis. We thus have the
following reformulation of Theorem \ref{rot1}

\begin{Theorem}\label{rot2}
For a Schr\"odinger cocycle $(\alpha,S_E^v)$ associated with a type I operator,
there does not exist $B\in C^{\omega}(\T,SL(2,\R))$ and $\psi\in
C^\omega(\T,\R)$ with $\int_\T \psi(x)dx=0$ such that
$$
B^{-1}(x+\alpha)S_E^v(x)B(x)=R_{\psi(x)}.
$$
\end{Theorem}

\subsection{Further results and remarks}
Since the set of type I operators is open, Theorem \ref{main11}
implies certain robustness of the ten martini problem.

\begin{Definition}
 {\rm We say that $H_{v,\alpha,x}$ has {\it robust} property $A$ in a subspace $X$ if $H_{v+w,\alpha,x}$ has property $A$ for any sufficiently small $w\in X$.}
\end{Definition}

Theorem \ref{main11} implies that type I operators have
robust Cantor spectrum for all irrational $\alpha$ in the space of even trigonometric
polynomials. A stronger than Cantor spectrum property is that all gaps
prescribed by the gap labeling theory are open for all irrational
$\alpha$, known (originally for the almost
Mathieu operator) as the dry ten
martini problem \footnote{With the name coined, of course, by Barry
  Simon \cite{barry}.}. It has been recently solved for the non-critical
almost Mathieu
case \cite{dryAYZ}. It has not yet been proved for any non-almost-Mathieu
operator. However it is of course also a very natural question
whether the dry ten martini problem is robust.

If all gaps are open for a certain $\alpha$ we
will say that corresponding operator has {\it dry Cantor
spectrum}, (as previously established for the almost
Mathieu family for
various sets of $\alpha$
in \cite{cey,Puig,aj1,ly}), and robustness of that weaker statement  is also interesting. Other than for the
almost Mathieu, dry Cantor spectrum has been proved for
nonperturbatively small operators \eqref{sch} with Diophantine $\alpha$ in
  \cite{aj1} but only for
  ``typical'' $v$ (outside a set of infinite codimension) \footnote{A
    similar result in the 
perturbative (that is smallness depending on the Diophantine
constants) regime of small couplings has been earlier established by Puig \cite{puig} using 
Eliasson's reducibility theorem \cite{Eli92}.}. The only
  reason the set of infinite codimension had to be excluded in the
  argument of \cite{aj1} was the
  issue of simplicity of point spectrum for the dual
  operator.

  Let \begin{equation}\label{beta}
\beta(\alpha)=\limsup\limits_{n\rightarrow\infty}\frac{\ln q_{n+1}}{q_n},\end{equation}
where $q_n$ are denominators of the continued fraction
     approximants to $\alpha.$
     By the almost reducibility of subcritical cocycles \cite{arc2,ge}
     we get into the regime of Eliasson \cite{Eli92}, thus combining
with Theorem \ref{gjy-puig} we obtain

  \begin{Theorem}\label{main114}
Subcritical type I operators with trigonometric potentials $v$ and
with $\beta(\alpha)=0,$ have dry Cantor spectrum. It is
    robust in the space of 
    trigonometric potentials.
  \end{Theorem}
 \begin{Remark}   {\rm Strictly speaking, Eliasson's theorem in
     \cite{Eli92} was proved for a slightly stronger, classical polynomial
     Diophantine condition, although there are now many ways to see
     that it holds under the condition $\beta(\alpha)=0,$ as well. It
     can also be replaced by a
     combination of Theorems 4.1 in \cite{ge} and  \cite{aj1} with,
     again, Theorem 4.1 in \cite{aj1} improved to $\beta(\alpha)=0$ by
     the arguments of \cite{a1}.}
 \end{Remark}
  \begin{Remark} {\rm As before, we will 
    upgrade to all analytic $v$ as well as full analytic robustness in \cite{gjy2}.}
  \end{Remark}
For the particular case of the neighborhood of the
   almost Mathieu operator dry robustness can also be proved directly,
   without appealing to almost reducibility conjecture. Recall that
  $H^\delta_{\lambda,\alpha,x}$ is given by
\begin{equation}\label{hdelta}
 (H^\delta_{\lambda,\alpha,x}u)_n=u_{n+1}+u_{n-1}+(2\lambda\cos2\pi(x+n\alpha)+\delta f(x+n\alpha))u_n,\ \ n\in\Z,
\end{equation}
  We have
 \begin{Theorem}\label{main14}For
   $0<|\lambda|<1$ and 
   trigonometric polynomial $f$ there exists
   $\delta_0(\lambda,\|f\|_0)$ such that all gaps are open for
   $H^\delta_{\lambda,\alpha,x}$ provided
   $\beta(\alpha)=0,\;|\delta|<\delta_0$.
 \end{Theorem}
 \begin{Remark}
 {\rm Proof of Theorem \ref{main14} can also be viewed as a detailed proof
 of Theorem \ref{main114}, where one only needs to replace Corollary
 5.1 of \cite{LYZZ} by the almost reducibility \cite{arc2,ge}.}
 \end{Remark}


\begin{Remark}
{\rm It is also possible to extend Theorem \ref{main14}  to
  analytic or even $C^k$ perturbations and to any irrational $\alpha$.}
\end{Remark}
\begin{Remark}
{\rm Recently, it was proved by Ge-Wang-Xu \cite{gwx} that for
  Diophantine $\alpha$ and sufficiently large $\lambda$ operators
  $H^\delta_{\lambda,\alpha,x}$ have dry Cantor spectrum \footnote{The
    authors actually can prove such results for more general $C^2$
    cosine-like functions.}. Together with Theorem \ref{main14} it
  implies that the dry Cantor spectrum  is robust for almost Mathieu operators with small/large couplings.}
\end{Remark}
\begin{Remark}
{\rm By quantitative almost reducibility and quantitative Aubry
  duality, estimates on the spectral gap lengths are also possible in Theorem \ref{main14}.}
\end{Remark}

We conjecture the following more general result is true, making the
dry statement local\\

\noindent
\textbf{Conjecture 1:} {For any irrational $\alpha$  and $v\in
  C^\omega(\T,\R),$ each  type I energy $E_k\in
  \Sigma_{v, \alpha}$ satisfying $N_{v, \alpha}(E_k)= k\alpha (\mod
  \Z)$   is a boundary of an open gap.}
\begin{Remark}
{\rm Checking the proof, one will find that we actually did prove
  Conjecture 1 under the assumptions of Theorem \ref{main114}. The proof for the general case is in progress \cite{gjy2}.}
\end{Remark}
\begin{Remark}
{\rm Many of our bounds can be obtained in a quantitative way. This
  will be pursued in \cite{gjy2} where we will also use it to obtain quantitative lower bounds of the lengths of  spectral gaps.}
\end{Remark}

Of course, as far as the  conjectures go, the most natural one in this
context, after the present work, is that ten martini {\it always} holds for operators
\eqref{sch} with analytic $v.$ Indeed, we are not aware of any
counterexamples. A weaker conjecture would be that ten martini is
generic (in a variety of senses of varying strength) or, at least, dense. The
latter two would follow from another conjecture, interesting in its
own right: that type I operators are generic/dense.

\subsection{Structure of the rest of the paper.} 
In Section \ref{start} we discuss the ideas and strategy of the proof in more
detail. Section \ref{pre} contains the preliminaries, and in
Section \ref{tio}, we present basic properties and some typical
examples of type I operators.

In Section \ref{puig1} we involve the symplectic orthogonality
property of different eigenfunctions to  prove the simplicity of point
spectrum, and Section \ref{puig2}  is devoted to our quantitative
and all-frequency version of Puig's argument.

In Section \ref{kotanis1}, we establish Kotani theory and thus
$L^2$-reducibility for general minimal operators whose cocycles are
partially hyperbolic with two-dimensional center. As a consequence, we
prove that absence of Cantor spectrum implies improved
$C^\omega$-rotations reducibility for the associated quasiperiodic
finite-range cocycle in an interval.

Sections \ref{puig2} and, especially, \ref{kotanis1}
  represent the main hard analysis 
and contain the key contributions of this paper beyond the main
result. Both are of
independent interest.

Sections \ref{m12} and \ref{m14} are devoted to the proofs of Theorem
\ref{main12} and Theorem \ref{main14}, based on a combination of
generalized Kotani theory (developed in Section \ref{kotanis1}), and
all-frequency Puig's argument (developed in Sections \ref{puig1} and \ref{puig2}) 

\section{The strategy}\label{start}
The self-duality of the almost Mathieu family plays a key role in
the ten martini proof in several aspects, one of which is that the ($x$-independent)
spectra of $H_{\lambda,\alpha,x}$ and $H_{\frac 1{\lambda},\alpha,x},$
given by \eqref{amo}, coincide up to scaling by $\lambda.$
Therefore, it is sufficient to work in the (sub)critical regime
$|\lambda|\leq 1.$

For general type I operators \eqref{sch}, self-similarity is, of course, lost, so
we have to develop different arguments for sub and super critical
regions. However, Aubry duality (see Section \ref{aubry}) remains
a crucial tool, and the central object for us will be the dual
operator \eqref{fi}.  Operator \eqref{fi} has $d$ non-negative
Lyapunov exponents, which we denote as $\gamma_d(E)\geq \cdots\geq \gamma_1(E)\geq 0$.

The foundation of this work is quantitative global theory developed
recently in \cite{gjyz}, and in particular, the partial
hyperbolicity of dual cocycles established there. In particular, an
immediate corollary of \cite{gjyz} is

\begin{Theorem}\label{ths}
  $E$ is of type I for operator \eqref{sch} with trigonometric potential $v$ if
  and only if $\gamma_1(E)$ is simple.
\end{Theorem}
In other words, cocycles of the dual operator $L_{v,\alpha,\theta}$
are partially hyperbolic with a two dimensional center which, as we will
prove, is precisely the feature that allows to extend many of the
techniques developed for Schr\"odinger operators.

Theorem \ref{ths} will be repeated as Proposition \ref{simple le}
which will be given a more detailed but still a one-line proof. More specifically, as will be shown in Section \ref{secdual} the results of
\cite{gjyz} immediately imply that for operators \eqref{sch} of type I
we have
\begin{center}{\large
\begin{tabular}{c|c|c}
\hline\hline
Regime& {$H_{v,\alpha,x}$}& $L_{v,\alpha,\theta}$\\
\hline
subcritical&$L(E)=\omega(E)=0$& $L(E)=0$ and $\gamma_1(E)>0$ is simple  \\
\hline
critical &$L(E)=0,\omega(E)=1$&$L(E)=0$ and $\gamma_1(E)=0$ is simple \\
\hline
supercritical &$L(E)>0,\omega(E)=1$ & $L(E)>0$ and $\gamma_1(E)=0$ is simple\\
\hline
\end{tabular}}
\end{center}

Before we discuss the details of our strategy in each regime, let us
recall the proof of the ten martini problem in \cite{aj}.

As mentioned in
\cite{aj}, the two key breakthroughs, \cite{cey} and
\cite{puig}, already led to the proof of Cantor spectrum for
\eqref{amo} for an explicit set of a.e. $\alpha$, covering
correspondingly the Liouville and Diophantine regimes. The Diophantine
approach of
Puig was based on localization for completely resonant phases
\cite{j,jks}. He showed that it leads to dual reducibility to a
parabolic matrix, implying by a Moser-P\"oshel argument that the corresponding energy is a gap edge. Both
the required localization result and the parabolic reducibilty
implication are based on the almost Mathieu specifics, in particular
that both \eqref{amo} and its Aubry dual are second-difference
operators. It is conjectured however in \cite{aj,solving} that  localization
for completely resonant phases does not hold, and thus this approach
cannot work at all, for the most difficult arithmetic
mid-range of parameters (see \cite{liuresonant} for a recent
development).

Even more dramatic was the situation with the Liouville side
proof of Choi-Eliott-Yui \cite{cey}, that was used in its core and only brought to its
technical limits in \cite{aj}. The proof of
\cite{cey} is nothing
short of a $C^*$-algebraic miracle, utilizing remarkable properties of
Gauss polynomials, and is very specific for the irrational rotation
algebra (thus the almost Mathieu). Indeed, it has not yet been
extended even in a weak way to any other model, including for example, the extended
Harper's model, to which most other almost-Mathieu specific results
usually extend.

The complete proof in
\cite{aj}, that fully transpires in the (sub)critical regime, is based on four key ingredients,
\begin{enumerate}
\item Kotani theory for ergodic Schr\"odinger operators  and
  (fictitious) improved
  regularity of the $m$-function;
\item Puig's  almost Mathieu argument;
\item Fixed frequency localization for the dual model;
  \item The magic of \cite{cey}.
  \end{enumerate}
  The proof of \cite{aj} consists of the Diophantine side, developing
  (1),(2),(3) in a rather elaborate way in an argument by contradiction, and the
Liouville side, using (1) to bring (4) to its technical limits.

The fact that the two above approaches did meet in the middle has been
viewed as a miracle by the authors of \cite{aj}, with no
rational explanation. At the time, it has been unclear whether the
arithmetic dependence of the proof of \cite{aj} is something
intrinsically required. 

This has changed with the non-critical dry ten martini proof for the
almost Mathieu operator in
\cite{dryAYZ}. While the proof of
\cite{dryAYZ}, requiring delicate estimates, focuses only on the arithmetic
range not previously covered by \cite{aj1} and therefore does not fully
bypass the algebraic argument of \cite{cey},  the key
  idea works for all frequencies, as the authors found a way to run a Moser-P\"oschel
type argument based on quantitative {\it almost} reducibility to identity, rather
than reducibility to the identity for which there are Diophantine
obstructions. 

While inspired by \cite{dryAYZ}, we instead replace an argument through
 reducibility to the identity by the one through {\it reducibility to rotations with zero rotation number}. The latter holds for all irrational
$\alpha$ (\cite{afk,hy}).  We do it not with Moser-P\"oschel but
with Puig's duality approach itself (Theorems \ref{rot1}), replacing
the localized eigenfunctions in his argument by the  {\it almost
  localized} ones. Moreover, our argument works for
  higher-dimensional cocycles, thus allows to use duality while going beyond
  the almost Mathieu family. 
Remarkably, our  method to prove  Cantor
spectrum works for all irrational
frequencies and in a uniform way (for both sub and supercritical situations).

Of the ingredients (1)-(4) above, only (1) did not require the almost
Mathieu specifics, while (2)-(4) did, in a big way. Currently, we
don't have a good argument to extend (4) beyond the almost Mathieu
family. While the analogue of (3) for all type I operators is
forthcoming \cite{gj}, in view of the lack of (4), a proof
\'a la 
\cite{aj} would be missing some frequencies anyway. Thus, we don't use (3) or (4)
at all and develop instead a simpler unified argument, based entirely
on our extensions of (1) (Kotani-Simon for $PH2$ operators, Corollary
\ref{C0reducibility}), and (2) (the all-frequency Puig argument, Theorem \ref{rot1}).

More specifically, in the {\it critical or subcritical} regime, the main difficulty lies in the Puig's argument.
Indeed, the original Kotani theory for ergodic Schr\"odinger operators
still holds, so the absence of Cantor spectrum implies rotations reducibility.
However, Puig's original argument is no longer effective since the
dual operator is  long-rang, which does not directly imply
the simplicity of point spectrum.  Using the fact that, by \cite{gjyz}, the dual
operator is $PH2,$ Theorem \ref{tsimp}  can be invoked to obtain the simplicity of
dual eigenvalues and therefore the conclusion that type I cocycles cannot be
analytically reduced to the identity matrix, Theorem \ref{gjy-puig}.
To deal with all irrational frequencies in a uniform way, we replace
the localized eigenfunctions by the almost localized eigenfunctions in
the simplicity argument to obtain Theorem \ref{rot2}, which implies that
a type I cocycle cannot be analytically rotations reducible  in an
interval. This gives  the desired contradiction.

In the {\it supercritical regime}, the main difficulty lies in the lack of Kotani theory.
Indeed, one should not expect rotations reducibility anymore since the
Lyapunov exponent is positive. However, since the spectrum is
invariant under Aubry duality, we can instead prove Cantor spectrum
for the dual operator $L_{v,\alpha,\theta}$. However, the dual operator is now
long-range, so the original Kotani theory does not work, and neither
does the Kotani-Simon extension \cite{ks}. 
Here, we invoke our extension of Kotani theory for operators
with partially hyperbolic cocycles, Theorem \ref{L2 reducibility},
which, coupled again with the fact that, by \cite{gjyz}, the dual
of type I operators are $PH2,$ leads to the rotations reducibility result for dual finite-range cocycles, in the analytic category.
We then need to further develop our all-frequency Puig's argument,
making it work for dual finite-range operators, to prove that the dual
finite-range cocycle cannot be {\it partially} analytically rotations reducible in an
interval, and then completing the argument as in the subcritical case.

Other than the classical preliminaries and the basics of Avila's
global theory \cite{avila0,ajs} and its quantitative version \cite{gjyz}, our proof of Theorems
\ref{main11}, \ref{main12} is fully self-contained. The crucial for us
fact that dual cocycles of type I operators are $PH2$ is essentially
contained in \cite{gjyz}, but we also give a proof of this statement in the present paper, for completeness.

\section{Preliminaries}\label{pre}

\subsection{Continued fraction expansion} Let $\alpha\in (0,1)\backslash\Q$, $a_0:=0$ and $\alpha_0:=\alpha$. Inductively, for $k\geq 1$, we define
$$
a_k:=[\alpha_{k-1}^{-1}], \  \ \alpha_k=\alpha_{k-1}^{-1}-a_k.
$$
Let $p_0:=0$, $p_1:=1$, $q_0:=1$, $q_1:=a_1$. Again inductively, set
$p_k:=a_kp_{k-1}+p_{k-2}$, $q_k:=a_kq_{k-1}+q_{k-2}$. Then $q_n$ are
the  denominators of the best rational approximamts of $\alpha,$ since
we have $\|k\alpha\|_{\R/\Z}\geq \|q_{n-1}\alpha\|_{\R/\Z}$ for all
$k$ satisfying $\forall 1\leq k< q_n$. We also have 
$$
\frac{1}{2q_{n+1}}\leq \|q_n\alpha\|_{\R/\Z}\leq \frac{1}{q_{n+1}}.
$$

\subsection{Cocycles and the Lyapunov exponents}\label{3.2}
Let ${\rm M}(m,\C)$ be the set of all $m\times m$ matrices, $T:\Omega\rightarrow \Omega$ be a minimal homeomorphisim and $(\Omega,T,\mu)$ be ergodic. Given $A \in C^0(\Omega,{\rm M}(m,\C))$, we define the complex minimal cocycle $(T,A)$ by:
$$
(\alpha,A)\colon \left\{
\begin{array}{rcl}
\Omega \times \C^{m} &\to& \Omega \times \C^{m}\\[1mm]
(\omega,v) &\mapsto& (T\omega,A(\omega)\cdot v)
\end{array}
\right.  .
$$
The iterates of $(T,A)$ are of the form $(T,A)^n=(T^n,A_n)$, where
$$
A_n(\omega):=
\left\{\begin{array}{l l}
A(T^{n-1}\omega) \cdots A(T\omega) A(\omega),  & n\geq 0\\[1mm]
A^{-1}(T^n\omega) A^{-1}(T^{n+1}\omega) \cdots A^{-1}(T^{-1}\omega), & n <0
\end{array}\right.    .
$$
We denote by $L_1(T, A)\geq L_2(T,A)\geq...\geq L_m(T,A)$ the Lyapunov exponents of $(T,A)$ repeatedly according to their multiplicities, i.e.,
$$
L_k(T,A)=\lim\limits_{n\rightarrow\infty}\frac{1}{n}\int_{\Omega}\ln\sigma_k(A_n(\omega))d\mu,
$$
where  $\sigma_1(A_n)\geq...\geq \sigma_m(A_n)$ denote its singular values (eigenvalues of $\sqrt{A_n^*A_n}$). Since the k-th exterior product $\Lambda^kB$ of any $B\in M(m,\C)$ satisfies $\sigma_1(\Lambda^kB)=\|\Lambda^kB\|$, $L^k(T, A)=\sum\limits_{j=1}^kL_j(T,A)$ satisfies

$$
L^k(T,A)=\lim\limits_{n\rightarrow \infty}\frac{1}{n}\int_{\Omega}\ln\|\Lambda^kA_n(\omega)\|d\mu.
$$
\begin{Remark}
{\rm For $A\in C^0(\Omega,Sp(2d,\C)),$ where $Sp(2d,\C)$ is the set of complex symplectic matrices,  the Lyapunov exponents of $(T,A)$ come in pairs $\{\pm L_i(T,A)\}_{i=1}^d$.}
\end{Remark}

An important for us example is the minimal finite-range cocycle $(T,L_E^f)$ with
\begin{align}\tiny\label{1111}
L_{E}^{f}(\omega)=\frac{1}{a_d}
\begin{pmatrix}
-a_{d-1}&\cdots&-a_1&E-f(\omega)-a_0&-a_{-1}&\cdots&-a_{-d+1}&-a_{-d}\\
a_d& \\
& &  \\
& & & \\
\\
\\
& & &\ddots&\\
\\
\\
& & & & \\
& & & & & \\
& & & & & &a_{d}&
\end{pmatrix}.
\end{align}
Let
\begin{equation}\label{sdef}
C=\begin{pmatrix}
a_d&\cdots&a_1\\
0&\ddots&\vdots\\
0&0&a_d
\end{pmatrix},\ \ S=\begin{pmatrix}0&-C^*\\
C&0\end{pmatrix}.
\end{equation}
One can check that $L^f_E(\omega)$ is complex symplectic with respect
to $S,$ that is 
$$
(L^f_E(\omega))^*S L^f_E(\omega)=S,
$$ 
when $E\in\R$. Thus we can denote its non-negative Lyapunov exponents by $L_f^i(E)=L_i(T,L_{E}^{f})$ for $1\leq i\leq d$ for short.

\subsection{Uniform hyperbolicity and dominated splitting}\label{3.3}
For  $A\in C^0(\Omega,Sp(2d,\C))$,  we say the  cocycle $(T, A)$ is {\it uniformly hyperbolic} if for every $\omega \in \Omega$, there exists a continuous splitting $\C^{2}=E^s(\omega)\oplus E^u(\omega)$ such that for some constants $C>0,c>0$, and for every $n\geqslant 0$,
$$
\begin{aligned}
\lvert A_n(\omega)v\rvert \leqslant Ce^{-cn}\lvert v\rvert, \quad & v\in E^s(\omega),\\
\lvert A_n(\omega)^{-1}v\rvert \leqslant Ce^{-cn}\lvert v\rvert,  \quad & v\in E^u(T^n\omega).
\end{aligned}
$$
This splitting is invariant by the dynamics, which means that for every $\omega \in \Omega$,
$$
A(\omega)E^{\ast}(\omega)=E^{\ast}(T\omega),
$$
for $\ast=s,u$.   The set of uniformly hyperbolic cocycles is open in the $C^0$-topology.

For complex minimal cocycle $(T,A)\in C^0(\Omega,{\rm M}(m,\C))$, a
  related property is called {\it dominated splitting}. Recall that  Oseledets theorem provides us with  strictly decreasing
sequence of Lyapunov exponents $L_j(T,A) \in [-\infty,\infty)$ of multiplicity $m_j\in\N$, $1\leq j \leq \ell$ with $\sum_{j}m_j=m$, and for $\mu$ a.e. $\omega$, there exists
a measurable invariant decomposition $$\C^m=E_\omega^1\oplus E_\omega^2\oplus\cdots\oplus E_\omega^\ell$$ with $\dim E_\omega^j=m_j$ for $1\leq j\leq \ell$ such that $$
\lim\limits_{n\rightarrow\infty}\frac{1}{n}\ln\|A_n(\omega)v\|=L_j(T,A),\ \  \forall v\in E_\omega^j\backslash\{0\}.
$$
An invariant decomposition $\C^m=E_\omega^1\oplus
E_\omega^2\oplus\cdots\oplus E_\omega^\ell$  is  {\it dominated} if there exists $n$ such that  for any unit vector $v_j\in E_\omega^j\backslash \{0\}$, we have $$\|A_n(\omega)v_j\|>\|A_n(\omega)v_{j+1}\|.$$
Recall that Oseledets decomposition is a priori only measurable, however if an invariant decomposition $\C^m=E_\omega^1\oplus E_\omega^2\oplus\cdots\oplus E_\omega^\ell$  is  {\it dominated}, then $E_\omega^j$ depends  continuous on $\omega$ \cite{bdv}.

We also recall that  $(T,A)$     is called $k$-dominated (for some $1\leq k\leq m-1$) if there exists a dominated decomposition $\C^m=E^+_\omega \oplus E_\omega^- $    with $\dim E^+_\omega = k.$  It follows from the definitions that the Oseledets splitting is dominated if and only if $(T,A)$ is $k$-dominated for each $k$ such that  $L_k(T,A)> L_{k+1}(T,A)$.

\subsection{Global theory of one-frequency quasiperiodic cocycles} For
$\Omega=\T:=\R/\Z,$  $T:x\rightarrow x+\alpha,$ where
$\alpha\in\R\backslash\Q$, and $A:\T\to {\rm M}(m,\C),$ we call $(\alpha,A)$  a one-frequency quasiperiodic cocycle.
Global  theory of analytic one-frequency quasiperiodic cocycles was first developed for $SL(2,\C)$-cocycles \cite{avila0}, and
later  generalized to any ${\rm M}(m,\C)$-cocycles \cite{ajs}.  The key  concept for the global theory is the acceleration.  If $A\in C^{\omega}(\T,{\rm M}(m,\C))$ admits a holomorphic extension to $|\Im z|<\delta$, then for $|\e|<\delta$ we can define $A_\e\in C^{\omega}(\T,M(m,\C))$ by $A_\e(x)=A(x+i\e)$.
The accelerations of $(\alpha,A)$  are defined as
$$
\omega^k(\alpha,A)=\lim\limits_{\e\rightarrow 0^+}\frac{1}{2\pi\e}(L^k(\alpha,A_\e)-L^k(\alpha,A)), \qquad  \omega_k(\alpha,A)= \omega^k(\alpha,A)-\omega^{k-1}(\alpha,A).
$$
The key ingredient of the global theory is that the acceleration  is quantized.

\begin{Theorem}[\cite{avila0,ajs}]\label{ace}
There exists $1\leq l\leq m$, $l\in\N$, such that $l\omega^k$ and $l \omega_k$ are integers. In particular, if $A\in C^\omega(\T, SL(2,\C))$, then $\omega^1(\alpha,A)$ is an integer.
\end{Theorem}
\begin{Remark}\label{rem1}
{\rm If $L_j(\alpha,A)>L_{j+1}(\alpha,A)$, then $\omega^j(\alpha,A)$
  is an integer, as follows from the proof of Theorem 1.4 in \cite{ajs}, see also footnote 17 in \cite{ajs}.}
\end{Remark}

By subharmonicity, we know $L^k(\alpha,A(\cdot+i\e))$ is a convex function of $\e $  in a neighborhood of $0$,  unless it is identically equal to $-\infty$.
We say that $(\alpha,A)$ is {\it $k$-regular} if $\e\rightarrow L^k(\alpha,A(\cdot+i\e))$ is an affine function of $\e$ in a neighborhood of $0$.  In general, one can relate regularity and dominated splitting as follows.

\begin{Theorem}[\cite{avila0,ajs}]\label{t2.1}
Let $\alpha\in\R\backslash\Q$ and $A\in C^\omega(\T,M(m,\C))$. If $1\leq j\leq m-1$ is such that $L_j(\alpha,A)>L_{j+1}(\alpha,A)$, then $(\alpha,A)$ is $j$-regular if and only if $(\alpha,A)$ is $j$-dominated. In particular, if  $A\in C^\omega(\T, SL(2,\C))$ with $L(\alpha,A)>0$, then $(\alpha,A)$ is $1$-regular (or regular) if and only if $(\alpha,A)$ is uniformly hyperbolic.
\end{Theorem}

\subsection{One-frequency quasiperiodic $SL(2,\R)$-cocycles: rotation
  number and the IDS}\label{secrot}
For  one-frequency quasiperiodic $SL(2,\R)$-cocycle $(\alpha,A)$ with $A \in C^0(\T, {\rm SL}(2, \R))$,
assume that $A$ is homotopic to the identity. Then $(\alpha, A)$ induces the projective skew-product $F_A\colon \T \times \mathbb{S}^1 \to \T \times \mathbb{S}^1$
$$
F_A(x,w):=\left(x+\a,\, \frac{A(x) \cdot w}{|A(x) \cdot w|}\right),
$$
which is also homotopic to the identity. Lift $F_A$ to a map $\widetilde{F}_A\colon \T \times \R \to \T \times \R$ of the form $\widetilde{F}_A(x,y)=(x+\alpha,y+\psi_x(y))$, where for every $x \in \T$, $\psi_x$ is $\Z$-periodic.
Map $\psi\colon\T \times \R  \to \R$ is called a {\it lift} of $A$. Let $\mu$ be any probability measure on $\T \times \R$ which is invariant by $\widetilde{F}_A$, and whose projection on the first coordinate is given by Lebesgue measure.
The number
$$
\rho(\alpha,A):=\int_{\T \times \R} \psi_x(y)\ d\mu(x,y) \ {\rm mod} \ \Z
$$
 depends  neither on the lift $\psi$ nor on the measure $\mu$, and is called the \textit{fibered rotation number} of $(\alpha,A)$ (see \cite{H,johonson and moser} for more details).

Given $\theta\in\T$, let $
R_\theta:=
\begin{pmatrix}
\cos2 \pi\theta & -\sin2\pi\theta\\
\sin2\pi\theta & \cos2\pi\theta
\end{pmatrix}$.
If $A\colon \T\to{\rm PSL}(2,\R)$ is homotopic to $x \mapsto R_{nx/2}$ for some $n\in\Z$,
then we call $n$ the {\it degree} of $A$ and denote it by $\deg A$.
The fibered rotation number is invariant under real conjugacies which are homotopic to the identity. More generally, if $(\alpha,A_1)$ is conjugated to $(\alpha, A_2)$, i.e., $B(x+\alpha)^{-1}A_1(x)B(x)=A_2(x)$, for some $B \colon \T\to{\rm PSL}(2,\R)$ with $\deg{B}=n$, then
\begin{equation}\label{rotation number}
\rho(\alpha, A_1)= \rho(\alpha, A_2)+ \frac{n\alpha}{2}.
\end{equation}

In particular, for quasiperiodic Schr\"odinger cocycle
$(\alpha,S_E^v)$ where $S_E^v$ is given by \eqref{S},
we denote the rotation number $\rho(E):=\rho(\alpha,S_E^v)$.

The  {\it integrated density of states} (IDS)
$N_{v,\alpha}:\R\rightarrow [0,1]$ of $H_{v,\alpha,x}$ is defined by
$$
N_{v,\alpha}(E):=\int_{\T}\mu_{v,\alpha,x}(-\infty,E]dx,
$$
where $\mu_{v,\alpha,x}$ is the spectral measure of $H_{v,\alpha,x}$
and vector $\delta_0.$.

It is well known that $\rho(E)\in[0,\frac{1}{2}]$ is related to the integrated density of states $N=N_{v,\alpha}$ as follows:
\begin{equation}\label{relation}
N(E)=1-2\rho(E).
\end{equation}

\subsection{Aubry duality}\label{aubry}

Consider the fiber direct integral,
$$
\mathcal{H}:=\int_{\T}^{\bigoplus}\ell^2(\Z)dx,
$$
which, as usual, is defined as the space of $\ell^2(\Z)$-valued, $L^2$-functions over the measure space $(\T,dx)$.  The extensions of the
Sch\"odinger operators  and their long-range duals to  $\mathcal{H}$ are given in terms of their direct integrals, which we now define.
Let $\alpha\in\T$ be fixed. Interpreting $H_{v,\alpha,x}$ as fibers of the decomposable operator,
$$
H_{v,\alpha}:=\int_{\T}^{\bigoplus}H_{v,\alpha,x}dx,
$$
the family $\{H_{v,\alpha,x}\}_{x\in\T}$ naturally induces an operator on the space $\mathcal{H}$, i.e.,
$$
(H_{v,\alpha} \Psi)(x,n)= \Psi(x,n+1)+ \Psi(x,n-1) +  v(x+n\alpha) \Psi(x,n).
$$

Similarly,  the direct integral of long-range operators
$L_{v,\alpha,\theta}$ given by \eqref{fi},
denoted as $L_{v,\alpha}$, is given by
$$
(L_{v,\alpha}  \Psi)(\theta,n)=  \sum\limits_{k\in\Z} \hat{v}_k \Psi(\theta,n+k)+  2\cos2\pi (\theta+n\alpha) \Psi(\theta,n),
$$
where $\hat{v}_k$ is the $k$-th Fourier coefficient of $v(x)$.

Let  $U$ be the following operator on $\mathcal{H}:$
$$
(\mathcal{U}\phi)(\eta,m)=\sum_{n\in\Z}\int_{\T}e^{2\pi imx}e^{2\pi in(m\alpha+\eta)}\phi(x,n)dx.
$$
Then direct computations show that $U$ is unitary and satisfies
$$U H_{v,\alpha} U^{-1}=L_{v,\alpha}.$$ $U$ represents the so-called
{\it Aubry duality} transformation.
The quasiperiodic  long-range operator  $L_{v,\alpha,\theta}$ is called the dual operator of $H_{v,\alpha,x}$ \cite{gjls}.
More generally, let
$$
v(x)=\sum_{k=-d}^d\hat{v}_ke^{2\pi ikx},\ \ w(x)=\sum_{k=-l}^l \hat{w}_k e^{2\pi i kx}
$$
be two trigonometric polynomials.  We define a quasiperiodic
finite-range operator on $\ell^2(\Z)$ by
\begin{equation}\label{lvw}
(L^w_{v,\alpha,x}u)_n=\sum\limits_{k=-d}^d \hat{v}_ku_{n+k}+w(x+n\alpha)u_n, \ \ n\in\Z.
\end{equation}
So that we have $L_{v,\alpha,x}=L^{2\cos}_{v,\alpha,x}.$

As above, we can interpret $L^w_{v,\alpha,x}$ as fibers of the decomposable operator,
$$
L^w_{v,\alpha}:=\int_{\T}^{\bigoplus}L^w_{v,\alpha,x}dx,
$$
thus $L^w_{v,\alpha}$ acts on the space $\mathcal{H}$, by
$$
(L^w_{v,\alpha} \Psi)(x,n)= \sum\limits_{k=-d}^d \hat{v}_k\Psi(x,n+k) +  v(x+n\alpha) \Psi(x,n).
$$
We then have $$U L^w_{v,\alpha} U^{-1}=L^v_{w,\alpha},$$ so
finite-range operator  $L^v_{w,\alpha,\theta}$ is the dual of operator of $L^w_{v,\alpha,x}$ \cite{haro}.

\section{Type I operators}\label{tio}
In this section, we give the detailed definitions of  type I cocycles and operators and prove that all type I operators form an open set. We will give some natural examples of  type I operators. At the end, we give the dual characterization of type I operators based on quantitative global theory developed in \cite{gjyz}, which plays a crucial role in our further arguments.
\subsection{Type I cocycles}
Given $A \in C^\omega(\T,{\rm SL}(2,\R))$ one can extend it to the
band $\left\{z||\Im z|<h\right\}.$ For $\alpha\in\R\backslash\Q$ and
$\e<h,$ we can define the Lyapunov exponent of complexified cocycle
$(\alpha,A)$ just like we did for Schr\"dinger cocycles:
$$
L_\e(\alpha, A)=\lim\limits_{n\rightarrow\infty}\frac{1}{n}\int_{\T}\ln\|A(x+i\e+(n-1)\alpha)\cdots A(x+i\e+\alpha)A(x+i\e)\|dx.
$$

T-acceleration can also be defined for general analytic cocycles in the
same way as for Schr\"odinger cocycles, that is
$$
\bar{\omega}(\alpha,A):=\lim\limits_{\e\rightarrow \e_1^+}\frac{L_\e(\alpha,A)-L_{\e_1}(\alpha,A)}{\e-\e_1}
$$
where $\e_1$ is the first turning point of the piecewise affine
function $L_\e(\alpha,A),$ and $\bar{\omega}(\alpha,A)=1$ if
$L_\e(\alpha,A)\equiv 0, \e<h.$
When $\alpha$ is fixed through the argument, we will often write, for convenience,
$
L_\e(A):=L_\e(\alpha, A)$ and $\bar{\omega}(A):=\bar{\omega}(\alpha,A).$

An important fact is that  T-acceleration one cocycles are stable under analytic perturbations.
\begin{Lemma}\label{con}
For $(\alpha,A_0)\in \R\backslash\Q\times C^\omega_h(\T,\R)$ for some $h>0$ with $\bar{\omega}(\alpha,A_0)=1$, there is $\delta(\alpha,A_0)>0$ such that if $(\alpha_1,A)\in \R\backslash\Q\times C^\omega_h(\T,\R)$ with $\max\{\|A-A_0\|_h,|\alpha-\alpha_1|\}<\delta$, then $\bar{\omega}(\alpha_1,A)=1$.
\end{Lemma}
\begin{Remark}
{\rm This result is not true if one replaces  T-acceleration one by acceleration one.}
\end{Remark}
\begin{pf} Follows immediately by continuity and convexity of $L_\e.$ \end{pf}



\subsection{Type I operators}
Recall that analytic one-frequency quasiperiodic  Schr\"odinger
operator are given by \eqref{sch}, and
the corresponding Schr\"odinger cocycles are ($\alpha,S_E^v$)  where
$S_E^v(x)$ is given by \eqref{S}.
The ($x$-independent \cite{as}) spectrum of $H_{v,\alpha,x}$ is
denoted by $\Sigma_{v,\alpha}$. Recall that operators \eqref{sch} are
called Type I if $\bar{\omega}(E)=\bar{\omega}(S_E^v)=1$ for all $E\in \Sigma_{v,\alpha}$.

The property of being type I is stable,
more precisely
\begin{Corollary}\label{sgao}
Given $\alpha\in \R\backslash\Q$ and $v\in C_h^\omega(\T,\R)$ and
assume $\left\{H_{v,\alpha,x}\right\}_{x\in\T}$ is a type I operator,
then there is $\delta_0(v)>0$ such that if $w\in C^\omega_h(\T,\R)$ and $\alpha_1\in \R\backslash\Q$
are such that $\max\{\|w-v\|_h,|\alpha-\alpha_1|\}<\delta_0$  , then $\left\{H_{w,\alpha_1,x}\right\}_{x\in\T}$ is also a type I operator.
\end{Corollary}
which follows immediately from Lemma \ref{con}  by compactness.
The prime example of  type I operators is analytic perturbations of
the almost Mathieu operators (that were called PAMO in \cite{gjz}), see Example \ref{exa1},
$$
(H^\delta_{\lambda,\alpha,x}u)_n=u_{n+1}+u_{n-1}+\left(2\lambda\cos2\pi(x+n\alpha)+\delta f(x+n\alpha)\right)u_n, \ \ n\in\Z,
$$
where $f$ is a 1-periodic real analytic function.
\begin{Corollary}\label{deigen_pamo}
For $\alpha\in\R\backslash\Q$, $\lambda\neq 0$ and $f\in
C^\omega(\T,\R)$, there is $\delta_0(\lambda, \|f\|_h)>0$ such that for $|\delta|\leq \delta_0$, operator $\left\{H^\delta_{\lambda,\alpha,x}\right\}_{x\in\T}$ is of type I.
\end{Corollary}
\begin{pf}
By Corollary \ref{sgao}  we only need to prove the almost Mathieu
operator is of type I. This follows directly from a computation of
complexified Lyapunov exponent in Theorem 19 of \cite{avila0}: for any $E\in \R$, $\lambda\neq 0$ and $\e>0$, we have
$$
L(\alpha,S_E^{2\lambda\cos}(\cdot+i\e))=\max\{L(\alpha,S_E^{2\lambda\cos}),2\pi\e+\ln|\lambda|\}.
$$

\end{pf}
The proofs for Examples \ref{exa2}, \ref{exa3} follow from the fact  that
GPS model and supercritical generalized Harper's model are of  type I
\cite{wwxyz,gjyz} and a similar argument.
\subsection{The duality characterization}\label{secdual}Thi
Throughout this subsection, let
$$
v(x)=\sum\limits_{k=-d}^d \hat{v}_k e^{2\pi ikx}
$$
be a trigonometric polynomial of degree $d$. We will involve the Aubry
duality to study  type I operators \eqref{sch} with   trigonometric
polynomial potentials $v.$ The dual operator $L_{v,\alpha,\theta}$ is then defined by  the \eqref{fi} (see Section \ref{aubry} for details),
We denote the associated cocycle of the eigenequation $L_{v,\alpha,\theta}u=Eu$ by $(\alpha,L_{E,v})$. Its Lyapunov exponents are denoted by $\pm\gamma_{1}(E), \cdots, \pm\gamma_{d}(E)$. We assume that
$$
0\leq \gamma_{1}(E)\leq \gamma_{2}(E)\leq \cdots\leq \gamma_{d}(E).
$$

An important basis of our proof is the following duality characterization of  type I operators

\begin{Proposition}\label{simple le}
For $\alpha\in \R\backslash\Q$ and $E\in\R$, $\bar{\omega}(E)=1$ if and only if $\gamma_1(E)$ is simple.
\end{Proposition}
\begin{pf}
By Theorem 1 in \cite{gjyz},
$
\bar{\omega}(E)$ is equal to the multiplicity of $\gamma_1(E).$
Thus $\bar{\omega}(E)=1$ if and only if $\gamma_1(E)$ is simple.
\end{pf}
More importantly, we have
\begin{Theorem}\label{dominate1}
For $\alpha\in \R\backslash\Q$ and $E\in\R$, the cocycle $(\alpha, L_{E,v}^1)$ is $(d-1)$ and $(d+1)$-dominated.
\end{Theorem}
\begin{pf}
We let
$$
C=\begin{pmatrix}
\hat{v}_d&\cdots&\hat{v}_1\\
0&\ddots&\vdots\\
0&0&\hat{v}_d
\end{pmatrix},\ \ B(\theta)=\begin{pmatrix}
2\cos2\pi(\theta_{d-1})&\hat{v}_{-1}&\cdots&\hat{v}_{-d+1}\\
\hat{v}_1&\ddots&\ddots&\vdots\\
\vdots&\ddots&2\cos2\pi(\theta_1)&\hat{v}_{-1}\\
\hat{v}_{d-1}&\cdots&\hat{v}_1&2\cos2\pi(\theta)
\end{pmatrix}
$$
where $\theta_j=\theta+j\alpha$.  Then one can check that
\begin{equation}\label{strip2}
L_{E,v}(\theta+(d-1)\alpha)\cdots L_{E,v}(\theta)=:L_{d,E,v}(\theta)=\begin{pmatrix}
C^{-1}(EI-B(\theta))& -C^{-1}C^*\\
I_d&O_d
\end{pmatrix}
\end{equation}
where $I_d$ and $O_d$ are the $d$-dimensional identity and zero matrices, respectively.

Notice that \eqref{strip2} implies that we always have
$$
dL^{d-1}(\alpha,L_{E,v})=L^{d-1}(d\alpha,L_{d,E,v}).
$$
Thus by the definition of regularity, $(\alpha,L_{E,v})$ is
($d-1$)-regular if and only if  $(d\alpha,L_{d,E,v})$ is
($d-1$)-regular. Let $\left(\ell_{ij}\right)_{1\leq i,j\leq
  2d}:=(L_{d,E,v})_n(\theta)$. It is easy to check that each
$\ell_{ij}$ is a polynomial of $\cos2\pi(\theta)$ with degree $\leq
n$. Similarly,  let $L_{ij}$ be the $ij$-th entry of
$\Lambda^{d-1}(L_{d,E,v})_n(\theta).$ By the definition of wedge
product, each $L_{ij}$ is a polynomial of $\cos2\pi(\theta)$ of degree $\leq n(d-1)$. Hence
\begin{align*}
&|\omega^{d-1}(d\alpha,L_{d,E,v})|=\left|\lim\limits_{\e\rightarrow 0^+}\frac{1}{2\pi\e}(L^{d-1}(d\alpha,L_{d,E,v}(\cdot+i\e))-L^{d-1}(d\alpha,L_{d,E,v})\right|\\
=&\frac{1}{2\pi}\left|\lim\limits_{n\rightarrow \infty}\frac{1}{n}\int_{\T}\ln(\|\Lambda^{d-1}\left(L_{d,E,v}\right)_n(\theta+i\e)\|)d\theta-\lim\limits_{n\rightarrow \infty}\frac{1}{n}\int_{\T}\ln(\|\Lambda^{d-1}\left(L_{d,E,v}\right)_n(\theta)\|)d\theta\right|\\
\leq& d-1.
\end{align*}
It follows that
$$
|\omega^{d-1}(\alpha,L_{E,v}^{1})|=\left|\frac{\omega^{d-1}(d\alpha,L_{d,E,v}^{1})}{d}\right|\leq \frac{d-1}{d}<1.
$$
By Proposition \ref{simple le}, $\gamma_{1}(E)<\gamma_{2}(E)$, thus by
Remark  \ref{rem1}, $\omega^{d-1}(\alpha,L_{E,v}^{1})$ is an
integer. Thus, since $|\omega^{d-1}(\alpha,L_{E,v}^{1})|$ is strictly smaller than $1$, we have $\omega^{d-1}(\alpha,L_{E,v}^{1})=0$.  This implies that
$$
L^{d-1}(\alpha,L_{E,v}(\cdot+i\e))=L^{d-1}(\alpha,L_{E,v})
$$
for sufficiently small $\e>0$. A similar  argument  works for  $\e<0$.
This means  $(\alpha,L_{E,v})$ is ($d-1$)-regular, hence, by Theorem
\ref{t2.1}, $(\alpha,L_{E,v})$ is ($d-1$)-dominated, which implies
$(\alpha,L_{d,E,v})$ is ($d-1$)-dominated. Since $(\alpha,L_{d,E,v})$
is complex symplectic, we also have $(\alpha,L_{d,E,v})$ is ($d+1$)-dominated, thus $(\alpha,L_{E,v})$ is ($d+1$)-dominated.
\end{pf}
The next corollary follows directly from the definition of dominated splitting,
\begin{Corollary}\label{partial}
For $\alpha\in \R\backslash\Q$, and $E\in\R$ with $\bar{\omega}(E)=1$,
there exists a continuous invariant decomposition
$$
\C^{2d}=E^s(\theta)\oplus E^c(\theta)\oplus E^u(\theta).
$$
Moreover,  for any $\theta\in\T$, we have
\begin{equation}\label{eq10}
\limsup\limits_{n\rightarrow\infty}\frac{1}{n}\ln\|(L_{E,v})_n(\theta)v\|>0,\ \  \forall v\in E^s(\theta)\backslash\{0\},
\end{equation}
\begin{equation}\label{eq11}
\limsup\limits_{n\rightarrow\infty}\frac{1}{n}\ln\|(L_{E,v})_n(\theta)v\|<0,\ \  \forall v\in E^u(\theta)\backslash\{0\}.
\end{equation}
\begin{equation}\label{eq12}
{\rm dim} E^c(\theta)=2,
\end{equation}
\end{Corollary}
\begin{pf}
It follows from Theorem \ref{dominate1} and the definition of dominated splitting.
\end{pf}
We therefore have
\begin{Corollary} \label{Iph2}
If $E$ is of type I for $H_{v,\alpha,x}$ then $(\alpha, L_{E,v})$ is $PH2$.
\end{Corollary}

\section{Simplicity of point spectra of minimal $PH2$
    operators. Proof of Theorems \ref{tsimp} and \ref{gjy-puig}}\label{puig1}

Let finite-range operator $L^w_{v,\alpha,x}$ be given by
\eqref{lvw}. We denote  the cocycle induced by the eigenequation
$L^w_{v,\alpha,x}u=Eu$ by $(\alpha,L_{E,v}^w),$ so
\begin{align*}
L_{E,v}^{w}(x)=\frac{1}{\hat{v}_{d}}
\begin{pmatrix}
-\hat{v}_{d-1}&\cdots&-\hat{v}_{1}&E-w(x)-\hat{v}_0&-\hat{v}_{-1}&\cdots&-\hat{v}_{-d+1}&-\hat{v}_{-d}\\
\hat{v}_{d}& \\
& &  \\
& & & \\
\\
\\
& & &\ddots&\\
\\
\\
& & & & \\
& & & & & \\
& & & & & &\hat{v}_{d}&
\end{pmatrix}
\end{align*}
With
$$
S=\begin{pmatrix}
0&-C^*\\
C&0
\end{pmatrix},\ \  C=\begin{pmatrix}
\hat{v}_d&\cdots&\hat{v}_1\\
0&\ddots&\vdots\\
0&0&\hat{v}_d
\end{pmatrix}.
$$
we have, by the discussion in Section \ref{3.2}, that $L_{E,v}^w$  is complex symplectic with respect
to $S.$

We denote the non-negative Lyapunov exponents of $(\alpha,L_{E,v}^w)$ by $\{L^i(E)\}_{i=1}^d$. 
Let $\Sigma_{v,\alpha}^w$ be the
($x$-independent) spectrum of of $L_{v,\alpha,x}^w$.

With $PH2$ property as in Definition \ref{defph2}, Theorem \ref{gjy-puig} follows directly from the following slightly
more general version

\begin{Theorem}\label{contra1}
For $\alpha\in\R\backslash\Q$, if there exist $H\in
C^\omega(\T,Sp_{2l\times 2}(\R))$ such that
\begin{equation}\label{spe}
L^v_{E,w}(x)H(x)=H(x+\alpha),
\end{equation}
then the dual cocycle $(\alpha,L_{E,v}^w)$ is not $PH2.$
\end{Theorem}
Proof of Theorem \ref{contra1} will be split into the following three
subsections. But first we list two important corollaries.
\begin{Corollary}
Type I cocycle cannot be (analytically) reduced to the identity. I.e., if $H_{v,\alpha,x}$ is a type I operator, there does not exist $B\in C^{\omega}(\T,SL(2,\R))$ such that
$$
B^{-1}(x+\alpha)S_E^v(x)B(x)=Id.
$$
\end{Corollary}
\begin{pf}
By Theorem \ref{dominate1}, cocycles corresponding to the  duals of  type I operators are $PH2$.
\end{pf}

\begin{Corollary}
For $\alpha\in\R\backslash\Q$, there does not exist $F\in C^\omega(T,Sp_{2d\times 2}(\R))$ such that
\begin{equation}\label{spe}
L^{2\cos}_{E,v}(x)F(x)=F(x+\alpha).
\end{equation}
\end{Corollary}
\begin{pf}
The dual of $L_{v,\alpha,\theta}^{2\cos}$ is  Schr\"odinger operator
$H_{v,\alpha,x}$, and corresponding $SL(2,\R)$ Schr\"odinger cocycles are automatically $PH2.$
\end{pf}

\subsection{Symplectic orthogonality of the eigenpairs}
It is well known that for second-difference operators
$$
 (Hu)_n=u_{n-1}+u_{n+1}+V(n)u(n),\ \ n\in\Z.
 $$
 point spectrum is simple. Indeed,
 If $u,v\in\ell^2(\Z)$ satisfy $Hu=Eu$ and $Hv=Ev$, then, by the
 constancy of Wronskian,
 \begin{equation}\label{lg}
 u(n+1)v(n)-u(n)v(n+1)=0, \ \ \forall n\in\Z
 \end{equation}
For the   finite-range  operators
\begin{equation}\label{finite_difference}
(Lu)(n)=\sum\limits_{k=-d}^{d} a_k u_{n+k}+b(n)u_n, \ \ n\in\Z,
\end{equation}
where $a_{-k}=\overline{a}_k$ and $\{b(n)\}_{n\in\Z}\subset \R^{\Z}$
is a  bounded sequence of real numbers, and $d>1$ this of course no
longer works, as  Wronskians of pairs of eigenfunctions are no longer constant.

However, one can rewrite \eqref{lg} as
\begin{equation}\label{so}
\left\langle\begin{pmatrix} u(n+1)\\ u(n)\end{pmatrix}, \begin{pmatrix} 0&-1\\ 1&0\end{pmatrix}\begin{pmatrix} v(n+1)\\ v(n)\end{pmatrix}\right\rangle=0, \ \ \forall n\in\Z.
\end{equation}
and therefore view the simplicity of the point spectrum for
Schr\"odinger operator as a corollary of  {\it symplectic
  orthogonality} \eqref{so} of  eigenfunctions. It turns out
symplectic orthogonality still holds in the finite-range case.

Let $S,C$ be defined by \eqref{sdef}, and let
$$\vec{u}(n)=\begin{pmatrix} u(nd+d-1) \\ \vdots \\ u(nd)\end{pmatrix},\ \  \vec{v}(n)=\begin{pmatrix}v(nd+d-1) \\ \vdots \\ v(nd)\end{pmatrix}.
$$

\begin{Lemma}\label{symplectic}
For any  two eigenfunctions $u,v$ of $L$, corresponding to the same eigenvalue
$E,$ vectors $\begin{pmatrix}\vec{u}(1)\\ \vec{u}(0)\end{pmatrix}$ and $\begin{pmatrix}\vec{v}(1)\\ \vec{v}(0)\end{pmatrix}$ are symplectic orthogonal with respect to $S.$
\end{Lemma}


\begin{pf}
The eigenequation $Lu=Eu$ can be rewritten  as a second-order $2d$-dimensional difference equation by introducing the auxiliary variables
$$
\vec{u}(n)=\left(u(nd+d-1)\ \ \cdots\ \ u(nd+1)\ \ u(nd)\right)^T
$$
for $n\in \Z$. It is easy to check that $(\vec{u}(n))_{n\in\Z}$ satisfies
\begin{equation}\label{ef}
C\vec{u}(n+1)+T(n)\vec{u}(n)+C^*\vec{u}(n-1)=E\vec{u}(n),
\end{equation}
where $T(n)$ is the Hermitian matrix
$$
T(n)=\begin{pmatrix}
b(nd+d-1)&a_{-1}&\cdots&a_{-d+1}\\
a_1&\ddots&\ddots&\vdots\\
\vdots&\ddots&b(nd+1)&a_{-1}\\
a_{d-1}&\cdots&a_1&b(nd)
\end{pmatrix}.
$$
Note that equation \eqref{ef} is an eigenequation of the following vector-valued Schr\"odinger operator
$$
(L_d\vec{u})(n)=C\vec{u}(n+1)+T(n)\vec{u}(n)+C^*\vec{u}(n-1),
$$
acting on $\ell^2(\Z,\C^d)$.

To obtain a first-order system and the corresponding cocycle we use the fact that $C$ is invertible (since $a_d\neq 0$) and write
$$
\begin{pmatrix}
\vec{u}(n+1)\\
\vec{u}(n)
\end{pmatrix}
=\begin{pmatrix}
C^{-1}(EI-T(n))& -C^{-1}C^*\\
I_d&O_d
\end{pmatrix}
\begin{pmatrix}
\vec{u}(n)\\
\vec{u}(n-1)
\end{pmatrix},
$$
where $I_d$ and $O_d$ are the d-dimensional identity and zero matrices, respectively.
Set
$$
L_{d,E}(n)=\begin{pmatrix}
C^{-1}(EI_d-T(n))& -C^{-1}C^*\\
I_d&O_d
\end{pmatrix}
$$
Then for real $E$, matrix $L_{d,E}(n)$ is  complex-symplectic with respect to the complex-symplectic structure
$S$ given by \eqref{sdef}
i.e.,
\begin{equation}\label{s_eq}
(L_{d,E}(n))^{*}S L_{d,E}(n)=S, \  \ n\in\Z,
\end{equation}

Since $u,v\in \ell^2(\Z),$ 
we have
\begin{equation}\label{limit}
\lim\limits_{n\rightarrow \infty}\left\langle \begin{pmatrix}\vec{u}(n+1)\\ \vec{u}(n)\end{pmatrix}, \begin{pmatrix}0&-C^*\\
C&0\end{pmatrix}\begin{pmatrix}\vec{v}(n+1)\\ \vec{v}(n)\end{pmatrix}\right\rangle=0.
\end{equation}
On the other hand,
$$
\begin{pmatrix}\vec{u}(n+1)\\ \vec{u}(n)\end{pmatrix}=L_{d,E}(n)\cdots L_{d,E}(1)\begin{pmatrix}\vec{u}(1)\\ \vec{u}(0)\end{pmatrix}.
$$
By \eqref{s_eq}, for any $n\in\Z$,
\begin{align*}
&\left\langle \begin{pmatrix}\vec{u}(n+1)\\ \vec{u}(n)\end{pmatrix}, \begin{pmatrix}0&-C^*\\
C&0\end{pmatrix}\begin{pmatrix}\vec{v}(n+1)\\ \vec{v}(n)\end{pmatrix}\right\rangle\\
=&\left\langle L_{d,E}(n)\cdots L_{d,E}(1)\begin{pmatrix}\vec{u}(1)\\ \vec{u}(0)\end{pmatrix}, \begin{pmatrix}0&-C^*\\
C&0\end{pmatrix}L_{d,E}(n)\cdots L_{d,E}(1)\begin{pmatrix}\vec{v}(1)\\ \vec{v}(0)\end{pmatrix}\right\rangle\\
=&\left\langle \begin{pmatrix}\vec{u}(1)\\ \vec{u}(0)\end{pmatrix}, (L_{d,E}(n)\cdots L_{d,E}(1))^*\begin{pmatrix}0&-C^*\\
C&0\end{pmatrix}L_{d,E}(n)\cdots L_{d,E}(1)\begin{pmatrix}\vec{v}(1)\\ \vec{v}(0)\end{pmatrix}\right\rangle\\
=&\left\langle \begin{pmatrix}\vec{u}(1)\\ \vec{u}(0)\end{pmatrix}, \begin{pmatrix}0&-C^*\\
C&0\end{pmatrix}\begin{pmatrix}\vec{v}(1)\\ \vec{v}(0)\end{pmatrix}\right\rangle.
\end{align*}
By \eqref{limit}, we obtain
$$
\left\langle \begin{pmatrix}\vec{u}(1)\\ \vec{u}(0)\end{pmatrix}, \begin{pmatrix}0&-C^*\\
C&0\end{pmatrix}\begin{pmatrix}\vec{v}(1)\\ \vec{v}(0)\end{pmatrix}\right\rangle=0.
$$
\end{pf}
\subsection{Simplicity of point spectrum. Proof of Theorem \ref{tsimp}}
Let $\Sigma_f$ be the ($\omega$-independent spectrum of a minimal $PH2$
operator $L_{f,\omega}$ given by  \eqref{finite operator11}. Fix $E\in\Sigma_f$.
By the definitions of $PH2$ and dominated splitting, there exist continuous invariant decompositions
$$
\C^{2d}=E^s(\omega)\oplus E^c(\omega)\oplus E^u(\omega).
$$
and $C(E),\delta(E)>\delta'(E)>0$, such that for any $\omega\in\Omega$ and $n\geq 1$, we have
\begin{align}\label{eq10}
\left\|(L_{E}^{f})_{-n}(\omega)v\right\|>C^{-1}e^{\delta n},\ \  \forall v\in E^s(\omega)\backslash\{0\},\ \ \|v\|=1,
\end{align}
\begin{align}\label{eq11}
\left\|(L_{E}^{f})_{n}(\omega)u\right\|>C^{-1}e^{\delta n},\ \  \forall u\in E^u(\omega)\backslash\{0\},\ \ \|u\|=1.
\end{align}
\begin{align}\label{eq13}
\left\|(L_{E}^{f})_{\pm n}(\omega)u\right\|<Ce^{\delta' n},\ \  \forall w\in E^c(\omega)\backslash\{0\},\ \ \|w\|=1.
\end{align}
\begin{align}\label{eq12}
{\rm dim} E^c(\omega)=2.
\end{align}

Clearly, by \eqref{eq10}-\eqref{eq13} and invariance,  if $u(\omega)$ is an
$\ell^2$ eigenfunction, vector $\begin{pmatrix}\vec{u}(1,\omega)\\
  \vec{u}(0,\omega)\end{pmatrix}$ cannot have nonzero components in
either $E^s(\omega)$ or $E^u(\omega),$ for otherwise there would be
exponential growth at either $-\infty$ or $\infty.$

We now prove by contradiction. Assume $L_{f,\omega}$ has two linearly
independent eigenfunctions $u(\omega),v(\omega)$ corresponding to the
same eigenvalue $E.$ We then have
\begin{equation}\label{key}
\begin{pmatrix}\vec{u}(1,\omega)\\ \vec{u}(0,\omega)\end{pmatrix},\begin{pmatrix}\vec{v}(1,\omega)\\ \vec{v}(0,\omega)\end{pmatrix}\in E^c(\omega).
\end{equation}

On the other hand, by \eqref{eq12}, for any $\begin{pmatrix}\vec{x}(1,\omega)\\ \vec{x}(0,\omega)\end{pmatrix}\in E^c(\omega)$, we have
$$
\begin{pmatrix}\vec{x}(1,\omega)\\ \vec{x}(0,\omega)\end{pmatrix}=c_1(\omega)\begin{pmatrix}\vec{u}(1,\omega)\\ \vec{u}(0,\omega)\end{pmatrix}+c_2(\omega)\begin{pmatrix}\vec{v}(1,\omega)\\ \vec{v}(0,\omega)\end{pmatrix}.
$$
By Lemma \ref{symplectic}, we have
\begin{equation}\label{c22}
\left\langle \begin{pmatrix}\vec{u}(1,\omega)\\ \vec{u}(0,\omega)\end{pmatrix}, \begin{pmatrix}0&-C^*\\
C&0\end{pmatrix}\begin{pmatrix}\vec{v}(1,\omega)\\ \vec{v}(0,\omega)\end{pmatrix}\right\rangle=0.
\end{equation}
It follows that for any $\begin{pmatrix}\vec{x}(1,\omega)\\ \vec{x}(0,\omega)\end{pmatrix}\in E^c(\omega)$,
$$
\left\langle \begin{pmatrix}\vec{u}(1,\omega)\\ \vec{u}(0,\omega)\end{pmatrix}, \begin{pmatrix}0&-C^*\\
C&0\end{pmatrix}\begin{pmatrix}\vec{x}(1,\omega)\\ \vec{x}(0,\omega)\end{pmatrix}\right\rangle=0.
$$
This contradicts the non-degeneracy of the symplectic form.
\qed
\subsection{Proof of Theorem \ref{contra1}} We proceed by contradiction. Assume $(\alpha,L_{E,v}^w)$ is $PH2$ and  there exists $H\in C^\omega(\T,Sp_{2l\times 2}(\R))$ such that
\begin{equation}\label{form11}
L^v_{E,w}(x)H(x)=H(x+\alpha).
\end{equation}
Let
$$
H=\begin{pmatrix}
h_{1,1}&h_{1,2}\\
h_{2,1}&h_{2,2}\\
\vdots&\vdots\\
h_{2l,1}&h_{2l,2}
\end{pmatrix}\in C^\omega(\T,Sp_{2l\times 2}(\R)).
$$
By the definition of $L^v_{E,w}(x)$ and \eqref{form11}, one has for $i=1,2$,
\begin{equation}\label{e2e1}
-\frac{1}{\hat{w}_{l}}\left(\sum\limits_{k=1}^{2l} \hat{w}_{l-k}h_{k,i}(x)+(E-v(x))h_{l,i}(x)\right)-h_{1,i}(x+\alpha)=0,
\end{equation}
\begin{equation}\label{e3e1}
h_{k,i}(x)=h_{k+1,i}(x+\alpha), \ \ \forall 1\leq k\leq 2l-1.
\end{equation}
It follows from \eqref{e2e1} and \eqref{e3e1} that
\begin{equation}\label{e41}
\sum\limits_{k=-l}^l \hat{w}_{k}h_{l,i}(x+k\alpha)+(E-v(x))h_{l,i}(x)=0.
\end{equation}
Let $h_{l,i}(x)=\sum_k\hat{h}_i(k)e^{2\pi i kx}$ be the Fourier
expansion. Taking the Fourier transform of \eqref{e41}, we get
\begin{equation*}
\sum\limits_{k=-l}^l \hat{w}_{k}e^{2\pi ikn\alpha}\hat{h}_i(n)+\sum\limits_{k=-d}^d(E-\hat{v}_k)\hat{h}_i(n-k)=0.
\end{equation*}
Thus $\left\{\hat{h}_1(k)\right\}_{k\in\Z}$ and $\left\{\hat{h}_2(k)\right\}_{k\in\Z}$ are two linearly independent  eigenfunctions of $L^w_{v,\alpha,0}$ corresponding to the same eigenvalue $E$.

On the other hand, since $(\alpha,L_{E,v}^w)$ is $PH2$, by Theorem
\ref{tsimp}, $L^w_{v,\alpha,0}$ has simple point spectrum, a
contradiction. \qed

\section{An all-frequency Puig's argument. Proof of Theorem \ref{rot1}}\label{puig2}
Approximants of the Theorem \ref{rot1} can be equivalently reformulated as
\begin{Theorem}\label{contra2}
For any $\alpha\in\R\backslash\Q$, if there exist $H\in C^\omega(\T,Sp_{2\ell\times 2}(\R))$ and $\psi\in C^\omega(\T,\R)$ with $\int_\T \psi(x)dx=0$ such that
\begin{equation}\label{spe}
L^v_{E,w}(x)H(x)=H(x+\alpha)R_{\psi(x)},
\end{equation}
then the dual cocycle $(\alpha,L_{E,v}^w)$ is not $PH2$.
\end{Theorem}
\begin{Remark}
{\rm Roughly speaking, while Theorem \ref{contra1} can be used to prove Cantor
  spectrum for type I operators with  {\it Diophantine}
  frequencies, Theorem \ref{contra2} will be used to prove Cantor
  spectrum for type I operators with {\it all}  irrational frequencies.}
\end{Remark}
Proof of Theorem \ref{contra2} will be split into the following two
subsections. We first list two important corollaries.
\begin{Corollary}
If operator $H_{v,\alpha,x}$ is of type I, then there does not exist $B\in C^{\omega}(\T,SL(2,\R))$  and $\psi\in C^\omega(\T,\R)$ with $\int_\T \psi(x)dx=0$ such that
$$
B^{-1}(x+\alpha)S_E^v(x)B(x)=R_{\psi(x)}.
$$
\end{Corollary}

\begin{Corollary}
For any $\alpha\in\R\backslash\Q$, there does not exist $F\in C^\omega(\T,Sp_{2d\times 2}(\R))$  and $\psi\in C^\omega(\T,\R)$ with $\int_\T \psi(x)dx=0$ such that
\begin{equation}\label{spe}
L^{2\cos}_{E,v}(x)F(x)=F(x+\alpha)R_{\psi(x)}.
\end{equation}
\end{Corollary}

\subsection{Quantitative almost reducibility via rotations
  reducibility} In this subsection, we derive quantitative almost
reducibility from rotations reducibility for quasiperiodic
finite-range operators. Let $p_n/q_n$ be the approximants of the continued fraction expansion of $\alpha$ 
By the definition \eqref{beta} of $\beta(\alpha)$, for any $0<\e<\frac{\beta}{100}$, there is a subsequence $q_{n_k}$ of $q_n$ such that
\begin{equation}\label{naaa}
q_{n_{k}+1}>e^{(\beta-\e)q_{n_k}}.
\end{equation}
The key technical fact is
\begin{Theorem}\label{th1}
For all $\alpha\in\R\backslash\Q$, if there exist $F\in C_h^\omega(\T,Sp_{2\ell\times 2}(\R))$ and $\psi\in C_h^\omega(\T,\R)$ with $\int_\T \psi(x)dx=0$ such that
\begin{equation}\label{spe}
L^v_{E,w}(x)F(x)=F(x+\alpha)R_{\psi(x)},
\end{equation}
then we have
\begin{enumerate}
\item if $\beta(\alpha)<h$, then there exists $H\in C^\omega(\T,Sp_{2\ell\times 2}(\R))$ such that
$$
L^v_{E,w}(x)H(x)=H(x+\alpha),
$$
\item if $\beta(\alpha)\geq h$, then for every $k\ge 1$ there exist $H^k\in C_{h/2}^\omega(\T,Sp_{2\ell\times 2}(\R))$ and $\e^k\in C_{h/2}^\omega(\T,\R)$ such that
\begin{equation}\label{al}
L^v_{E,w}(x)H^k(x)=H^k(x+\alpha)R_{\e^k(x)},
\end{equation}
with
$$
|H^k|_{\frac{h}{2}}\leq |F|_he^{8(q_{n_k}+e^{-\frac{h}{2}q_{n_k}}q_{n_{k}+1})|\psi|_h},
$$
$$
|\e^k|_{\frac{h}{2}}\leq e^{-\frac{1}{20}q_{n_k+1}h}|\psi|_h.
$$
\end{enumerate}
\end{Theorem}
\begin{pf}
Case (1): If $\beta(\alpha)<h$, the argument is standard. Define
$$
\phi(x)=\sum\limits_{k\in\Z\backslash\{0\}}\frac{1}{e^{2\pi ik\alpha}-1}\hat{\psi}(k)e^{2\pi i kx}.
$$
It's easy to check that $\phi\in C^\omega(\T,\R)$ and
$$
\phi(x+\alpha)-\phi(x)=\psi(x).
$$
For
$$
H(x)=F(x)R_{\phi(x)},
$$
since $F\in C^\omega(\T,Sp_{2\ell\times 2}(\R))$, one can check that $H\in C^\omega(\T,Sp_{2\ell\times 2}(\R))$ and
$$
L^v_{E,w}(x)H(x)=H(x+\alpha)
$$
\qed

Case (2): If $\beta(\alpha)\geq h$, we need the following lemma.
\begin{Lemma}\label{small}
For $\alpha$ with $\beta(\alpha)\geq h$ and $f\in C_h^\omega(\T,\R)$ with $\int_\T f(x)dx=0$, there exist sequences of $g_k\in C_{h/2}^\omega(\T,\R)$ such that
$$
|g_k|_{\frac{h}{2}}\leq 8(q_{n_k}+e^{-\frac{h}{2}q_{n_k}}q_{n_{k}+1})|f|_h,
$$
$$
|f(x)-(g_k(x+\alpha)-g_k(x))|_{\frac{h}{2}}\leq e^{-\frac{1}{20}q_{n_k+1}h}|f|_h.
$$
\end{Lemma}
\begin{pf}
First we observe that
\begin{Proposition}\label{naaa1}
For any $0<|k|<\frac{q_{n+1}}{6}$, if $k\notin R_n=\{\ell q_n: \ell\in\Z\}$ and $q_{n+1}>100q_n$, then
$$
\|k\alpha\|_{\R/\Z}\geq \frac{1}{4q_n}.
$$
\end{Proposition}
\begin{pf}
Since $k\notin R_n=\{\ell q_n: \ell\in\Z\}$, we have that
$$
k=\ell_0 q_n+r, \ \ 0< r\leq q_n-1.
$$
On the other hand,
$$
|\ell_0|=\left|\frac{k-r}{q_n}\right|<\frac{q_{n+1}}{6q_n}+1.
$$
Thus
\begin{align*}
\|k\alpha\|_{\R/\Z}&\geq \|r\alpha\|_{\R/\Z}-|\ell_0|\|q_n\alpha\|_{\R/\Z}\\
&\geq \frac{1}{2q_n}-\frac{1}{6q_n}-\frac{1}{q_{n+1}}\geq \frac{1}{4q_n}.
\end{align*}
\end{pf}
Define $N_k=\left[\frac{q_{n_{k}+1}}{6}\right]$ and
$$
g_k(x)=\sum\limits_{j=-N_k}^{-1}\frac{\hat{f}(j)}{e^{2\pi i j\alpha}-1}e^{2\pi ijx}+\sum\limits_{j=1}^{N_k}\frac{\hat{f}(j)}{e^{2\pi i j\alpha}-1}e^{2\pi ijx}.
$$
In view of \eqref{naaa} and Proposition \ref{naaa1}, we distinguish two cases:

Case 1: $0< |j|< \frac{q_{n_k+1}}{6}$, $j\notin R_{n_k}$, then $|e^{2\pi i j\alpha}-1|\geq \frac{1}{4q_{n_k}}$,

Case 2: $q_{n_k}\leq j< \frac{q_{n_k+1}}{6}$, $j\in R_{n_k}$, then $|e^{2\pi i j\alpha}-1|\geq \frac{1}{2q_{n_k+1}}$.\\
It follows that
$$
|g_k|_{\frac{h}{2}}\leq 4q_{n_k}|f|_h+4e^{-\frac{h}{2}q_{n_k}}q_{n_{k}+1}|f|_h.
$$
Moreover, we have
$$
f(x)-(g_k(x+\alpha)-g_k(x))=\sum\limits_{|j|\geq q_{n_k+1}/6}\hat{f}(j)e^{2\pi ijx},
$$
which implies that
$$
|f(\cdot)-(g_k(\cdot+\alpha)-g_k(\cdot))|_{\frac{h}{2}}\leq e^{-\frac{1}{20}q_{n_k+1}h}|f|_h.
$$

\end{pf}

Thus, by Lemma \ref{small}, there are  $\phi^k\in C_{h/2}^\omega(\T,\R)$ such that
\begin{equation}\label{sd1}
|\phi^k|_{\frac{h}{2}}\leq 8(q_{n_k}+e^{-\frac{h}{2}q_{n_k}}q_{n_{k}+1})|\psi|_h,
\end{equation}
\begin{equation}\label{sd2}
|\psi-(\phi^k(\cdot+\alpha)-\phi^k(\cdot))|_{\frac{h}{2}}\leq e^{-\frac{1}{20}q_{n_k+1}h}|\psi|_h.
\end{equation}
Define
$$
H^k(x)=F(x)R_{\phi^k(x)}.
$$
We have
\begin{equation}\label{al}
L^v_{E,w}(x)H^k(x)=H^k(x+\alpha)R_{\e^k(x)},
\end{equation}
where $\e^k(x)=\psi(x)-(\phi^k(x+\alpha)-\phi^k(x))$ and $H^k\in C^\omega(\T,Sp_{2\ell\times 2}(\R))$. Moreover,
$$
|H^k|_{\frac{h}{2}}\leq |F|_he^{8(q_{n_k}+e^{-\frac{h}{2}q_{n_k}}q_{n_{k}+1})|\psi|_h},
$$
$$
|\e^k|_{\frac{h}{2}}\leq |\psi(\cdot)-(\phi^k(\cdot+\alpha)-\phi^k(\cdot))|_{\frac{h}{2}}\leq e^{-\frac{1}{20}q_{n_k+1}h}|\psi|_h.
$$

\end{pf}

\subsection{Proof of Theorem \ref{contra2}}
In case  $\beta(\alpha)<h$, Theorem \ref{contra2} follows
immediately from Theorem \ref{contra1}. In case  $\beta(\alpha)\geq
h$, we prove Theorem \ref{contra2} via quantitative almost
reducibility and quantitative Aubry duality. Essentially, we need to
establish a {\it quantitative version} of Puig's argument for
finite-range operators. The proof proceeds by contradiction. Given $\alpha\in\R\backslash\Q$, we assume
\begin{enumerate}
\item There exist $F\in C^\omega(\T,Sp_{2\ell\times 2}(\R))$ and $\psi\in C^\omega(\T,\R)$ with $\int_\T \psi(x)dx=0$ such that
\begin{equation}\label{spe}
L^v_{E,w}(x)F(x)=F(x+\alpha)R_{\psi(x)},
\end{equation}
\item $(\alpha,L_{E,v}^w)$ is $PH2.$
\end{enumerate}

By Theorem \ref{th1},  there exist  $H^k\in C^\omega(\T,Sp_{2\ell\times 2}(\R))$ and $\e^k\in C^\omega(\T,\R)$ such that
\begin{equation}\label{al}
L^v_{E,w}(x)H^k(x)=H^k(x+\alpha)R_{\e^k(x)},
\end{equation}
with
\begin{align}\label{ese1}
|H^k|_{\frac{h}{2}}\leq |F|_he^{8(q_{n_k}+e^{-\frac{h}{2}q_{n_k}}q_{n_{k}+1})|\psi|_h},
\end{align}
\begin{align}\label{ese2}
|\e^k|_{\frac{h}{2}}\leq e^{-\frac{1}{20}q_{n_k+1}h}|\psi|_h.
\end{align}
Let
\begin{equation}\label{hk}
H^k=\begin{pmatrix}
h_{1,1}&h_{1,2}\\
h_{2,1}&h_{2,2}\\
\vdots&\vdots\\
h_{2l,1}&h_{2l,2}
\end{pmatrix}\in C^\omega(\T,Sp_{2l\times 2}(\R)).
\end{equation}
Involving the form of $L^v_{E,w}(x)$ and \eqref{al}, one has for $j=1,2$,
\begin{equation}\label{e22}
-\frac{1}{\hat{w}_{l}}\left(\sum\limits_{k=1}^{2l} \hat{w}_{l-k}h_{k,j}(x)+(E-v(x))h_{d,j}(x)\right)-h_{1,j}(x+\alpha)=g_{1,j}(x),
\end{equation}
\begin{equation}\label{e33}
h_{m,j}(x)=h_{m+1,j}(x+\alpha)+g_{m+1,j}(x), \ \ \forall 1\leq m\leq 2l-1,
\end{equation}
where
$$
g_{m,1}(x)=(\cos2\pi(\e^k(x)-1)h_{m,1}(x+\alpha)+\sin2\pi(\e^k(x))h_{m,2}(x+\alpha),
$$
$$
g_{m,2}(x)=(\cos2\pi(\e^k(x)-1)h_{m,2}(x+\alpha)-\sin2\pi(\e^k(x))h_{m,1}(x+\alpha).
$$
It follows from \eqref{e22} and \eqref{e33} that
\begin{equation}\label{e466}
\sum\limits_{k=-l}^l \hat{w}_{k}h_{l,j}(x+k\alpha)+(E-v(x))h_{\ell,j}(x)=e_j(x),
\end{equation}
where $e_j$ is a linear combination of $\{g_{m,j}\}_{m=1}^{2l-1}$ of at most $4l^2$ terms. Hence by \eqref{ese1} and \eqref{ese2}, for $k$ sufficiently large depending on $E,v,w$, we have
\begin{equation}\label{ee1}
|e_j|_{\frac{h}{2}}\leq C(v,w)l^2|F|_he^{8(q_{n_k}+e^{-\frac{h}{2}q_{n_k}}q_{n_{k}+1})|\psi|_h}|\e^k|_{\frac{h}{2}}\leq e^{-\frac{1}{40}q_{n_k+1}h}.
\end{equation}

Let $h_{l,i}(x)=\sum_k\hat{h}_i(k)e^{2\pi i kx}$ be the Fourier
expansion. By \eqref{e466}, $\left\{\hat{h}_1(n)\right\}_{n\in\Z}$ and $\left\{\hat{h}_2(n)\right\}_{n\in\Z}$ are two approximate solutions of $L^w_{v,\alpha,0}u=Eu$, i.e. they satisfy
\begin{equation}\label{new4}
\left((L^w_{v,\alpha,0}-E)\hat{h}_1\right)(n)=\hat{e}_1(n),
\end{equation}
\begin{equation}\label{new4-1}
\left((L^w_{v,\alpha,0}-E)\hat{h}_2\right)(n)=\hat{e}_2(n),
\end{equation}
where $\{\hat{e}_j(n)\}_{n\in\Z}$ are the Fourier coefficients of $e_j$.

We denote
$$
I_{k}=\left[-e^{-\frac{h}{8}q_{n_k}}q_{n_k+1},e^{-\frac{h}{8}q_{n_k}}q_{n_{k}+1}\right].
$$
Then by \eqref{ese1}, for $n\notin I_k$, we have
\begin{align}\label{eeee1}
|\hat{h}_{1}(n)|\leq Ce^{e^{-\frac{h}{4}q_{n_k}}q_{n_k+1}} e^{-\frac{h}{2}e^{-\frac{h}{8}q_{n_k}}q_{n_k+1}}
\leq e^{-\frac{h}{2}e^{-\frac{h}{4}q_{n_k}}q_{n_k+1}}.
\end{align}
On the other hand, we have
\begin{Lemma}\label{initial}
There exists $n_0\in I_k$, such that
\begin{equation}\label{z1-estimate-21}
|\hat{h}_{1}(n_0)|\geq e^{\frac{h}{16}q_{n_k}}q^{-1}_{n_k+1}.
\end{equation}
\end{Lemma}
\begin{pf}
Denote
$$
\overrightarrow{u}=\begin{pmatrix}
h_{1,1}\\
h_{2,1}\\
\vdots\\
h_{2l,1}
\end{pmatrix},\ \ \overrightarrow{v}=\begin{pmatrix}
h_{1,2}\\
h_{2,2}\\
\vdots\\
h_{2l,2}
\end{pmatrix}
$$
By \eqref{hk} and the fact that $H^k\in C^\omega(\T,Sp_{2\ell\times 2}(\R))$, it follows that
$$
\overrightarrow{u}^*S\overrightarrow{v}=1.
$$
Thus
$$
\|\overrightarrow{u}\|_{L^2}\|S\overrightarrow{v}\|_{L^2}\geq 1
$$
which implies that
$$\|\overrightarrow{u}\|_{L^2}\geq \frac{1}{\|S\overrightarrow{v}\|_{L^2}}>\frac{1}{C\|H^k\|_{C^0}}.$$
By \eqref{e33}, one has
\begin{align}\label{b111}
2l\|\hat{h}_{1}\|_{\ell^2}\geq \|\overrightarrow{u}\|_{L^2}-4l^2\sup_{m,j}\|g_{m,j}\|_{L^2}\geq (C\|H^k\|_{C^0})^{-1} \geq (C\|F\|_{C^0})^{-1}\geq c.
\end{align}
By \eqref{eeee1}, we have that
\begin{align*}
\sum\limits_{n\notin I_k}|\hat{h}_{1}(n)|^2
\leq 2e^{-he^{-\frac{h}{4}q_{n_k}}q_{n_k+1}}.
\end{align*}
By \eqref{b111} and the fact that $|I_k|\leq 2e^{-\frac{h}{8}q_{n_k}}q_{n_k+1}$, it follows that there exists $n_0\in I_k$, such that
\begin{align*}
2e^{-\frac{h}{8}q_{n_k}}q_{n_k+1}|\hat{h}_{1}(n_0)|^2 \geq \sum\limits_{n\in I_k}|\hat{h}_{1}(n)|^2 &=\|\hat{h}_{1}\|_{\ell^2}^2-\sum\limits_{n\notin I_k}|\hat{h}_{1}(n)|^2\\
 &\geq  cl^{-1}-2e^{-he^{-\frac{h}{4}q_{n_k}}q_{n_k+1}} .
\end{align*}
Hence  there exists $K_0>0$, such that
$$
|\hat{h}_{1}(n_0)|\geq  e^{\frac{h}{16}q_{n_k}}q^{-1}_{n_k+1},
$$
provided $k>K_0$.
\end{pf}

Define
$$
\vec{h}_1(n_0)=\begin{pmatrix}
\hat{h}_{1}(n_0+d-1)\\
\hat{h}_{1}(n_0+d-2)\\
\vdots\\
\hat{h}_{1}(n_0-d)
\end{pmatrix},\ \ \vec{h}_2(n_0)=\begin{pmatrix}
\hat{h}_{2}(n_0+d-1)\\
\hat{h}_{2}(n_0+d-2)\\
\vdots\\
\hat{h}_{2}(n_0-d)
\end{pmatrix}.
$$
Notice that by \eqref{new4} and \eqref{new4-1}, for $j=1,2$, we have
\begin{align}\label{vech}
\begin{pmatrix}\hat{h}_j(n+d-1)\\ \hat{h}_j(n+d-2)\\ \vdots\\ \hat{h}_j(n-d)\end{pmatrix}&=(L^w_{E,v})_{n-n_0}(n_0\alpha)\vec{h}_j(n_0)+\sum\limits_{k=n_0+1}^n (L_{E,v}^w)_{n-k}((k+1)\alpha)p_j(k)
\end{align}
where
$$
p_j(k)=\begin{pmatrix}\hat{e}_j(k+d-1)\\0\\ \vdots\\ 0\end{pmatrix}.
$$
Let
\begin{align}\label{vecu}
\vec{u}_j(n)=\begin{pmatrix}u_j(n+d-1)\\ u_j(n+d-2)\\ \vdots\\ u_j(n-d)\end{pmatrix}=(L^w_{E,v})_{n-n_0}(n_0\alpha)\vec{h}_j(n_0),
\end{align}

We have

\begin{Lemma}\label{key1}
 There exists $K_1(E,v,w,\alpha)>0$ such that if  $k>K_1$, then
\begin{equation}\label{fwn}
\left|\vec{u}^*_2(n_0)S\vec{u}_1(n_0)\right|\geq e^{-e^{-q^{\frac{1}{2}}_{n_k}}q_{n_k+1}}.
\end{equation}
\end{Lemma}
\begin{pf}
Note that by \eqref{vech}, \eqref{vecu} and \eqref{ee1},  for $j=1,2$ and $|n|\leq e^{-q^{\frac{1}{4}}_{n_k}}q_{n_k+1}$, we have
\begin{align}\label{eeee22}
|u_j(n)-\hat{h}_j(n)|\leq q_{n_k+1}e^{Ce^{-q^{\frac{1}{4}}_{n_k}}q_{n_k+1}}e^{-\frac{1}{40}q_{n_k+1}h}\leq e^{-\frac{1}{80}q_{n_k+1}h}.
\end{align}
Combining \eqref{eeee1} and \eqref{eeee22}, we have for $e^{-\frac{h}{20}q_{n_k}}q_{n_k+1}\leq|n|\leq e^{-\frac{h}{100}q_{n_k}}q_{n_k+1}$,
\begin{equation}\label{esu}
|u_j(n)|\leq 2e^{-\frac{h}{2}e^{-\frac{h}{4}q_{n_k}}q_{n_k+1}}.
\end{equation}
Fix $E\in\Sigma_{v,\alpha}^w.$ Since $(\alpha,L_{E,v}^w)$ is $PH2$, there exist continuous invariant decompositions
$$
\C^{2d}=E^s(x)\oplus E^c(x)\oplus E^u(x).
$$
Moreover, there are $C(E),\delta(E)>\delta'(E)>0$, such that for any $x\in\T$ and $n\geq 1$, we have
\begin{align}\label{eq10a1}
\left\|(L_{E,v}^{w})_{-n}(x)v\right\|>C^{-1}e^{\delta n},\ \  \forall v\in E^s(x)\backslash\{0\},\ \ \|v\|=1,
\end{align}
\begin{align}\label{eq11a1}
\left\|(L_{E,v}^{w})_{n}(x)v\right\|>C^{-1}e^{\delta n},\ \  \forall v\in E^u(x)\backslash\{0\},\ \ \|u\|=1.
\end{align}
\begin{align}\label{eq13a1}
\left\|(L_{E,v}^{w})_{\pm n}(x)v\right\|<Ce^{\delta' n},\ \  \forall v\in E^c(x)\backslash\{0\},\ \ \|w\|=1.
\end{align}
\begin{align}\label{eq12a1}
{\rm dim} E^c(x)=2.
\end{align}
Thus, for $j=1,2$, there exist
$$
\vec{u}^s_j(n_0)\in E^s(T^{n_0}x), \ \ \vec{u}^c_j(n_0)\in E^c(T^{n_0}x),\ \ \vec{u}^u_j(n_0)\in E^u(T^{n_0}x)
$$
such that
\begin{align*}
\vec{u}_j(n_0)=\vec{u}^s_j(n_0)+\vec{u}^c_j(n_0)+\vec{u}^u_j(n_0).
\end{align*}
By \eqref{esu} and Lemma \ref{lee2}, we have
\begin{equation}\label{err}
\left\|\vec{u}^s_j(n_0)\right\|,\ \ \left\|\vec{u}^u_j(n_0)\right\|\leq Ce^{-(\delta-\delta') e^{-\frac{h}{100}q_{n_k}}q_{n_k+1}},\ \ j=1,2.
\end{equation}
\begin{Proposition}\label{ff5}
We have
$$
\vec{u}_2(n_0)=\frac{\langle\vec{u}_2(n_0),\vec{u}_1(n_0)\rangle}{\|\vec{u}_1(n_0)\|^2}\vec{u}_1(n_0)+\frac{\langle\vec{u}_2(n_0),S\vec{u}_1(n_0)\rangle}{\|S\vec{u}_1(n_0)\|^2}S\vec{u}_1(n_0)+\vec{e}_2(n_0)
$$
with $\|\vec{e}_2(n_0)\|\leq Ce^{- e^{-\frac{h}{200}q_{n_k}}q_{n_k+1}}$.
\end{Proposition}
\begin{pf}
Note that $\{\vec{u}_1^c(n_0),S\vec{u}^c_1(n_0)\}$ is an orthogonal
basis for $E^c(T^{n_0}x)$, since, by the $PH2$ property, ${\rm
  dim}E^c(T^{n_0}x)=2$. Thus we have
\begin{equation}\label{ff11}
\vec{u}^c_2(n_0)=\frac{\langle\vec{u}^c_2(n_0),\vec{u}^c_1(n_0)\rangle}{\|\vec{u}^c_1(n_0)\|^2}\vec{u}^c_1(n_0)+\frac{\langle\vec{u}^c_2(n_0),S\vec{u}^c_1(n_0)\rangle}{\|S\vec{u}^c_1(n_0)\|^2}S\vec{u}^c_1(n_0).
\end{equation}
By Lemma \ref{initial} and the fact that $\|H^k\|_{C^0}\leq \|F\|_{C^0}\leq C$, one has
\begin{equation}\label{ff4}
e^{\frac{h}{16}q_{n_k}}q^{-1}_{n_k+1}\leq |\vec{u}_{1}(n_0)| \leq C,\ \  e^{\frac{h}{16}q_{n_k}}q^{-1}_{n_k+1}\leq |S\vec{u}_{1}(n_0)|\leq C.
\end{equation}
Thus by \eqref{err} and the above inequality,
\begin{equation}\label{ff22}
\left\|\frac{\langle\vec{u}^c_2(n_0),\vec{u}^c_1(n_0)\rangle}{\|\vec{u}^c_1(n_0)\|^2}\vec{u}^c_1(n_0)-\frac{\langle\vec{u}_2(n_0),\vec{u}_1(n_0)\rangle}{\|\vec{u}_1(n_0)\|^2}\vec{u}_1(n_0)\right\|\leq Ce^{-e^{-\frac{h}{200}q_{n_k}}q_{n_k+1}},
\end{equation}
\begin{equation}\label{ff33}
\left\|\frac{\langle\vec{u}^c_2(n_0),S\vec{u}^c_1(n_0)\rangle}{\|S\vec{u}^c_1(n_0)\|^2}\vec{u}^c_1(n_0)-\frac{\langle\vec{u}_2(n_0),S\vec{u}_1(n_0)\rangle}{\|S\vec{u}_1(n_0)\|^2}S\vec{u}_1(n_0)\right\|\leq Ce^{-e^{-\frac{h}{200}q_{n_k}}q_{n_k+1}}.
\end{equation}
By \eqref{err}, \eqref{ff11}, \eqref{ff22} and \eqref{ff33}, we have
$$
\left\|\vec{u}_2(n_0)-\frac{\langle\vec{u}_2(n_0),\vec{u}_1(n_0)\rangle}{\|\vec{u}_1(n_0)\|^2}\vec{u}_1(n_0)+\frac{\langle\vec{u}_2(n_0),S\vec{u}_1(n_0)\rangle}{\|S\vec{u}_1(n_0)\|^2}S\vec{u}_1(n_0)\right\|\leq Ce^{-e^{-\frac{h}{200}q_{n_k}}q_{n_k+1}}.
$$
\end{pf}
We now prove \eqref{fwn} by contradiction. If
$$
\left|\vec{u}^*_2(n_0)S\vec{u}_1(n_0)\right|< e^{-e^{-q^{\frac{1}{2}}_{n_k}}q_{n_k+1}},
$$
then by \eqref{ff4} and Proposition \ref{ff5}, we have
\begin{equation}\label{error}
\left\|\vec{u}_2(n_0)-\frac{\langle\vec{u}_2(n_0),\vec{u}_1(n_0)\rangle}{\|\vec{u}_1(n_0)\|^2}\vec{u}_1(n_0)\right\|\leq e^{-\frac{1}{10}e^{-q^{\frac{1}{2}}_{n_k}}q_{n_k+1}},
\end{equation}
which means the orthogonal projection of $\vec{u}_2(n_0)$ to the vector  $S\vec{u}_1(n_0)$ is small.  By Lemma \ref{initial}, we have
\begin{equation}\label{gne5}
\left\|\frac{\langle\vec{u}_2(n_0),\vec{u}_1(n_0)\rangle}{\|\vec{u}_1(n_0)\|^2}\right\|\leq C\|F\|^2_{C^0}e^{hq_{n_k}}q^2_{n_k+1},
\end{equation}
provided $k$ is sufficiently large.

In the following, we consider
\begin{align*}
b(x)=\overrightarrow{u}^*(x)S \overrightarrow{v}(x)=\overrightarrow{u}^*(x)S\left(\overrightarrow{v}(x)-\frac{\langle\vec{u}_2(n_0),\vec{u}_1(n_0)\rangle}{\|\vec{u}_1(n_0)\|^2} \overrightarrow{u}(x)\right)\end{align*}
and aim to estimate $b(x)$.   First, we will show as a consequence of
\eqref{error}, that
$
\overrightarrow{v}(x)-\frac{\langle\vec{u}_2(n_0),\vec{u}_1(n_0)\rangle}{\|\vec{u}_1(n_0)\|^2} \overrightarrow{u}(x)
$
is small. To this end, we only need to estimate the
Fourier coefficients
$$
\hat{h}_1(n)-\frac{\langle\vec{u}_2(n_0),\vec{u}_1(n_0)\rangle}{\|\vec{u}_1(n_0)\|^2} \hat{h}_2(n).
$$
We distinguish two cases:\\

\textbf{Case I:} If $|n|\geq e^{-\frac{h}{100}q_{n_k}}q_{n_k+1}$, then
by \eqref{eeee1} and \eqref{gne5}, we have
\begin{align}\label{new10}
\left|\hat{h}_1(n)-\frac{\langle\vec{u}_2(n_0),\vec{u}_1(n_0)\rangle}{\|\vec{u}_1(n_0)\|^2} \hat{h}_2(n)\right|\leq e^{-\frac{h}{4}e^{-\frac{h}{2}q_{n_k}}q_{n_k+1}}
\end{align}
provided $k$ is sufficiently large.\\

\textbf{Case II:} If $|n|\leq e^{-\frac{h}{100}q_{n_k}}q_{n_k+1}$, we set
$$
\tilde{p}_n=\begin{pmatrix}\hat{e}_2(n+d-1)-\frac{\langle\vec{u}_2(n_0),\vec{u}_1(n_0)\rangle}{\|\vec{u}_1(n_0)\|^2} \hat{e}_1(n+d-1)\\0\\ \vdots\\ 0\end{pmatrix},
$$
$$
\tilde{y}_n=\begin{pmatrix}\hat{h}_2(n+d-1)-\frac{\langle\vec{u}_2(n_0),\vec{u}_1(n_0)\rangle}{\|\vec{u}_1(n_0)\|^2}  \hat{h}_1(n+d-1)\\ \hat{h}_2(n+d-2)-\frac{\langle\vec{u}_2(n_0),\vec{u}_1(n_0)\rangle}{\|\vec{u}_1(n_0)\|^2}  \hat{h}_1(n+d-2)\\
\vdots\\ \hat{h}_2(n-d)-\frac{\langle\vec{u}_2(n_0),\vec{u}_1(n_0)\rangle}{\|\vec{u}_1(n_0)\|^2}  \hat{h}_1(n-d)\end{pmatrix}.
$$
Then as a result of \eqref{new4} and \eqref{new4-1}, we have
$$\tilde{y}_n=\hat{L}^w_{E,v}(n\alpha)\tilde{y}_{n-1}+\tilde{p}_n,$$
which implies that
$$\tilde{y}_n=(L_{E,v}^w)_{n-n_0}(n_0\alpha)\tilde{y}_{n_0}+\sum\limits_{j=n_0+1}^n (L_{E,v}^w)_{n-j}(j\alpha)\tilde{p}_j$$
where $(n\alpha,(L_{E,v}^w)_n):= (\alpha, L_{E,v}^w)^n$, are the iterates of the dual cocycle.

To give an estimate of $\tilde{y}_n$,  first note that by assumption and  \eqref{ee1}, we have
\begin{eqnarray*}
 |\tilde{p}_n| & \leq&  e^{-\frac{1}{40}q_{n_k+1}h},\\
  |\tilde{y}_{n_0}|&\leq &\left\|\vec{u}_2(n_0)-\frac{\langle\vec{u}_2(n_0),\vec{u}_1(n_0)\rangle}{\|\vec{u}_1(n_0)\|^2}\vec{u}_1(n_0)\right\|\leq e^{-\frac{1}{10}e^{-q^{\frac{1}{2}}_{n_k}}q_{n_k+1}}.
\end{eqnarray*}
On the other hand, we have
\begin{equation}\label{ub}
\|(L_{E,v}^w)_n\|_{C^0}\leq C^n, \quad \forall n\in \Z.
\end{equation}
As a result, if $k$ is sufficiently large depending on $v,w$, then one can estimate
\begin{eqnarray}\label{new102}
 \nonumber |\tilde{y}_n|&\leq& C^{e^{-\frac{h}{100}q_{n_k}}q_{n_k+1}} e^{-\frac{1}{10}e^{-q^{\frac{1}{2}}_{n_k}}q_{n_k+1}}+2e^{-\frac{h}{100}q_{n_k}}q_{n_k+1}C^{e^{-\frac{h}{100}q_{n_k}}q_{n_k+1}} e^{-\frac{1}{40}q_{n_k+1}h}\\
&\leq&e^{-\frac{1}{100}e^{-q^{\frac{1}{2}}_{n_k}}q_{n_k+1}}.
\end{eqnarray}

Consequently,  by \eqref{new10} and \eqref{new102}, there is $K_1>0$ such that if $|k|\geq K_1$,
$$
\left\|h_{\ell,2}-\frac{\langle\vec{u}_2(n_0),\vec{u}_1(n_0)\rangle}{\|\vec{u}_1(n_0)\|^2} h_{\ell,1}\right\|_{C^0}\leq e^{-\frac{1}{100}e^{-q^{\frac{1}{2}}_{n_k}}q_{n_k+1}}.
$$
As a consequence, by \eqref{e33}
\begin{align*}
|b(x)|\leq 2l\left(\left\|h_{\ell,2}-\frac{\langle\vec{u}_2(n_0),\vec{u}_1(n_0)\rangle}{\|\vec{u}_1(n_0)\|^2} h_{\ell,1}\right\|_{C^0}+ 2e^{-\frac{1}{40}q_{n_k+1}h}\right)\leq  e^{-\frac{1}{200}e^{-q^{\frac{1}{2}}_{n_k}}q_{n_k+1}}.
\end{align*}
This contradicts 
\begin{align*}
|b(x)|= |\overrightarrow{u}^*(x)S\overrightarrow{v}(x)|=1.
\end{align*}
\end{pf}

Finally, by \eqref{eeee1} and \eqref{eeee22}, taking $n_1=[e^{-q^{\frac{1}{4}}_{n_k}}q_{n_k+1}]$, for $j=1,2$, we have
$$
for |u_j(n_1)|\leq e^{-e^{-\frac{1}{2}q^{\frac{1}{4}}_{n_k}}q_{n_k+1}}.
$$
It follows that
$$
\left|\vec{u}^*_2(n_1)S\vec{u}_1(n_1)\right|\leq Ce^{-e^{-\frac{1}{2}q^{\frac{1}{4}}_{n_k}}q_{n_k+1}}.
$$
Note that by the symplectic invariance,
$$
\vec{u}^*_2(n_1)S\vec{u}_1(n_1)=\vec{u}^*_2(n_0)S\vec{u}_1(n_0),
$$
which contradicts  \eqref{fwn}. Thus \eqref{spe} is not compatible
with $(\alpha,L_{E,v}^w)$ being $PH2$. \qed

\section{Kotani theory for {minimal $PH2$} operators. Proof
of Theorem \ref{L2 reducibility}}\label{kotanis1}
This section contains the technically most difficult part of this paper.
\subsection{$C^0$ reducibility.}
\cite{ak1,dcj,sim83,kot}.
Let $\mathbb{H}$ be the upper half plane.  Given an ergodic dynamical system
$(\Omega,\mu,T)$ and $f:\Omega\to\R,$
it is well known that there exists a continuous function
$m=m_{v,T}: \mathbb{H} \times \R/\Z
\to \mathbb{H}$ such that for the $S^f_{E}$ given by \eqref{S} we have $S^f_{E}(x) \cdot m(E,x)=m(E,Tx)$, thus
defining an invariant section for the Schr\"odinger cocycle
$(T,S_{E}^f)$ corresponding to  ergodic Schr\"odinger operator \eqref{ergs}.
\begin{equation} \label {minvariantsection}
(\alpha,S_{E}^f) (x,m(E,x))=(Tx,m(E,Tx)).
\end{equation}
Moreover, $E \mapsto m(E,x)$ is holomorphic on $\mathbb{H}$.

Let $L_f(E)$ be the  Lyapunov exponent of the Schr\"odinger cocycle $(T,S_E^f).$ 

The following result is an important consequence of the the classical Kotani theory 
\begin{Theorem}\label{kotani}
Assume $(\Omega,T)$ is minimal, $f:\Omega\rightarrow\R$ is continuous, and $L_f(E)=0$ in an open interval $J\subset\R$. Then for every $\omega\in\Omega$, the function $E\rightarrow m(E,\omega)$ admits a holomorphic extension to $\C\backslash(\R\backslash J)$, with values in $\mathbb{H}$. The function $m: \C\backslash(\R\backslash J)\rightarrow \mathbb{H}$ is continuous in both variables.
\end{Theorem}
Theorem \ref{kotani} played an important role in solving the almost
Mathieu ten martini problem \cite{aj}, since it implies the
$C^0$-reducibility of the Schr\"odinger cocycle $(T,S_E^f)$.  The aim
of this section is to present an analogue of such $C^0$-reducibility
result for  finite-range operators \eqref{finite operator11}
which allows existence of positive Lyapunov exponents. 

Let 
$\{L^i_f(E)\}_{i=1}^d$ be the
non-negative Lyapunov exponents of the complex symplectic cocycle
$(T,L_E^f)$ associated with operator $L_{f,\omega}$ given by \eqref{finite operator11}. The
matrix version of Theorem \ref{kotani}, i.e., Theorem \ref{kotani} for
operators \eqref{finite operator11} was proved by Kotani-Simon \cite{ks}
(see also Xu \cite{xu} for the monotonic case), assuming
$L_f^1(E)=\cdots=L_f^d(E)=0$.  Removing this restriction was stated as
a problem in \cite{ks}, and it has seen no serious progress until this
work. Here we
establish $C^0$-rotations reducibility in the $2$-dimensional center
for $PH2$ operators, thus solving the Kotani-Simon problem under the $PH2$
condition.


\begin{Theorem}\label{C0reducibility1}
Assume $(\Omega,T)$ is minimal, $(\Omega,\mu,T)$ is ergodic, $L_{f,w}$
is $PH2$, and $L^d_f(E)=0$ in an interval $I\subset \R$, then there exist  $(U_E,V_E)\in C^0(\Omega, Sp_{2d\times 2}(\R))$, and $R_E\in C^0(\Omega,SO(2,\R))$, depending analytically on $\mathbb{C}_{\delta'}\backslash(\R\backslash I)$ for some $\C_{\delta'}\subset\C_\delta$, such that
$$
L_E^f(\omega)(U_E(\omega),V_E(\omega))=(U_E(T\omega), V_E(T\omega))R_E(\omega).
$$
\end{Theorem}
Clearly Theorem \ref{C0reducibility1} directly implies Corollary
\ref{C0reducibility}. In order to ptove
Theorem \label{C0reducibility1} we start with some preparations.

\subsection{An extension of Johnson-Moser's theorem}\label{6.2}
For any $z\in\mathbb{H}$, the Green's function of Schr\"odinger
operator $H_{f,\omega}$ is defined as
$$
g_f(z,\omega)=\langle\delta_0,(H_{f,\omega}-z)^{-1}\delta_0\rangle.
$$
For Schr\"odinger operators, Johnson and Moser \cite{johonson and moser} (see also \cite{cs} for strip case) proved the following
relation between the Lyapunov exponent and the Green's function
$$
L'_f(z)=\int_\Omega g_f(z,\omega)d\mu,
$$

Johnson-Moser's theorem plays an important role in the proof of the classical Kotani theory.  In this subsection, we extend Johnson-Moser's theorem to minimal finite-range operators $L_{f,\omega}$, satisfying the $PH2$ condition. Notice that the eigenvalue equations $L_{f,\omega}u=Eu$ can be written as a second-order $2d$-dimensional difference equation by introducing the auxiliary variables
$$
\vec{u}_n=(u_{nd+d-1}\ \ \cdots\ \ u_{nd+1}\ \ u_{nd})^T \in \C^d
$$
for $n\in \Z$. By the proof of Lemma \ref{symplectic}  we have that $(\vec{u}_n)_n$ satisfies
\begin{equation}\label{ef100}
C\vec{u}_{n+1}+B(T^{dn}\omega)\vec{u}_n+C^*\vec{u}_{n-1}=E\vec{u}_n,
\end{equation}
where
$$
C=\begin{pmatrix}
a_d&\cdots&a_1\\
0&\ddots&\vdots\\
0&0&a_d
\end{pmatrix},
$$
and $B(\omega)$ is the Hermitian matrix
$$
B(\omega)=\begin{pmatrix}
f(T^{d-1}\omega)&a_{-1}&\cdots&a_{-d+1}\\
a_1&\ddots&\ddots&\vdots\\
\vdots&\ddots&f(T\omega)&a_{-1}\\
a_{d-1}&\cdots&a_1&f(\omega)
\end{pmatrix}.
$$
Moreover, equation \eqref{ef100} is an eigenequation of the following vector-valued Schr\"odinger operator
\begin{equation}\label{ldfo}
(L_{d,f,\omega}\vec{u})_n=C\vec{u}_{n+1}+B(T^{dn}\omega)\vec{u}_n+C^*\vec{u}_{n-1},
\end{equation}
acting on $\ell^2(\Z,\C^d)$.

Let $\Sigma_f$ be the spectrum of $L_{f,\omega}$. By  Definition
\ref{defph2} of the 
$PH2$ property and continuity of dominated splitting \cite{bdv}, there
is $\delta(f)>0$ such that if operator $L_{f,\omega}$ is $PH2$, then
every $E\in \mathbb{C}_\delta$ where $\mathbb{C}_\delta$ is a small
open neighborhood of $\Sigma_f,$ is $PH2$ for $L_{f,\omega}.$ 
It is known that for any $z\in
\mathbb{H}_\delta=\C_\delta\cap\mathbb{H}$, the cocycle $(T,L_z^f)$ is uniformly hyperbolic, thus $d$-dominated. Hence $(T,L_z^f)$ is $(d-1)$, $d$, $(d+1)$-dominated. As a consequence of dominated splitting, for any $z\in\mathbb{H}_\delta$, there exist continuous invariant decompositions
$$
\C^{2d}=E_z^s(\omega)\oplus E_z^+(\omega)\oplus E_z^-(\omega)\oplus E_z^u(\omega),\ \ \forall \omega\in\Omega,
$$
which implies that there are linearly independent $\{\vec{u}_z^i(\omega)\}_{i=1}^{d-1}\in E_z^s(\omega)$, $\vec{u}_z^+(\omega)\in E_z^+(\omega)$, $\vec{u}_z^-(\omega)\in E_z^-(\omega)$ and $\{\vec{v}_z^i(\omega)\}_{i=1}^{d-1}\in E_z^u(\omega)$ depending continuously on $\omega$ and analytically on $z$, such that
\begin{equation}\label{r1}
\begin{pmatrix}F_z^+(0,\omega)\\ F_z^+(-1,\omega)\end{pmatrix}:=\begin{pmatrix}\vec{u}_z^1(\omega),\cdots, \vec{u}_z^{d-1}(\omega),\vec{u}^+_z(\omega)\end{pmatrix}
\end{equation}
\begin{equation}\label{r2}
\begin{pmatrix}F_z^-(0,\omega)\\ F_z^-(-1,\omega)\end{pmatrix}:=\begin{pmatrix}\vec{u}^-_z(\omega), \vec{v}_z^{d-1}(\omega), \cdots, \vec{v}_z^1(\omega)\end{pmatrix}
\end{equation}
satisfy
$$
\sum\limits_{k=0}^\infty\|F_z^+(k,\omega)\|^2<\infty,\ \ \sum\limits_{k=0}^{-\infty}\|F_z^-(k,\omega)\|^2<\infty,
$$
where
$$
\begin{pmatrix}F_z^\pm(k,\omega)\\ F_z^\pm(k-1,\omega)\end{pmatrix}=(L_z^f)_{dk}(\omega)\begin{pmatrix}F_z^\pm(0,\omega)\\ F_z^\pm(-1,\omega)\end{pmatrix}.
$$
Moreover, for any $\omega\in\Omega$, we have
$$
\limsup\limits_{k\rightarrow \infty}\frac{1}{2k}\ln\left(\|\vec{u}_z^-(k,\omega)\|^2+\|\vec{u}_z^-(k+1,\omega)\|^2\right)= dL_f^d(z),
$$
\begin{align}\label{eeeq1}
\limsup\limits_{k\rightarrow \infty}\frac{1}{2k}\ln\left(\|\vec{u}_z^+(k,\omega)\|^2+\|\vec{u}_z^+(k+1,\omega)\|^2\right)=- dL_f^d(z).
\end{align}
where
\begin{equation}\label{u(n)}
\begin{pmatrix}u_z^\pm(k,\omega)\\ u_z^\pm(k-1,\omega)\end{pmatrix}=(L_z^f)_{dk}(\omega)u_z^\pm(\omega).
\end{equation}
and $L_f^d(z)$ is the smallest positive Lyapunov exponent of the
cocycle $(T, L_z^f).$
Once we have $F_z^{\pm}(k,\omega)$, one can define M matrices by
$$
M_+(z,\omega)=F_z^+(1,\omega)(F_z^{+}(0,\omega)^{-1},
$$
$$
M_-(z,\omega)=F_z^-(-1,\omega)(F_z^-(0,\omega))^{-1}.
$$
just as in  \cite{ks}, and note that $M_\pm$  satisfy  the following Ricatti equations.
\begin{Lemma}For any $z\in \mathbb{H}_\delta$, we have
\begin{equation}\label{rica}
CM_+(z,\omega)+C^*M^{-1}_+(z,T^{-d}\omega)+(B(\omega)-z)=0.
\end{equation}
\begin{equation}\label{rica1}
C^*M_-(z,\omega)+CM^{-1}_-(z,T^{d}\omega)+(B(\omega)-z)=0.
\end{equation}
\end{Lemma}
\begin{pf}
Note that
$$
CF_z^\pm(1,\omega)+C^*F_z^\pm(-1,\omega)+(B(\omega)-z)F_z^\pm(0,\omega)=0.
$$
The results follow from the definition of $M_\pm$.
\end{pf}

Similar to \cite{ks}, one can define the Green's matrix by
$$G(z,\omega)= \langle \vec{\delta}_0,(L_{d,f,\omega}-z)^{-1}\vec{\delta}_0\rangle,
$$
where
$$
\vec{\delta}_j(n)=
\begin{cases}
0& n\neq j\\
I_d &n=j
\end{cases}.
$$
The Green's matrix $G(z,\omega)$ can be then expressed as:

\begin{Lemma}\label{Green_Matrix}
For any $z\in\mathbb{H}_\delta$, we have
\begin{align*}
G(z,\omega)=(CM_+(z,\omega)+C^*M_-(z,\omega)+B(\omega)-z)^{-1}
\end{align*}
\end{Lemma}
\begin{pf}
It is easy to check that
\begin{align*}
&\langle \vec{\delta}_m,(H_{d,f,\omega}-z)^{-1}\vec{\delta}_n\rangle \\
=&\begin{cases}
F_z^+(m,\omega)(CF_z^+(n+1,\omega)+C^*F_z^-(n-1,\omega)+(B(\omega)-z)F_z^+(n,\omega))^{-1}& m\geq n\\
F_z^-(m,\omega)(CF_z^+(n+1,\omega)+C^*F_z^-(n-1,\omega)+(B(\omega)-z)F_z^+(n,\omega))^{-1}& m< n
\end{cases}.
\end{align*}
%
\end{pf}

The following proposition gives the relation between $M_\pm$ and the Green's matrix.
\begin{Proposition}\label{prop123} For any $z\in\mathbb{H}_\delta$, the following relations hold:
\begin{eqnarray}
\nonumber G(z,\omega)&=&(-C^*M^{-1}_+(z,T^{-d}\omega)+C^*M_-(z,\omega))^{-1}, \\
\nonumber G(z,T^{-d}\omega)&=& (CM_+(z,T^{-d}\omega)-CM^{-1}_-(z,\omega))^{-1},\\
\label{g3} G(z,\omega)C^*M^{-1}_+(z,T^{-d}\omega) &=& M_+(z,T^{-d}\omega)G(z,T^{-d}\omega) C - I_d.
\end{eqnarray}
\end{Proposition}
\begin{pf}
By Lemma \ref{Green_Matrix} and \eqref{rica}, one has
\begin{align*}
G(z,T^{-d}\omega)=(CM_+(z,T^{-d}\omega)-CM^{-1}_-(z,\omega))^{-1},
\end{align*}
\begin{align*}
G(z,\omega)= (-C^*M_+^{-1}(z,T^{-d}\omega)+C^*M_-(z,T\omega))^{-1}.
\end{align*}
Consequently, we have the following
\begin{align*}
G(z,\omega)C^*M^{-1}_+(z,T^{-d}\omega)
=&(-I_d+M_+(z,T^{-d}\omega)M_-(z,\omega))^{-1}\\ \nonumber
=&M_-^{-1}(z,\omega)(-M_-^{-1}(z,\omega)+M_+(z.T^{-d}\omega))^{-1}\\ \nonumber
=&M_+(z,T^{-d}\omega)(-M_-^{-1}(z,\omega)+M_+(z.T^{-d}\omega))^{-1}- I_d\\ \nonumber
=&M_+(z,T^{-d}\omega)G(z,T^{-d}\omega) C - I_d.
\end{align*}
\end{pf}

We can now formulate the $PH2$ extension of the Johnson-Moser's theorem.
\begin{Theorem}\label{mg}
For any $z\in \mathbb{H}_\delta$, we have 
$$ \frac{\partial L_f^d}{\partial \Im z}(z)=-\frac{1}{d}\Im\int_{\Omega}g(z,\omega) d\omega.$$
where
$$
g(z,\omega):=\langle \delta_d, (F_z^+(0,\omega))^{-1}G(z,\omega)F_z^+(0,\omega)\delta_d\rangle.
$$
\end{Theorem}
\begin{pf}
Let $(T^d,L_{d,z}^f)$ be the cocycle corresponding to the
eigenequation $L_{d,f,\omega}u=zu$ where $L_{d,f,\omega}$ is defined
by \eqref{ldfo}. By invariance, there is $\tau(z,\omega)$ depending continuously on $\omega$ and analytically on $z$ such that
$$
L_{d,z}^f(\omega)\begin{pmatrix}\vec{u}^+_z(0,\omega)\\ \vec{u}^+_z(-1,\omega)
\end{pmatrix}=\begin{pmatrix}\vec{u}^+_z(0,T^d\omega)\\ \vec{u}^+_z(-1,T^d\omega)
\end{pmatrix}\frac{1}{\tau(z,\omega)}.
$$
By \eqref{eeeq1}, we have
\begin{align*}
 d\cdot L_f^d(z)=\int_\Omega \ln \tau(z,\omega)d\omega.
\end{align*}
It suffices for us to prove
\be\label{gi}
\frac{\partial \tau(z,\omega)}{\partial z}\frac{1}{\tau(z,\omega)}=d\frac{\partial L_f^d}{\partial z}(z)=\int_{\Omega}g(z,\omega)d\omega,
\ee
Once we have this, then the result follows from  the Cauchy-Riemann equations.

Again by invariance and the definition of $\{F_z^+(k,\omega)\}_{k\in\Z}$, we have
\begin{align}\label{tm}
\tau(z,\omega)&=\langle \delta_d, (F_z^+(1,\omega))^{-1}F_z^+(0,T^d\omega)\delta_d\rangle.
\end{align}
\begin{Lemma}\label{wewe1}
We have that
\begin{align*}
\frac{\partial \tau(z,\omega)}{\partial z}\frac{1}{\tau(z,\omega)}=&\langle\delta_d,(F_z^+(0,\omega))^{-1}\frac{\partial  M^{-1}_+(z,\omega)}{\partial z}  M_+(z,\omega)F_z^+(0,\omega)\delta_d\rangle\\
&-h(z,\omega)+h(z,T^d\omega).
\end{align*}
where
$h(z,\omega)=\langle \delta_d, (F_z^+(0,\omega))^{-1}\frac{\partial  F_z^+(0,\omega)}{\partial z}\delta_d\rangle$.
\end{Lemma}
\begin{pf}
By invariance, we have  for some $U$
\begin{align*}
(F_z^+(1,\omega))^{-1}F_z^+(0,T^d\omega)=:\widetilde{M}^{-1}_+(z,\omega)={\rm diag}\{U(z,\omega),\tau(z,\omega)\}
\end{align*}
It follows,
\begin{align*}
& (F_z^+(0,\omega))^{-1}\frac{\partial  M^{-1}_+(z,\omega)}{\partial z}  M_+(z,\omega)F_z^+(0,\omega)\\
=&(F_z^+(0,\omega))^{-1}\frac{\partial  F_z^+(0,\omega) \widetilde{M}^{-1}_+(z,\omega)(F_z^+(0,T^d\omega))^{-1}}{\partial z}  M_+(z,\omega)F_z^+(0,\omega)\\
=&(F_z^+(0,\omega))^{-1}\frac{\partial  F_z^+(0,\omega)}{\partial z} \widetilde{M}^{-1}_+(z,\omega)(F_z^+(0,T^d\omega))^{-1} M_+(z,\omega)F_z^+(0,\omega)\\
&+(F_z^+(0,\omega))^{-1} F_z^+(0,\omega) \frac{\partial\widetilde{M}^{-1}_+(z,\omega)}{\partial z}(F_z^+(0,T^d\omega))^{-1}  M_+(z,\omega)F_z^+(0,\omega)\\
&+(F_z^+(0,\omega))^{-1} F_z^+(0,\omega) \widetilde{M}^{-1}_+(z,\omega)\frac{\partial(F_z^+(0,T^d\omega))^{-1}}{\partial z} M_+(z,\omega)F_z^+(0,\omega)\\
=& \frac{\partial \widetilde M_+^{-1} (z,\omega)}{\partial z}\widetilde M_+(z,\omega) + E(z,\omega)
\end{align*}
where  we set
\begin{align}\label{e2}
E(z,\omega):=&(F_z^+(0,\omega))^{-1}\frac{\partial  F_z^+(0,\omega)}{\partial z}\\ \nonumber
&-\widetilde{M}^{-1}_+(z,\omega)(F_z^+(0,T^d\omega))^{-1}\frac{\partial F_z^+(0,T^d\omega)}{\partial z}\widetilde{M}_+(z,\omega).
\end{align}
Here we used that for any invertible matrix $A$, we have
$\frac{\partial A^{-1}}{\partial z}A=-A^{-1}\frac{\partial A}{\partial
  z}.$

The result follows since $\widetilde{M}_+(z,\omega)$ is a block diagonal matrix.
\end{pf}
Finally, we introduce the  auxiliary function
 $$f(z,\omega)=\langle \delta_d, (F_z^+(0,\omega))^{-1}G(z,\omega)\frac{\partial CM_+(z,\omega)}{ \partial z}F_z^+(0,\omega)\delta_d\rangle.
$$
\begin{Lemma}
 We have that
\begin{align}\label{ft}
&\frac{\partial \tau(z,T^{-d}\omega)}{\partial z}\frac{1}{\tau(z,T^{-d}\omega)} -g(z,\omega)\\ \nonumber
=&-f(z,\omega)+f(z,T^{-d}\omega)+h(z,\omega)-h(z,T^{-d}\omega).
\end{align}
\end{Lemma}
\begin{pf}
By \eqref{rica}, we have
\begin{align*}
&\frac{\partial CM_+(z,\omega)}{ \partial z}=-\frac{\partial C^*M^{-1}_+(z,T^{-d}\omega)}{ \partial z}+I_d\\
=&C^*M^{-1}_+(z,T^{-d}\omega)\frac{\partial M_+(z,T^{-d}\omega)}{\partial z}M^{-1}_+(z,T^{-d}\omega)+I_d.
\end{align*}
Then by Proposition \ref{prop123},
\begin{align*}
&(F_z^+(0,\omega))^{-1}G(z,\omega)\frac{\partial CM_+(z,\omega)}{ \partial z}F_z^+(0,\omega)\\
=&(F_z^+(0,\omega))^{-1}G(z,\omega)\left(C^*M^{-1}_+(z,T^{-d}\omega)\frac{\partial M_+(z,T^{-d}\omega)}{\partial z}M^{-1}_+(z,T^{-d}\omega)+I_d\right)F_z^+(0,\omega)\\
=&(F_z^+(0,\omega))^{-1}G(z,\omega)C^*M^{-1}_+(z,T^{-d}\omega)\frac{\partial M_+(z,T^{-d}\omega)}{\partial z}M^{-1}_+(z,T^{-d}\omega)F_z^+(0,\omega)\\
&+(F_z^+(0,\omega))^{-1}G(z,\omega)F_z^+(0,\omega) \\
=&(F_z^+(0,\omega))^{-1}\left(M_+(z,T^{-d}\omega)G(z,T^{-d}\omega) C - I_d\right)\frac{\partial M_+(z,T^{-d}\omega)}{\partial z}M^{-1}_+(z,T^{-d}\omega)F_z^+(0,\omega)\\
&+(F_z^+(0,\omega))^{-1}G(z,\omega)F_z^+(0,\omega) \\
=&\widetilde{M}_+(z,T^{-d}\omega)(F_z^+(0,T^{-d}\omega))^{-1}G(z,T^{-d}\omega)\frac{\partial CM_+(z,T^{-d}\omega)}{ \partial z}F_z^+(0,T^{-d}\omega)\widetilde{M}^{-1}_+(z,T^{-d}\omega)\\
&- (F_z^+(0,T^{-d}\omega))^{-1}\frac{\partial  M^{-1}_+(z,T^{-d}\omega)}{\partial z}  M_+(z,T^{-d}\omega)F_z^+(0,T^{-d}\omega)\\
&+(F_z^+(0,\omega))^{-1}G(z,\omega)F_z^+(0,\omega).
\end{align*}
By Lemma \ref{wewe1},
$$
f(z,\omega)-f(z,T^{-d}\omega)=-\frac{\partial \tau(z,T^{-d}\omega)}{\partial z}\frac{1}{\tau(z,T^{-d}\omega)}-h(z,T^{-d}\omega)+h(z,\omega)+g(z,\omega).
$$
\end{pf}
Finally, the integral of both sides of \eqref{ft} over $\Omega$
\begin{align*}
\int_\Omega\frac{\partial \tau(z,T^{-d}\omega)}{\partial z}\frac{1}{\tau(z,T^{-d}\omega)}d\mu =\int_\Omega g(z,\omega)d\mu
\end{align*}
leads to \eqref{gi} and thus we get the desired result.
\end{pf}

\subsection{Kotani theoretic estimates} Recall that for any $z\in
\mathbb{H}_\delta$, there are non-zero solutions to
$L_{f,\omega}u=zu$,  $(u_z^\pm(n,\omega))_{n\in\Z}\in E_z^\pm(\omega)$
respectively, that are $\ell^2$ at $\pm \infty$. By \eqref{u(n)} we have 
$$
\vec{u}_z^\pm(n,\omega)=\begin{pmatrix} u_z^\pm(nd+d-1,\omega) \\ \vdots \\ u_z^\pm(nd,\omega)\end{pmatrix}.
$$
\begin{Lemma}\label{le1}
We have
\begin{align*}
&\int_\Omega \ln\left(1-\frac{\Im z\|\vec{u}^+_z(0,\omega)\|^2}{\Im \left(\vec{u}^+_z(0,\omega)\right)^*C\vec{u}^+_z(1,\omega)}\right)d\mu=2dL_f^d(z),
\end{align*}
\begin{align*}
&\int_\Omega \ln\left(1-\frac{\Im z\|\vec{u}^-_z(0,\omega)\|^2}{\Im \left(\vec{u}_z^-(0,\omega)\right)^*C^*\vec{u}_z^-(-1,\omega)}\right)d\mu=2dL_f^d(z).
\end{align*}
\end{Lemma}
\begin{pf}
Notice that $\{\vec{u}_z^\pm(n,\omega)\}_{n\in\Z}$ are solutions of the following equation
\begin{equation}\label{efsch}
C^*\vec{u}(n-1,\omega)+C\vec{u}(n+1,\omega)+B(T^{dn}\omega)\vec{u}(n,\omega)=z\vec{u}(n,\omega).
\end{equation}
It follows from \eqref{efsch} that
\begin{equation}\label{w1}
\left(\vec{u}_z^\pm(0,\omega)\right)^*C^*\vec{u}_z^\pm(-1,\omega)+\left(\vec{u}_z^\pm(0,\omega)\right)^*C\vec{u}_z^\pm(1,\omega)+\left(\vec{u}_z^\pm(0,\omega)\right)^*(B(\omega)-z)\vec{u}_z^\pm(0,\omega)=0.
\end{equation}
Let
\begin{equation}\label{w2}
c_+(\omega)=\left\|\left(\vec{u}_z^+(0,\omega)\right)^*C\vec{u}_z^+(0,T^{d}\omega)\right\|,
\end{equation}
\begin{equation}
c_-(\omega)=\left\|\left(\vec{u}_z^-(0,\omega)\right)^*C^*\vec{u}_z^-(0,T^{-d}\omega)\right\|,
\end{equation}
\begin{equation}\label{ww3}
m_+(\omega)=-\frac{\left(\vec{u}_z^+(0,\omega)\right)^*C\vec{u}_z^+(1,\omega)}{\left\|\left(\vec{u}_z^+(0,\omega)\right)C\vec{u}_z^+(0,T^{d}\omega)\right\|},
\end{equation}
\begin{equation}\label{w3}
m_-(\omega)=-\frac{\left(\vec{u}_z^-(0,\omega)\right)^*C^*\vec{u}_z^-(-1,\omega)}{\left\|\left(\vec{u}_z^-(0,\omega)\right)C^*\vec{u}_z^-(0,T^{-d}\omega)\right\|}.
\end{equation}
Multiplying both sides of equation \eqref{efsch} by $\vec{u}^*(n,\omega)$, taking the imaginary part and summing all the terms of each side, we get
\begin{equation*}
\Im \vec{u}^*(0,\omega)C\vec{u}(1,\omega)=-\Im z\sum\limits_{n=1}^{\infty}\|\vec{u}(n,\omega)\|^2,
\end{equation*}
\begin{equation*}
\Im \vec{u}^*(0,\omega)C\vec{u}(-1,\omega)=-\Im z\sum\limits_{n=-\infty}^{-1}\|\vec{u}(n,\omega)\|^2.
\end{equation*}
Thus for $\{\vec{u}_z^{\pm}(k,\omega)\}_{k\in\Z}$, we have
\begin{equation}\label{+}
\Im \left(\vec{u}_z^{+}(0,\omega)\right)^*C\vec{u}_z^+(1,\omega)=-\Im z\sum\limits_{n=1}^{\infty}\|\vec{u}_z^+(n,\omega)\|^2,
\end{equation}
\begin{equation}\label{-}
\Im \left(\vec{u}_z^{-}(0,\omega)\right)^*C^*\vec{u}_z^-(-1,\omega)=-\Im z\sum\limits_{n=-\infty}^{-1}\|\vec{u}_z^-(n,\omega)\|^2,
\end{equation}
It follows from \eqref{ww3},\eqref{w3},\eqref{+}, and \eqref{-} that
$$
\Im m_\pm(\omega)>0.
$$

On the other hand, by the invariance, there are $\tau_\pm(\omega)$ such that
$$
\vec{u}_z^\pm(1,\omega)=\vec{u}_z^\pm(0,T^{d}\omega)\tau_\pm(\omega),\ \ \vec{u}_z^\pm(-1,\omega)=\vec{u}_z^\pm(0,T^{-d}\omega)\tau_{\pm}^{-1}(T^{-d}\omega),
$$
which means that
\begin{equation}\label{w4}
m_+(\omega)=-\tau_+(\omega)\frac{\left(\vec{u}_z^+(0,\omega)\right)^*C\vec{u}_z^+(0,T^d\omega)}{\left\|\left(\vec{u}_z^+(0,\omega)\right)^*C\vec{u}_z^+(0,T^d\omega)\right\|},
\end{equation}
\begin{equation}\label{w4-}
m_-(\omega)=-\tau^{-1}_-(T^{-d}\omega)\frac{\left(\vec{u}_z^-(0,\omega)\right)^*C^*\vec{u}_z^-(0,T^{-d}\omega)}{\left\|\left(\vec{u}_z^-(0,\omega)\right)^*C^*\vec{u}_z^-(0,T^{-d}\omega)\right\|}.
\end{equation}
By \eqref{w1}-\eqref{w4-}, we have
\begin{equation*}
-c_+(T^{-d}\omega)m_+^{-1}(T^{-d}\omega)-c_+(\omega)m_+(\omega)+\left(\vec{u}_z^+(0,\omega)\right)^*(B(\omega)-z)\vec{u}_z^+(0,\omega)=0,
\end{equation*}
\begin{equation*}
-c_-(T^{d}\omega)m_-^{-1}(T^{d}\omega)-c_-(\omega)m_-(\omega)+\left(\vec{u}_z^-(0,\omega)\right)^*(B(\omega)-z)\vec{u}_z^-(0,\omega)=0.
\end{equation*}
Taking the imaginary part, one has
\begin{equation}\label{equal}
-c_\pm(T^{\mp d}\omega)\frac{\Im m_\pm(T^{\mp d}\omega)}{|m_\pm(T^{\mp d}\omega)|^2}+c_\pm(\omega)\Im m_\pm(\omega)+\Im z\|u_z^\pm(0,\omega)\|^2=0.
\end{equation}
Thus
\begin{align*}
&\ln\left(1+\frac{\Im z \|u_z^\pm(0,\omega)\|^2}{c_\pm(\omega)\Im m_\pm(z,\omega)}\right)=\ln c_\pm(T^{\mp d}\omega)-\ln c_\pm(\omega)\\
&+\ln \Im m_\pm(T^{\mp d}\omega)-\ln \Im m_\pm(\omega)-2\ln |m_\pm(T^{\mp d}\omega)|,
\end{align*}
and it follows that
\begin{align}\label{w7}
\ln\left(1-\frac{\Im z \|u_z^\pm(0,\omega)\|^2}{\Im \left(\vec{u}^+_z(0,\omega)\right)^*C\vec{u}^+_z(\pm 1,\omega)}\right)&=-2\int_\Omega \ln |m_\pm(T^{\mp d}\omega)|d\mu\\ \nonumber
=&\mp2\int_\Omega \ln |\tau_\pm(\omega)|d\mu.
\end{align}
Finally, by dominated splitting and the definition of $\{u_z^\pm(n,\omega)\}_{n\in\Z}$,  we have for $\mu$-almost every $\omega$,
$$
\lim\limits_{n\rightarrow\infty}\frac{1}{n}\ln \frac{\|u_z^\pm(n,\omega)\|}{\|u_z^\pm(0,\omega)\|}=\mp dL_d^f(z),
$$
By the invariance,
$$
\ln \frac{\|u_z^\pm(n,\omega)\|}{\|u_z^\pm(0,\omega)\|}=\sum\limits_{m=0}^{n-1}\ln|\tau_\pm(T^{md}\omega)|,
$$
and hence Birkhoff's ergodic theorem implies,
\begin{equation}\label{w8}
\int_\Omega \ln|\tau_\pm(\omega)|d\mu=\mp dL_f^d(z).
\end{equation}
\eqref{w7} and \eqref{w8} completes the proof.
\end{pf}
\begin{Lemma}\label{le2}
We have that
\begin{align*}
&\int_\Omega \frac{1}{-\Im \frac{\left(\vec{u}^+_z(0,\omega)\right)^*C\vec{u}^+_z(1,\omega)}{\|\vec{u}^+_z(0,\omega)\|^2}+\frac{1}{2}\Im z}d\mu\leq \frac{2dL_f^d(z)}{\Im z},
\end{align*}
\begin{align*}
&\int_\Omega \frac{1}{-\Im \frac{\left(\vec{u}^-_z(0,\omega)\right)^*C\vec{u}^-_z(-1,\omega)}{\|\vec{u}^-_z(0,\omega)\|^2}+\frac{1}{2}\Im z}d\mu\leq \frac{2dL_f^d(z)}{\Im z},
\end{align*}
\end{Lemma}
\begin{pf}
For $x\geq 0$, consider the function
$$
A(x)=\ln (1+x)-\frac{x}{1+\frac{x}{2}}.
$$
Clearly,
$$
A(0)=0,\ \ A'(x)=\frac{1}{1+x}-\frac{1}{1+x+\frac{x^2}{4}}\geq 0.
$$
Hence
\begin{equation}\label{ff3}
\ln(1+x)\geq \frac{x}{1+\frac{x}{2}},\  \ \forall x\geq 0,
\end{equation}
By \eqref{ff3} and Lemma \ref{le1}, we have
\begin{align*}
\int_\Omega \frac{1}{-\Im \frac{\left(\vec{u}^+_z(0,\omega)\right)^*C\vec{u}^+_z(1,\omega)}{\|\vec{u}^+_z(0,\omega)\|^2}+\frac{1}{2}\Im z}d\mu&=\frac{1}{\Im z}\int_\Omega \frac{-\frac{\Im z}{\Im \frac{\left(\vec{u}^+_z(0,\omega)\right)^*C\vec{u}^+_z(1,\omega)}{\|\vec{u}^+_z(0,\omega)\|^2}}}{1+\frac{\Im z}{-2\Im \frac{\left(\vec{u}^+_z(0,\omega)\right)^*C\vec{u}^+_z(1,\omega)}{\|\vec{u}^+_z(0,\omega)\|^2}}}d\mu\\
&\leq \int_\Omega \ln\left(1+\frac{-\Im z\|\vec{u}^+_z(0,\omega)\|^2}{\Im \left(\vec{u}^+_z(0,\omega)\right)^*C\vec{u}^+_z(1,\omega)}\right)d\mu\\
&=\frac{2dL_f^d(z)}{\Im z},
\end{align*}
The proof for $\vec{u}^-_z(\omega)$ follows in exactly the same way.
\end{pf}

Recall that by the $PH2$ property, for any $z\in\C_\delta$, there exist continuous invariant decompositions
$$
\C^{2d}=E_z^s(\omega)\oplus E_z^c(\omega)\oplus E_z^u(\omega),\ \ \forall \omega\in\Omega,
$$
where $E_z^c(\omega)$ is the two dimensional invariant subspace corresponding to the minimal Lyapunov exponent and moreover $E_z^s(\omega)$, $E_z^c(\omega)$, and $E_z^u(\omega)$ depend continuously on $\omega$ and analytically on $z$.

We are now ready to define the finite-range analogue of the
$m$-function. Note that when $\Im z=0$, involving  the complex
symplectic structure, we  actually have a {\it symplectic} continuous invariant decomposition of $\R^{2d}$.
\begin{Lemma}\label{uvec}
 For any $z\in \mathbb{C}_\delta$, there exist $u_z(\omega), v_z(\omega)\in E_z^c(\omega)$, such that
\begin{enumerate}
\item $u_z(\omega)$ and $v_z(\omega)$ depend continuously on $\omega$ and analytically on $z\in \C_\delta$,
\item $\begin{pmatrix}\vec{u}^*_E(0,\omega)&\vec{u}^*_E(-1,\omega)\end{pmatrix}S\begin{pmatrix}\vec{v}_E(0,\omega)\\ \vec{v}_E(-1,\omega)\end{pmatrix}\neq 0$ for any $\omega\in \Omega$ and $E\in\R$ where
$$
S=\begin{pmatrix}0&-C^*\\ C&0\end{pmatrix},\ \  \begin{pmatrix}\vec{u}_z(0,\omega)\\ \vec{u}_z(-1,\omega)\end{pmatrix}=u_z(\omega),\ \  \begin{pmatrix}\vec{v}_z(0,\omega)\\ \vec{v}_z(-1,\omega)\end{pmatrix}=v_z(\omega).
$$
\end{enumerate}
\end{Lemma}
\begin{pf}
  Let $s$ be a symplectic form on $\R^{2d}$ with symplectic inner product
$$
s(u,v)=u^*S v, \ \ u,v\in\R^{2d}.
$$
Given a subspace $V\subset \R^{2d}$, we denote its $s$-orthogonal complement $V^\bot$ which is defined by those vectors $u\in\R^{2d}$ such that $s(u,v)=0$ for any $v\in V$. $V$ is called a symplectic subspace  if $V\cap V^\bot=\{0\}$.
We therefore only need to prove that $E_z^c(\omega)$ is a symplectic
subspace for $\Im z=0$ that is $E_z^c(\omega)\cap
(E_z^c(\omega))^{\bot}=\{0\}$.

Given a nonzero vector $v_1\in E_z^c(\omega)$, there exists
$\bar{v}_1\in \R^{2d}$ such that $s(v_1,\bar{v}_1)\neq 0$ (otherwise
$v_1=0$). Clearly $\bar{v}_1\in E_z^c(\omega)$. We complete  the
symplectic basis of $E_z^c(\omega)$ obtaining $\{v_1,\bar{v}_1\}$. Let
$u\in E_z^c(\omega)\cap (E_z^c(\omega))^\bot.$ Then since $u\in
E_z^c(\omega)$, we have $u=c_1 v_1+c_2 \bar{v}_1$ for some
$c_1,c_2\in\R$. On the other hand, since $u\in (E_z^c(\omega))^\bot$,
we have $0=s(u,v_1)=s(c_1 v_1+c_2 \bar{v}_1,v_1)=c_2s(\bar{v}_1,v_1)$
which implies that $c_2=0$. Similarly, using that $s(u,\bar{v}_1)=0$, one has $c_1=0$. Thus $E_z^c(\omega)$ is a symplectic subspace.

Notice that $E_z^c(\omega)$ depends continuously on $\omega$ and $z$,
actually since $L_z^f$ depends analytically on $z$, the basis of $E_z^c(\omega)$ can be chosen to depend holomorphically on  $z\in\C_\delta$. i.e., there exist $u_z(\omega), v_z(\omega)\in E_z^c(\omega)$, depending continuously in $\omega$ and analytically on $z\in \C_\delta$, such that
$$
\begin{pmatrix}\vec{u}^*_E(0,\omega)&\vec{u}^*_E(-1,\omega)\end{pmatrix}S\begin{pmatrix}\vec{v}_E(0,\omega)\\ \vec{v}_E(-1,\omega)\end{pmatrix}\neq 0
$$
for any $\omega\in \Omega$ and $E\in\R$. 
\end{pf}

Note that for any $z\in \mathbb{H}_\delta$, there are $b^\pm_1(z,\omega)$ and $b^\pm_0(z,\omega)$ \footnote{Since $u_z(\omega)$ and $v_z(\omega)$ is a basis for $E_z^c(\omega)$.}, depending continuously on $\omega$ and analytically on $z$ such that
\begin{align}\label{ff1}
u_z^+(\omega)=b^+_1(z,\omega)u_z(\omega)+b^+_0(z,\omega)v_z(\omega),
\end{align}
\begin{align}\label{ff2}
u_z^-(\omega)=b^-_1(z,\omega)u_z(\omega)+b^-_0(z,\omega)v_z(\omega).
\end{align}
satisfy $u_z^\pm(\omega)\in E_z^\pm(\omega)$.
\begin{Definition}[m functions]\label{m+def}
{\rm We define
$$
m_+(z,\omega)=\frac{b^+_1(z,\omega)}{b_0^+(z,\omega)},\ \ m_-(z,\omega)=\frac{b^-_0(z,\omega)}{b_1^-(z,\omega)}.
$$
}
\end{Definition}
For simplicity, we also set
\begin{align}\label{ce}
c_E(\omega)&=\begin{pmatrix}\vec{u}^*_E(1,\omega)&\vec{u}^*_E(0,\omega)\end{pmatrix}S\begin{pmatrix}\vec{v}_E(1,\omega)\\ \vec{v}_E(0,\omega)\end{pmatrix}\\ \nonumber
&=-\vec{u}^*_E(1,\omega)C^*\vec{v}_E(0,\omega)+\vec{u}^*_E(0,\omega)C\vec{v}_E(1,\omega).
\end{align}
The following $PH2$ analogue of Kotani theory is key to this work.
\begin{Theorem}[Kotani-theoretic estimates]\label{keyth}
For $z=E+i\delta$, we have that
\begin{enumerate}
\item $m_\pm(E+i0,\omega):=\lim\limits_{\delta\rightarrow 0}m_\pm(z,\omega)$ exist for almost every $E$ and $\mu$-almost every $\omega$.
\item For almost every $E\in\Sigma_f$, we have
\begin{align*}\int_\Omega \frac{\left\|m_+(E+i0,\omega)\begin{pmatrix}
\vec{u}_E(0,\omega)\\
\vec{u}_E(-1,\omega)
\end{pmatrix}+\begin{pmatrix}
\vec{v}_E(0,\omega)\\
\vec{v}_E(-1,\omega)
\end{pmatrix}\right\|^2}{\Im \left(m_+(E+i0,\omega)c_E(\omega)\right)}d\mu<\infty,
\end{align*}
\begin{align*}\int_\Omega \frac{\left\|\begin{pmatrix}
\vec{u}_E(0,\omega)\\
\vec{u}_E(-1,\omega)
\end{pmatrix}+m_-(E+i0,\omega)\begin{pmatrix}
\vec{v}_E(0,\omega)\\
\vec{v}_E(-1,\omega)
\end{pmatrix}\right\|^2}{\Im \left(m_-(E+i0,\omega)c_E(\omega)\right)}d\mu<\infty.
\end{align*}
\item For $\mu$-almost every $\omega$, we have
\begin{align*}
m_+(E+i0,\omega)=\frac{1}{\overline{m_-(E+i0,\omega)}},
\end{align*}
for almost every $E\in\Sigma_f$.
\end{enumerate}
\end{Theorem}
\begin{pf}
For (1) and (2), we only give the proof of the results for
$m_+(z,\omega)$ since the proof for  $m_-(z,\omega)$ is exactly the
same.  Recall that for any $z\in\mathbb{H}_\delta$, there are linearly
independent $\{\vec{u}_z^i(\omega)\}_{i=1}^{d-1}\in E_z^s(\omega)$
depending continuously on $\omega$ and analytically on $z$. It is also
easy to check that for any $z\in\mathbb{H}_\delta$,
\begin{equation}\label{r100}
\widetilde{F}_z^+(\omega)=\begin{pmatrix}\vec{u}_z^1(\omega),\cdots, \vec{u}_z^{d-1}(\omega),m_+(z,\omega)\vec{u}_z(\omega)+\vec{v}_z(\omega)\end{pmatrix}
\end{equation}
satisfies
$$
\sum\limits_{k=0}^\infty\|\widetilde{F}_z^+(k,\omega)\|^2<\infty,
$$
where
\begin{align*}
&\begin{pmatrix}\vec{u}_z^1(k,\omega)&\cdots&\vec{u}_z^{d-1}(k,\omega)&m_+(z,\omega)\vec{u}_z(k,\omega)+\vec{v}_z(k,\omega)\\ \vec{u}_z^1(k,\omega)&\cdots&\vec{u}_z^{d-1}(k-1,\omega)&m_+(z,\omega)\vec{u}_z(k-1,\omega)+\vec{v}_z(k-1,\omega)\end{pmatrix}\\
=&\begin{pmatrix}\widetilde{F}_z^+(k,\omega)\\ \widetilde{F}_z^+(k-1,\omega)\end{pmatrix}=(L_z^f)_{dk}(\omega)\widetilde{F}_z^+(\omega).
\end{align*}
Hence
\begin{equation}\label{r102}
M_+(z,\omega)=\widetilde{F}_z^+(1,\omega)\left(\widetilde{F}_z^+(0,\omega)\right)^{-1}.
\end{equation}
Using the existence of $\lim_{\delta\rightarrow 0^+}M_+(z,\omega)$ for
almost $E$ and $\mu$-almost every $\omega$ (by the property of
Herglotz functions), by \eqref{r100} and \eqref{r102}, we have
$$
\langle \delta_d, M_+(z,\omega)\delta_d\rangle=\frac{a(z,\omega)m_+(z,\omega)+b(z,\omega)}{c(z,\omega)m_+(z,\omega)+d(z,\omega)}
$$
where $a(z,\omega), b(z,\omega),c(z,\omega), d(z,\omega)$ depend analytically on $z$ and, moreover,
$$
c(z,\omega)=\det{\begin{pmatrix}\vec{u}_z^1(0,\omega),\cdots, \vec{u}_z^{d-1}(0,\omega),\vec{u}_z(0,\omega)\end{pmatrix}}\neq 0,
$$
which implies that
$$
m_+(z,\omega)=\frac{b(z,\omega)c(z,\omega)-a(z,\omega)d(z,\omega)}{c^2(z,\omega)\langle \delta_d, M_+(z,\omega)\delta_d\rangle-a(z,\omega)c(z,\omega)}-\frac{d(z,\omega)}{c(z,\omega)}.
$$
Thus $m_+(E+i0,\omega)=\lim_{\delta\rightarrow 0}m_+(z,\omega)$ exists
for almost every $E$ and $\mu$-almost every $\omega$. This completes the proof of part (1).

For part (2), multiplying both sides of equation \eqref{efsch} by $\vec{u}^*(n,\omega)$, taking the imaginary part and summing all the terms on each side, we obtain
\begin{equation*}
\Im \vec{u}^*(0,\omega)C\vec{u}(1,\omega)=-\Im z\sum\limits_{n=1}^{\infty}\|\vec{u}(n,\omega)\|^2.
\end{equation*}
Therefore for
$\vec{u'}_z^{+}(k,\omega)=m_+(z,\omega)u_z(k,\omega)+v_z(k,\omega)$,
we have
\begin{equation}\label{q1}
\Im \left(\vec{u'}_z^{+}(0,\omega)\right)^*C\vec{u’}_z^+(1,\omega)=-\Im z\sum\limits_{n=1}^{\infty}\|\vec{u'}_z^+(n,\omega)\|^2,
\end{equation}
It follows that
\begin{align*}
&\left(m_+(z,\omega)\vec{u}_z(0,\omega)+\vec{v}_z(0,\omega)\right)^*C\left(m_+(z,\omega)\vec{u}_z(1,\omega)+\vec{v}_z(1,\omega)\right)\\
=&|m_+(z,\omega)|^2\vec{u}^*_z(0,\omega)C\vec{u}_z(1,\omega)+\vec{v}^*_z(0,\omega)C\vec{v}_z(1,\omega)\\
&+\overline{m_+(z,\omega)}\vec{u}^*_z(0,\omega)C\vec{v}_z(1,\omega)+m_+(z,\omega)\vec{v}^*_z(0,\omega)C\vec{u}_z(1,\omega)
\end{align*}
Taking the imaginary part and letting $\delta\rightarrow 0$, we have
\begin{align*}
&\Im\left(m_+(E+i0,\omega)\vec{u}_E(0,\omega)+\vec{v}_E(0,\omega)\right)^*C\left(m_+(E+i0,\omega)\vec{u}_E(1,\omega)+\vec{v}_E(1,\omega)\right)\\
=&\Im \left(m_+(E+i0,\omega)(\vec{v}^*_E(0,\omega)C\vec{u}_E(1,\omega)-\vec{u}^*_E(0,\omega)C\vec{v}_E(1,\omega))\right)\\
=&\Im \left[m_+(E+i0,\omega)(\vec{v}^*_E(0,\omega)C\vec{u}_E(1,\omega)-\vec{v}^*_E(1,\omega)C^*\vec{u}_E(0,\omega))\right].
\end{align*}
Note that
$$
\lim_{\delta\rightarrow 0^+}\frac{L_d^f(z)}{\Im z}=\lim_{\delta\rightarrow 0^+}\frac{\partial L_d^f(z)}{\partial \Im z}<\infty
$$
for almost every $E$. Thus by applying Lemma \ref{le2} to $\vec{u'}_z^{+}(\omega)$ and Fatou's lemma, for almost every $E$, we have
\begin{align}\label{ff2}
&\int_\Omega \frac{\|m_+(E+i0,\omega)\vec{u}_E(0,\omega)+\vec{v}_E(0,\omega)\|^2}{\Im \left(m_+(E+i0,\omega)c_E(\omega)\right)}d\mu<\infty.
\end{align}
By invariance, we can write
$$
\vec{u'}_E^+(0,\omega)=\vec{u'}_E^+(-1,T^{d}\omega)\tau'_+(\omega).
$$
By \eqref{equal}, for almost every $E$, we have
\begin{align}\label{ff3}
&\int_\Omega \frac{\|m_+(E+i0,\omega)\vec{u}_E(-1,\omega)+\vec{v}_E(-1,\omega)\|^2}{\Im \left(m_+(E+i0,\omega)c_E(\omega)\right)}d\mu\\  \nonumber
=&\int_\Omega \frac{\|m_+(E+i0,\omega)\vec{u}_E(-1,\omega)+\vec{v}_E(-1,\omega)\|^2|\tau'_+(T^{-d}\omega)|^2}{\Im \left(m_+(E+i0,T^{-d}\omega)c_E(T^{-d}\omega)\right)}d\mu\\ \nonumber
=&\int_\Omega \frac{\|m_+(E+i0,T^d\omega)\vec{u}_E(-1,T^d\omega)+\vec{v}_E(-1,T^d\omega)\|^2|\tau'_+(\omega)|^2}{\Im \left(m_+(E+i0,\omega)c_E(\omega)\right)}d\mu\\  \nonumber
=&\int_\Omega \frac{\|m_+(E+i0,\omega)\vec{u}_E(0,\omega)+\vec{v}_E(0,\omega)\|^2}{\Im \left(m_+(E+i0,\omega)c_E(\omega)\right)}d\mu
\end{align}
By \eqref{ff2} and \eqref{ff3}, we have
$$\int_\Omega \frac{\left\|m_+(E+i0,\omega)\begin{pmatrix}
\vec{u}_E(0,\omega)\\
\vec{u}_E(-1,\omega)
\end{pmatrix}+\begin{pmatrix}
\vec{v}_E(0,\omega)\\
\vec{v}_E(-1,\omega)
\end{pmatrix}\right\|^2}{\Im \left(m_+(E+i0,\omega)c_E(\omega)\right)}d\mu<\infty,
$$
completing the proof of (2).

Finally, we prove (3).  Taking the imaginary part on each side of \eqref{q1}, for almost every $E$ and $\mu$-almost every $\omega$, we have
$$
\Im m_+(E+i0,\omega)c_E(\omega)\geq 0,
$$
$$
\Im m_-(E+i0,\omega)c_E(\omega)\geq 0.
$$
Part (2) further implies
\begin{align}\label{ffff1}
\Im m_+(E+i0,\omega)c_E(\omega)> 0,\ \
\Im m_-(E+i0,\omega)c_E(\omega)>0.
\end{align}
for $\mu$-almost every $\omega\in\Omega$.

Let $A_E\in C^\omega(\T,GL(2,\R))$ be such that
\begin{equation}\label{ff6}
L_E^f(\omega)\begin{pmatrix}\vec{u}_E(0,\omega)&\vec{v}_E(0,\omega)\\ \vec{u}_E(-1,\omega)&\vec{v}_E(-1,\omega)\end{pmatrix}=\begin{pmatrix}\vec{u}_E(0,T\omega)&\vec{v}_E(0,T\omega)\\ \vec{u}_E(-1,T\omega)&\vec{v}_E(-1,T\omega)\end{pmatrix}A_E(\omega).
\end{equation}
Let $T_E(\omega)=\begin{pmatrix}0&c_E(\omega)\\
  -c_E(\omega)&0\end{pmatrix}$. Then
\begin{Proposition}\label{p1}
We have
$$
A_E(\omega)^*T_E(T\omega)A_E(\omega)=T_E(\omega).
$$
\end{Proposition}
\begin{pf}
Taking the transpose on each side of equation \eqref{ff6}, we have
$$
\begin{pmatrix}\vec{u}^*_E(0,\omega)&\vec{u}^*_E(-1,\omega)\\ \vec{v}^*_E(0,\omega)&\vec{v}^*_E(-1,\omega)\end{pmatrix}\left(L_E^f(\omega)\right)^*=A^*_E(\omega)\begin{pmatrix}\vec{u}^*_E(0,T\omega)&\vec{u}^*_E(-1,T\omega)\\ \vec{v}^*_E(0,T\omega)&\vec{v}^*_E(-1,T\omega)\end{pmatrix}.
$$
Multiplying $S$ on both sides of the above equation, one has
$$
\begin{pmatrix}\vec{u}^*_E(0,\omega)&\vec{u}^*_E(-1,\omega)\\ \vec{v}^*_E(0,\omega)&\vec{v}^*_E(-1,\omega)\end{pmatrix}\left(L_E^f(\omega)\right)^*S=A^*_E(\omega)\begin{pmatrix}\vec{u}^*_E(0,T\omega)&\vec{u}^*_E(-1,T\omega)\\ \vec{v}^*_E(0,T\omega)&\vec{v}^*_E(-1,T\omega)\end{pmatrix}S.
$$
Involving the fact that
$$
\left(L_E^f(\omega)\right)^*S=S\left(L_E^f(\omega)\right)^{-1},
$$
it follows
$$
\begin{pmatrix}\vec{u}^*_E(0,\omega)&\vec{u}^*_E(-1,\omega)\\ \vec{v}^*_E(0,\omega)&\vec{v}^*_E(-1,\omega)\end{pmatrix}S=A^*_E(\omega)\begin{pmatrix}\vec{u}^*_E(0,T\omega)&\vec{u}^*_E(-1,T\omega)\\ \vec{v}^*_E(0,T\omega)&\vec{v}^*_E(-1,T\omega)\end{pmatrix}SL_E^f(\omega).
$$
Multiplying $\begin{pmatrix}\vec{u}_E(0,\omega)&\vec{v}_E(0,\omega)\\ \vec{u}_E(-1,\omega)&\vec{v}_E(-1,\omega)\end{pmatrix}$ on the right of each side,
\begin{align*}
&\begin{pmatrix}
\vec{u}^*_E(0,\omega)&\vec{u}^*_E(-1,\omega)\\ \vec{v}^*_E(0,\omega)&\vec{v}^*_E(-1,\omega)\end{pmatrix}S\begin{pmatrix}\vec{u}_E(0,\omega)&\vec{v}_E(0,\omega)\\ \vec{u}_E(-1,\omega)&\vec{v}_E(-1,\omega)\end{pmatrix}\\
=&A_E(\omega)^*\begin{pmatrix}\vec{u}^*_E(0,T\omega)&\vec{u}^*_E(-1,T\omega)\\ \vec{v}^*_E(0,T\omega)&\vec{v}^*_E(-1,T\omega)\end{pmatrix}SL_E^f(\omega)\begin{pmatrix}\vec{u}_E(0,\omega)&\vec{v}_E(0,\omega)\\ \vec{u}_E(-1,\omega)&\vec{v}_E(-1,\omega)\end{pmatrix}\\
=&A_E(\omega)^*\begin{pmatrix}\vec{u}^*_E(0,T\omega)&\vec{u}^*_E(-1,T\omega)\\ \vec{v}^*_E(0,T\omega)&\vec{v}^*_E(-1,T\omega)\end{pmatrix}S\begin{pmatrix}\vec{u}_E(0,T\omega)&\vec{v}_E(0,T\omega)\\ \vec{u}_E(-1,T\omega)&\vec{v}_E(-1,T\omega)\end{pmatrix}A_E(\omega).
\end{align*}
Finally one can easily check that
$$
T_E(\omega)=\begin{pmatrix}
\vec{u}^*_E(0,\omega)&\vec{u}^*_E(-1,\omega)\\ \vec{v}^*_E(0,\omega)&\vec{v}^*_E(-1,\omega)\end{pmatrix}S\begin{pmatrix}\vec{u}_E(0,\omega)&\vec{v}_E(0,\omega)\\ \vec{u}_E(-1,\omega)&\vec{v}_E(-1,\omega)\end{pmatrix},
$$
thus completing the proof.
\end{pf}
We need the following lemma,
\begin{Lemma}\label{final1}
For almost every E,
\begin{align*}
&\lim\limits_{\Im z\rightarrow 0^+}\frac{\partial L^f_d(z)}{\partial\Im z}\\
=&-\frac{1}{d}\int_\Omega\Im \frac{\left(m_+(E+i0,\omega)\vec{u}_E(0,\omega)+\vec{v}_E(0,\omega)\right)^T\left(\vec{u}_E(0,\omega)+\vec{v}_E(0,\omega)m_-(E+i0,\omega)\right)}{c_E(\omega)(1-m_+(E+i0,\omega)m_-(E+i0,\omega))}d\mu.
\end{align*}
\end{Lemma}
\begin{pf}
As before, for any $z\in\mathbb{H}_\delta$, there exist continuous invariant decompositions
$$
\C^{2d}=E_z^s(\omega)\oplus E_z^+(\omega)\oplus E_z^-(\omega)\oplus E_z^u(\omega),\ \ \forall \omega\in\Omega,
$$
which implies that there are $\{\vec{u}_z^i(\omega)\}_{i=1}^{d-1}\in E_z^s(\omega)$ and $\{\vec{v}_z^i(\omega)\}_{i=1}^{d-1}\in E_z^u(\omega)$ depending continuously on $\omega$ and analytically on $z$, such that
\begin{equation}\label{r1}
\widetilde{F}_z^+(\omega)=\begin{pmatrix}\vec{u}_z^1(\omega),\cdots, \vec{u}_z^{d-1}(\omega),m_+(z,\omega)\vec{u}_z(\omega)+\vec{v}_z(\omega)\end{pmatrix}
\end{equation}
\begin{equation}\label{r2}
\widetilde{F}_z^-(\omega)=\begin{pmatrix}\vec{u}_z(\omega)+m_-(z,\omega)\vec{v}_z(\omega), \vec{v}_z^{d-1}(\omega),\cdots, \vec{v}_z^{1}(\omega)\end{pmatrix}
\end{equation}
satisfy
$$
\sum\limits_{k=0}^\infty\|\widetilde{F}_z^+(k,\omega)\|^2<\infty,\ \ \sum\limits_{k=0}^{-\infty}\|\widetilde{F}_z^-(k,\omega)\|^2<\infty,
$$
where
\begin{align*}
\begin{pmatrix}\vec{u}_z^1(k,\omega)&\cdots&\vec{u}_z^{d-1}(k,\omega)&m_+(z,\omega)\vec{u}_z(k,\omega)+\vec{v}_z(k,\omega)\\ \vec{u}_z^1(k,\omega)&\cdots&\vec{u}_z^{d-1}(k-1,\omega)&m_+(z,\omega)\vec{u}_z(k-1,\omega)+\vec{v}_z(k-1,\omega)\end{pmatrix}=(L_z^f)_{dk}(\omega)\widetilde{F}_z^+(\omega).
\end{align*}
By Theorem \ref{mg} and Lebesgue dominated convergence theorem,  for almost every E,
\begin{align*}
\lim\limits_{\Im z\rightarrow 0^+}\frac{L^f_d(z)}{\Im z}=-\frac{1}{d}\int_\Omega\langle \delta_d, (\widetilde{F}_{E+i0}^+(0,\omega))^{-1}G(E+i0,\omega)\widetilde{F}_{E+i0}^+(0,\omega)\delta_d\rangle d\mu.
\end{align*}
Let
$$
\Phi_{E}(\omega)=\begin{pmatrix}
\vec{u}_E^1(1,\omega)&\cdots &\vec{u}_E^{d-1}(1,\omega)& \vec{u}_E(1,\omega)&\vec{v}_E(1,\omega)&\vec{v}^{d-1}_E(1,\omega)&\cdots&\vec{v}_E^1(1,\omega)\\
\vec{u}_E^1(0,\omega)&\cdots &\vec{u}_E^{d-1}(0,\omega)& \vec{u}_E(0,\omega)&\vec{v}_E(0,\omega)&\vec{v}^{d-1}_E(0,\omega)&\cdots&\vec{v}_E^1(0,\omega)
\end{pmatrix}
$$
Using that $(L_E^f(\omega))^*S L_E^f(\omega)=S,$ one can check that for $i=1,\cdots,d-1$,
$$
\begin{pmatrix}
\vec{u}^*_E(1,\omega) &\vec{u}^*_E(0,\omega)
\end{pmatrix}S\begin{pmatrix}
\vec{v}^i_E(1,\omega) \\ \vec{v}^i_E(0,\omega)
\end{pmatrix}=\begin{pmatrix}
\vec{u}^*_E(1,\omega) &\vec{u}^*_E(0,\omega)
\end{pmatrix}S\begin{pmatrix}
\vec{u}^i_E(1,\omega) \\ \vec{u}^i_E(0,\omega)
\end{pmatrix}=0,
$$
$$
\begin{pmatrix}
\vec{v}^*_E(1,\omega) &\vec{v}^*_E(0,\omega)
\end{pmatrix}S\begin{pmatrix}
\vec{v}^i_E(1,\omega) \\ \vec{v}^i_E(0,\omega)
\end{pmatrix}=\begin{pmatrix}
\vec{v}^*_E(1,\omega) &\vec{v}^*_E(0,\omega)
\end{pmatrix}S\begin{pmatrix}
\vec{u}^i_E(1,\omega) \\ \vec{u}^i_E(0,\omega)
\end{pmatrix}=0.
$$
Thus there are $C_E^\pm(\omega)$ such that
$$
\Phi_E^*(\omega)S\Phi_E(\omega)=\begin{pmatrix}C_E^+(\omega)&&\\ &&T_E(\omega)&\\&&&C_E^-(\omega)\end{pmatrix}
$$
which implies
$$
\Phi^{-1}_E(\omega)=\begin{pmatrix}(C_E^+(\omega))^{-1}&&\\ &&T^{-1}_E(\omega)&\\&&&(C_E^-(\omega))^{-1}\end{pmatrix}\Phi_E^*(\omega)S.
$$
On the other hand, we have
\begin{align*}
&\begin{pmatrix}\left(\widetilde{F}_{E+i0}^+(1,\omega)-\widetilde{F}_{E+i0}^-(1,\omega)(\widetilde{F}_{E+i0}^-(0,\omega))^{-1}\widetilde{F}_{E+i0}^+(0,\omega)\right)^{-1}&*
\\ *&*
\end{pmatrix}\\
:=&\begin{pmatrix}
\widetilde{F}_{E+i0}^+(1,\omega)&\widetilde{F}_{E+i0}^-(1,\omega)\\
\widetilde{F}_{E+i0}^+(0,\omega)&\widetilde{F}_{E+i0}^-(0,\omega)
\end{pmatrix}^{-1}\\
=&\left(\Phi_E(\omega)\begin{pmatrix}I_{d-1}&&&\\ &&m_+(E+i0,\omega)&1&\\&&1&m_-(E+i0,\omega)&\\&&&&I_{d-1}\end{pmatrix}\right)^{-1}\\
=&\begin{pmatrix}I_{d-1}&&&\\ &&\frac{m_-(E+i0,\omega)}{m_+(E+i0,\omega)m_-(E+i0,\omega)-1}&-\frac{1}{m_+(E+i0,\omega)m_-(E+i0,\omega)-1}&\\&&-\frac{1}{m_+(E+i0,\omega)m_-(E+i0,\omega)-1}&\frac{m_+(E+i0,\omega)}{m_+(E+i0,\omega)m_-(E+i0,\omega)-1}&\\&&&&I_{d-1}\end{pmatrix}\Phi^{-1}_E(\omega).
\end{align*}
Notice that
\begin{align}\label{aaa}
&\langle \delta_d,(\widetilde{F}_{E+i0}^+(0,\omega))^{-1}G(E+i0,\omega)\widetilde{F}_{E+i0}^+(0,\omega)\delta_d\rangle\\ \nonumber
=&\langle\delta_d, \begin{pmatrix}
\widetilde{F}_{E+i0}^+(1,\omega)&\widetilde{F}_{E+i0}^-(1,\omega)\\
\widetilde{F}_{E+i0}^+(0,\omega)&\widetilde{F}_{E+i0}^-(0,\omega)
\end{pmatrix}^{-1}\begin{pmatrix}C^{-1}\widetilde{F}_{E+i0}^+(0,\omega)&\\ & I_d\end{pmatrix}\delta_d\rangle,
\end{align}
\begin{align}\label{aaa1}
&\Phi_E^*(\omega)S\begin{pmatrix}C^{-1}\widetilde{F}_{E+i0}^+(0,\omega)&\\ & I_d\end{pmatrix}\\ \nonumber
=&\begin{pmatrix}\begin{pmatrix}
\vec{u}_E^1(0,\omega)&\cdots &\vec{u}_E^{d-1}(0,\omega)& \vec{u}_E(0,\omega)\end{pmatrix}^*\widetilde{F}_{E+i0}^+(0,\omega)&*
\\ \begin{pmatrix}
\vec{v}_E(0,\omega)&\vec{v}^{d-1}_E(0,\omega)&\cdots&\vec{v}_E^1(0,\omega)
\end{pmatrix}^*\widetilde{F}_{E+i0}^+(0,\omega)&*
\end{pmatrix}
\end{align}
\begin{align}\label{aaa2}
\nonumber &\begin{pmatrix}I_{d-1}&&&\\ &&\frac{m_-(E+i0,\omega)}{m_+(E+i0,\omega)m_-(E+i0,\omega)-1}&-\frac{1}{m_+(E+i0,\omega)m_-(E+i0,\omega)-1}&\\&&-\frac{1}{m_+(E+i0,\omega)m_-(E+i0,\omega)-1}&\frac{m_+(E+i0,\omega)}{m_+(E+i0,\omega)m_-(E+i0,\omega)-1}&\\&&&&I_{d-1}\end{pmatrix}\\
&\cdot \begin{pmatrix}(C_E^+(\omega))^{-1}&&\\ &&T^{-1}_E(\omega)&\\&&&(C_E^-(\omega))^{-1}\end{pmatrix}\\ \nonumber
=& \begin{pmatrix}(C_E^+(\omega))^{-1}&&&\\ &-\frac{1}{c_E(\omega)(m_+(E+i0,\omega)m_-(E+i0,\omega)-1)}&-\frac{m_-(E+i0,\omega)}{c_E(\omega)(m_+(E+i0,\omega)m_-(E+i0,\omega)-1)}&\\
&\frac{m_+(E+i0,\omega)}{c_E(\omega)(m_+(E+i0,\omega)m_-(E+i0,\omega)-1)}&\frac{1}{c_E(\omega)(m_+(E+i0,\omega)m_-(E+i0,\omega)-1)}&\\&&&(C_E^-(\omega))^{-1}\end{pmatrix}
\end{align}
It follows from \eqref{aaa}-\eqref{aaa2} that
\begin{align*}
&\langle \delta_d,(\widetilde{F}_{E+i0}^+(0,\omega))^{-1}G(E+i0,\omega)\widetilde{F}_{E+i0}^+(0,\omega)\delta_d\rangle\\
=&\frac{\left(\vec{u}_E(0,\omega)+\vec{v}_E(0,\omega)m_-(E+i0,\omega)\right)^T\left(m_+(E+i0,\omega)\vec{u}_E(0,\omega)+\vec{v}_E(0,\omega)\right)}{c_E(\omega)(1-m_+(E+i0,\omega)m_-(E+i0,\omega))}.
\end{align*}
\end{pf}
We are now ready to prove part (3).  For almost every $E$, we define
$$
\tilde{m}_+=m_+(E+i0,\omega),\ \ \tilde{m}_-=\frac{1}{m_-(E+i0,\omega)}.
$$
Omitting $\omega$  for simplicity, notice that
\begin{align*}
&\frac{\left(m_+(E+i0)\vec{u}_E(0)+\vec{v}_E(0)\right)^T\left(\vec{u}_E(0)+\vec{v}_E(0)m_-(E+i0)\right)}{c_E(1-m_+(E+i0)m_-(E+i0))}\\
=&\frac{\left(\vec{u}_E(0)\tilde{m}_++\vec{v}_E(0)\right)^T\left(\vec{u}_E(0)\tilde{m}_-+\vec{v}_E(0)\right)}{c_E(\tilde{m}_- -\tilde{m}_+)}\\
=&\frac{\left(\vec{u}_E(0)\tilde{m}_++\vec{v}_E(0)\right)^T\left(\vec{u}_E(0)\tilde{m}_-+\vec{v}_E(0)\right)\overline{(\tilde{m}_- -\tilde{m}_+)}}{c_E|\tilde{m}_- -\tilde{m}_+|^2}.
\end{align*}
Moreover, a direct calculation shows
\begin{align*}
&\Im \left(\vec{u}_E(0)\tilde{m}_++\vec{v}_E(0)\right)^T\left(\vec{u}_E(0)\tilde{m}_-+\vec{v}_E(0)\right)\overline{(\tilde{m}_- -\tilde{m}_+)}\\
=&\left(\|\vec{u}_E(0)\tilde{m}_-+\vec{v}_E(0)\|^2\right)\Im \tilde{m}_+-\left(\|\vec{u}_E(0)\tilde{m}_++\vec{v}_E(0)\|^2\right)\Im \tilde{m}_-.
\end{align*}
It follows that
\begin{align*}
& \int_\Omega\frac{\|m_+(E+i0)\vec{u}_E(0)+\vec{v}_E(0)\|^2}{c_E\Im m_+(E+i0)}+ \frac{\|\vec{u}_E(0)+m_-(E+i0)\vec{v}_E(0)\|^2}{c_E\Im m_-(E+i0)}\\
&+4\Im\frac{\left(m_+(E+i0)\vec{u}_E(0)+\vec{v}_E(0)\right)^T\left(\vec{u}_E(0)+\vec{v}_E(0)m_-(E+i0)\right)}{c_E(1-m_+(E+i0)m_-(E+i0))}d\mu\\
=&\int_\Omega\left(\frac{\|\vec{u}_E(0)\tilde{m}_++\vec{v}_E(0)\|^2}{c_E\Im \tilde{m}_+}-\frac{\|\vec{u}_E(0)\tilde{m}_-+\vec{v}_E(0)\|^2}{c_E\Im \tilde{m}_-}\right)\\
&\cdot\left(\frac{\left(\Re \left(\tilde{m}_--\tilde{m}_+\right)\right)^2+\left(\Im\left(\tilde{m}_-+\tilde{m}_+\right)\right)^2}{|\tilde{m}_+-\tilde{m}_-|^2}\right)d\mu\\
\leq& 0,
\end{align*}
by \eqref{ffff1}, we have $\tilde{m}_-=\overline{\tilde{m}_+}$ for almost every $\omega$. This finishes the proof of (3).

\end{pf}

\subsection{$L^2$-reducibility and proof of Theorems \ref{L2
    reducibility} and \ref{C0reducibility1}} We first prove the
$L^2$-reducibility theorem \ref{L2 reducibility} that we slightly reformulate as
\begin{Theorem}\label{L2 reducibility1}
For $PH2$ cocycles $(T,L_{E}^f)$ with minimal $T,$ for almost every
$E\in \Sigma_f^0:=\{E: L^f_d(E)=0\}$, there exist $U_E,V_E\in
L^2(\Omega,\R^{2d})$ and $R_E(\omega)\in SO(2,\R) $ such that
$$
L_E^f(\omega)(U_E(\omega),V_E(\omega))=(U_E(T\omega), V_E(T\omega))R_E(\omega),
$$
with
$$
H_E=(U_E,V_E)\in L^2(\Omega,Sp_{2d\times2}(\R)).
$$
\end{Theorem}
\begin{pf}
  Let $\vec{u}_E(i,\omega), \vec{v}_E(i,\omega), i=-1,0,$ be given by
  Lemma \ref{uvec} and $m_+$ be as defined in Definition \ref{m+def}.
For almost every $E\in \Sigma_f^0$, we set
$$
C_E(\omega)=\begin{pmatrix}
0&\frac{|m_+(E,\omega)|}{(\Im m_+(E,\omega))^{1/2}}\\
-\frac{(\Im m_+(E,\omega))^{1/2}}{|m_+(E,\omega)|}&\frac{\Re m_+(E,\omega)}{|m_+(E,\omega)|(\Im m_+(E,\omega))^{1/2}}
\end{pmatrix}.
$$
Let
$$
(U_E(\omega),V_E(\omega))=\frac{1}{\sqrt{c_E(\omega)}}\begin{pmatrix}
\vec{u}_E(0,\omega)&\vec{v}_E(0,\omega)\\
\vec{u}_E(-1,\omega)&\vec{v}_E(-1,\omega)
\end{pmatrix}C_E(\omega),
$$
then for $R_E(\omega)\in SO(2,\R) $ defined by
$$
L_E^f(\omega)(U_E(\omega),V_E(\omega))=(U_E(T\omega), V_E(T\omega))R_E(\omega),
$$
we have that $R_E(\omega)\in SO(2,\R) $, and moreover
\begin{align*}
\|U_E\|^2_{L^2}+\|V_E\|^2_{L^2}=&\left\|\frac{|m_+(E,\cdot)|}{(\Im m_+(E,\cdot))^{1/2}}\frac{\vec{u}_E(0,\cdot)}{\sqrt{c_E(\cdot)}}+\frac{\Re m_+(E,\cdot)}{|m_+(E,\cdot)|(\Im m_+(E,\cdot))^{1/2}}\frac{\vec{v}_E(0,\cdot)}{\sqrt{c_E(\cdot)}}\right\|^2_{L^2}\\
&+\left\|\frac{|m_+(E,\cdot)|}{(\Im m_+(E,\cdot))^{1/2}}\frac{\vec{u}_E(-1,\cdot)}{\sqrt{c_E(\cdot)}}+\frac{\Re m_+(E,\omega)}{|m_+(E,\cdot)|(\Im m_+(E,\cdot))^{1/2}}\frac{\vec{v}_E(-1,\cdot)}{\sqrt{c_E(\cdot)}}\right\|^2_{L^2}\\
&+\left\|\frac{(\Im m_+(E,\cdot))^{1/2}}{|m_+(E,\cdot)|}\frac{\vec{v}_E(0,\cdot)}{\sqrt{c_E(\cdot)}}\right\|_{L^2}^2+\left\|\frac{(\Im m_+(E,\cdot))^{1/2}}{|m_+(E,\cdot)|}\frac{\vec{v}_E(-1,\cdot)}{\sqrt{c_E(\cdot)}}\right\|_{L^2}^2\\
=&\int_\Omega\frac{\left\|m_+(E+i0,\cdot)\begin{pmatrix}
\vec{u}_E(0,\cdot)\\
\vec{u}_E(-1,\cdot)
\end{pmatrix}+\begin{pmatrix}
\vec{v}_E(0,\cdot)\\
\vec{v}_E(-1,\cdot)
\end{pmatrix}\right\|^2}{\Im m_+(E+i0,\cdot)c_E(\cdot)}d\mu<\infty,
\end{align*}
where the last inequality follows from part (2) of Theorem
\ref{keyth}. Finally, by the definition \eqref{ce} of $c_E$ and direct
calculation, we have for $S$ defined in Lemma \ref{uvec},
$$
H^*_E(\omega)SH_E(\omega)=\begin{pmatrix}0&1\\ -1&0\end{pmatrix}.
$$
\end{pf}

\noindent {\bf Proof of Theorem \ref{C0reducibility1}:}
The main idea is to involve the Schwartz reflection principle for functions on $\mathbb{H}^2$. Let $D$ be a Jordan domain such that
$D_\pm\subset \mathbb{H}^\pm_\delta$ and $I\subset \partial D_\pm$. We define
$$
m(z,\omega)=\begin{cases}
m_+(z,\omega)& \text{if $z\in D_+\cap \mathbb{H}_\delta$}\\
\frac{1}{\overline{m_-(z,\omega)}}& \text{if $z\in D_-\cap \mathbb{H}^-_\delta$}
\end{cases},
$$
Note that $m_+(z,\omega)\in \mathbb{H}^2(D_+)$ and  $\overline{m_-(z,\omega)}\in \mathbb{H}^2(D_-)$. By part (3) of Theorem \ref{kotani} and Theorem 2 in \cite{sy}, for almost every $\omega$, $m(z,\omega)$ can be extended analytically to $D_+\cup D_-$.

On the other hand,  we see that $m(E,\omega)_{\omega}$ is a normal
family. Thus for any compact $\mathcal K \in D_+\cup D_-$,
$m(\cdot,\omega)$ is uniformly Lipschitz in $E$.  Namely, there exists
a constant $c=c(\mathcal K)$ that depends only on $\mathcal K$ such that
$$
|m(E_1,\omega)-m(E_2,\omega)|\leq c|E_1-E_2|,\ \ \forall  E_1, E_2\in\mathcal K.
$$
Since $supp(\mu)=\Omega$, we can, for any $\omega\in \Omega,$ pick a sequence $\omega_n\in \Omega$ converging to $\omega$. Then we get a holomorphic function $m(\cdot,\omega)=\lim\limits_{\omega_n\rightarrow \omega}m(\cdot,\omega_n)$ on $D_+\cup D_-$.

Finally, let
$$
(U_E(\omega),V_E(\omega))=\frac{1}{\sqrt{c_E(\omega)}}\begin{pmatrix}
\vec{u}_E(0,\omega)&\vec{v}_E(0,\omega)\\
\vec{u}_E(-1,\omega)&\vec{v}_E(-1,\omega)
\end{pmatrix}C_E(\omega).
$$
Then we have  $(U_E,V_E)\in C^0(\Omega, Sp_{2d\times 2}(\R))$, and
there is $R_E\in C^0(\Omega,SO(2,\R))$, depending analytically on
$E\in \mathbb{C}_{\delta'}\backslash(\R\backslash I)$ for some $\C_{\delta'}\subset\C_\delta$, such that
$$
L_E^f(\omega)(U_E(\omega),V_E(\omega))=(U_E(T\omega), V_E(T\omega))R_E(\omega).
$$

\section{Proof of Theorem \ref{main12}}\label{m12}
We will actually prove a more general theorem. Assume $v,w$ are both even functions. We denote  the non-negative Lyapunov exponents of the complex symplectic cocycle $(\alpha,L_{E,w}^v)$ by $\{\gamma^i(E)\}_{i=1}^\ell$ where $\gamma^1(E)\geq \cdots \geq \gamma^d(E)\geq 0$.
\begin{Theorem}\label{th71}
Given $\alpha\in\R\backslash\Q$ and  an open interval $I\subset \Sigma_{\alpha,v}^w$ , it is impossible that for all $E\in I$,
\begin{enumerate}
\item both $(\alpha,L_{E,v}^w)$ and $(\alpha,L_{E,w}^v)$  are $PH2$ and
\item  $\gamma^\ell(E)=0$.
\end{enumerate}
\end{Theorem}
\begin{pf}
  The proof is via contradiction and an improvement of Corollary
  \ref{C0reducibility}. 
Assume there is an open interval $I\subset \Sigma_{\alpha,v}^w$ such that  for any $E\in I$,
\begin{enumerate}
\item both $(\alpha,L_{E,v}^w)$ and $(\alpha,L_{E,w}^v)$  are $PH2$ and
\item $\gamma^\ell(E)=0$.
\end{enumerate}

\begin{Corollary}\label{Cwreducibility}
If  $\gamma^\ell(E)=0$ and $(\alpha,L_{E,w}^v)$ is $PH2$ for all $E$ in an interval $I\subset \R$, then there exist   $H_E\in C^\omega(\T,Sp_{2l\times 2}(\R))$, amd $\psi_E\in C^\omega(\T,\R)$, depending analytically on $E\in I$ such that
$$
L_{E,w}^v(x)H_E(x)=H_E(x+\alpha)R_{\psi_E(x)}.
$$
\end{Corollary}
\begin{pf}
The proof is based on a Lemma by Avila-Jitomirskaya \cite{aj}.
\begin{Lemma}\label{aj}
Let $W\subset \C$ be a domain, and let $f:W\times \R/\Z\rightarrow\C$ be a continuous function. If $z\rightarrow f(z,w)$ is holomorphic for all $w\in \R/\Z$ and $w\rightarrow f(z,w)$ is analytic for some nonpolar set $z\in W$, then $f$ is analytic.
\end{Lemma}
Note that by Corollary \ref{C0reducibility}, $m(z,x)$ is continuous, $z\rightarrow m(z,x)$ is holomorphic on $D_+\cup D_-$ for any $x\in\R/\Z$ and $x\rightarrow m(z,x)$ is analytic on $\R/\Z$ for any $z\in (D_+\cup D_-)\backslash I$. Thus by Lemma \ref{aj}, $m(z,x)$ is analytic.

As before, we define
$$
C_E(x)=\begin{pmatrix}
0&\frac{|m(E,x)|}{(\Im m(E,x))^{1/2}}\\
-\frac{(\Im m(E,x))^{1/2}}{|m(E,x)|}&\frac{\Re m(E,x)}{|m_+(E,x)|(\Im m(E,x))^{1/2}}
\end{pmatrix}.
$$
Let
$$
H_E(x)=(U_E(x),V_E(x))=\frac{1}{\sqrt{c_E(x)}}\begin{pmatrix}
\vec{u}_E(0,x)&\vec{v}_E(0,x)\\
\vec{u}_E(-1,x)&\vec{v}_E(-1,x)
\end{pmatrix}C_E(x).
$$
Then $H_E\in C^\omega(\T,Sp_{2l\times 2}(\R))$ and there exists
$\psi_E\in C^\omega(\T,\R)$, both $H_E$ and $\psi_E$  depending analytically on $E\in I$, such that
\begin{equation}\label{final}
L_{E,w}^v(x)H_E(x)=H_E(x+\alpha)R_{\psi_E(x)}.
\end{equation}
\end{pf}

Notice that for $\e$ sufficiently small, \eqref{final} also holds for $E+i\e$, thus one has $\gamma^\ell(E+i\e)=\Im \int_\T \psi_{E+i\e}(x)dx$. On the other hand, by Lemma \ref{final1} and (3) of Theorem \ref{keyth}, for almost every $E\in I$,
$$
\lim\limits_{\e\rightarrow 0^+}\frac{\partial \Im \int_\T \psi_{E+i\e}(x)dx}{\partial \e}=\frac{1}{4\ell}\int_\T \frac{\left\|m_+(E+i0,x)\begin{pmatrix}
\vec{u}_E(0,x)\\
\vec{u}_E(-1,x)
\end{pmatrix}+\begin{pmatrix}
\vec{v}_E(0,x)\\
\vec{v}_E(-1,x)
\end{pmatrix}\right\|^2}{\Im \left(m_+(E+i0,x)c_E(x)\right)}dx>0.
$$
Thus $\int_\T\psi_E(x)dx$ is not a constant, so there is $E_0\in I$ with $\int_\T
\psi_{E_0}(x)dx=k\alpha(\mod \Z),$ for some $k\in\Z.$ We now define
$F_{E_0}(x)=H_{E_0}(x)R_{kx}$, and obtain
$$
L_{E_0,w}^v(x)F_{E_0}(x)=F_{E_0}(x+\alpha)R_{\psi_{E_0}(x)-k\alpha}.
$$
Since $\int_\T \left(\psi_{E_0}(x)-k\alpha\right)dx=0$ and $(\alpha, L_{E_0,v}^w)$ is $PH2$, this contradicts  Theorem \ref{contra2}.
\end{pf}

\noindent{\bf Proof of Theorem \ref{main12}}:

Assume  there exists an interval $I\subset \Sigma^1_{v,\alpha}$.
Recall that the associated Schr\"odinger cocycle is denoted
$(\alpha,S_{E}^v),$  the dual cocycle by $(\alpha,L_{E,v})$, and their
non-negative Lyapunov exponents are correspondingly $L(E)$ and $\{\gamma^i(E)\}_{i=1}^d,$ respectively. We distinguish two cases:
\begin{enumerate}
\item There exists $E_0\in I$ such that $L(E_0)>0.$ By continuity of
  Lyapunov exponents then there is $I'\subset I$ such that $L(E)>0$
  and, by assumption, $\bar{\omega}(E)=1$ on $I'$.  Then for any $E\in
  I'$, we have
  \begin{enumerate}
\item by Theorem \ref{dominate1}, both  $(\alpha,S_{E}^v)$ and $(\alpha,L_{E,v})$  are $PH2$
\item by Theorem 1.2 in \cite{gjyz},  $\gamma^d(E)=0$ .
\end{enumerate}
This contradicts  Theorem  \ref{th71}.

\item If no such $E_0$ exists, then for all $E\in I$ we have $L(E)=0.$ Then
\begin{enumerate}
\item by Theorem \ref{dominate1}, both $(\alpha,L_{E,v})$ and $(\alpha,S_{E}^v)$  satisfy the $PH2$ condition;
\item $L(E)=0$.
\end{enumerate}
This again contradicts Theorem  \ref{th71}.\qed
\end{enumerate}



\section{Proof of Theorem \ref{main14}}\label{m14}
Let
$\Sigma^\delta_{\lambda,\alpha}$ be the spectrum of
$H_{\lambda,\alpha,x}^\delta$ given by \eqref{hdelta}.  Let $E_k\in
\Sigma^\delta_{\lambda,\alpha}$ be such that
$2\rho(E_k)=k\alpha(\mod\Z)$, where $\rho$ is the rotation number
defined in Section \ref{secrot}.  We start from the following reducibility theorem.
\begin{Theorem}[\cite{LYZZ}]\label{red}
Assume $\beta(\alpha)=0.$ For any $0<|\lambda|<1$ and  a real
1-periodic trigonometric polynomial $f$, there is
$\delta_0(\lambda,f)>0$ such that if $|\delta|\leq \delta_0$, then
for any $k\in\Z,$ $(\alpha,S_{E_k}^{2\lambda\cos+\delta f})$ is
reducible to a constant matrix $A$ of the form $A=\begin{pmatrix}1&c\\ 0&1\end{pmatrix}$ with $c\in\R,$ in the sense that there exists $B\in C^\omega(\T,PSL(2,\R))$ such that
$$
B_{E_k}^{-1}(x+\alpha)S_{E_k}^{2\lambda\cos+\e V}(x)B_{E_k}(x)=A.
$$

\end{Theorem}
\begin{pf}
Notice that  $\Sigma^\delta_{\lambda,\alpha}\subset [-|\lambda|-4, |\lambda|+4]$  for sufficiently small $\delta$.
By \cite{aj1}, $(\alpha,S_{E}^{2\lambda\cos})$ is almost reducible for
all $E\in\R,$ if $\beta(\alpha)=0$ and $|\lambda|<1$. Since almost
reducibility is an open property \cite{avila0}, for any $E\in
[-|\lambda|-4, |\lambda|+4]$, there exists $\delta_E$ such that if
$|\delta|<\delta_E$, $(\alpha,S_{E}^{2\lambda\cos+\delta f})$ is
almost reducible. By compactness then there is a $\delta_0(\lambda,f)$
such that  $(\alpha,S_{E}^{2\lambda\cos+ \delta f})$ is almost
reducible for all $E\in  [-|\lambda|-4, |\lambda|+4]$ provided
$|\delta|<\delta_0(\lambda,f)$. Theorem \ref{red} is then just a special case of Corollary 5.1 in \cite{LYZZ}.
\end{pf}

\begin{Lemma}\label{degree_psl}
For $A\in C^\omega(T,PSL(2,\R))$, we have
\begin{enumerate}
\item $A\in C^\omega(T,SL(2,\R))$ if $\deg{A}$ is even.
\item $R_{\pi x/2}A\in C^\omega(T,SL(2,\R))$ if $\deg{A}$ is odd.
\end{enumerate}
\end{Lemma}
\begin{pf}
Let $A(x)=A_o(x)+A_e(x)$ where
$$
A_o(x)=\sum\limits_{k\in\Z}a_{k}e^{2k\pi ix},\ \ A_e(x)=\sum\limits_{k\in\Z}a_{2k+1}e^{(2k+1)\pi ix}.
$$
It is obvious that $A_o(x+1)=A_o(x)$, $A_e(x+1)=-A_e(x)$. Thus
$$
A(x+1)=A_o(x)-A_e(x).
$$
Since $A\in C^\omega(T,PSL(2,\R))$, thus there either exists a positive measure set of $x$ such that 
\begin{equation}\label{eqs}
A_o(x)-A_e(x)=A(x+1)=A(x)=A_o(x)+A_e(x),
\end{equation}
or such that
\begin{equation}\label{eqt}
A_o(x)-A_e(x)=A(x+1)=-A(x)=-A_o(x)-A_e(x).
\end{equation}
Thus there exists a positive measure
set of $x$ such that $A_e(x)=0$ or a positive measure
set of $x$ such that $A_o(x)=0, $  so either $A_o(x)\equiv0$ or $A_e(x)\equiv0$ which means either $A(x)\equiv A_e(x)$ or $A(x)\equiv A_o(x)$. 
If $A(x)\equiv A_e(x)$, we have that $A\in C^\omega(T,SL(2,\R))$, so
$\deg{A}=2m$ . If $A(x)\equiv A_o(x)$, then $R_{\pi x/2}A\in C^\omega(T,SL(2,\R))$, so $\deg{A}=2m+1$.\\

\end{pf}
The following proposition follows directly from duality.
\begin{Proposition}\label{dual}
If $(\alpha,S_{E}^{v})$ is reducible to the identity and $\deg{B_E}=k$,
where $B_E$ is the reducibility matrix. We have
\begin{enumerate}
\item If $k=2m+1$, then $L^{2\cos}_{v,\alpha,\alpha/2}$ has two
  linearly independent eigenfunctions $u_{E}$ and $v_{E}$ with
  eigenvalue $E.$.
\item If $k=2m$,  then $L^{2\cos}_{v,\alpha,0}$ has two linearly
  independent eigenfunctions $u_{E}$ and $v_{E}$ with
  eigenvalue $E.$
\end{enumerate}
\end{Proposition}
\begin{pf}
By assumption, there exists $B_{E}\in C^\omega(\T, PSL(2,\R))$ such that
\begin{equation}\label{f1}
B_{E}^{-1}(x+\alpha)S_{E}^v(x)B_{E}(x)=\left(
\begin{array}{ccc}
 1 &  0\\
0 &   1
 \end{array}\right).
\end{equation}
Let $B_E(x)=:\begin{pmatrix}b_E^{11}(x)&b_E^{12}(x)\\ b_E^{21}(x)&b_E^{22}(x)\end{pmatrix}$. It follows that
\begin{equation}\label{f2}
b_E^{11}(x)=b_E^{21}(x+\alpha),
\end{equation}
\begin{equation}\label{f3}
(E-v(x))b_E^{11}(x)-b_E^{21}(x)=b_E^{11}(x+\alpha).
\end{equation}
\eqref{f2} and \eqref{f3} imply that
\begin{align}\label{4.1}
&(E-v(x))b_E^{11}(x)=b_E^{11}(x-\alpha)+b_E^{11}(x+\alpha).
\end{align}

If $k=2m$, by Lemma \ref{degree_psl}, we obtain that $B_{E}\in C^\omega(\T, SL(2,\R))$. Let $u_E(n)=\widehat{b_E^{11}}(n):=\int_\T b_E^{11}(x)e^{2\pi inx}dx$. Taking the Fourier expansion of \eqref{4.1}, we have
\begin{equation}\label{ei1}
\sum\limits_{k=-d}^du_E(n+k)\hat{v}_k+2\cos(2\pi n\alpha)u_E(n)= Eu_E(n),
\end{equation}
i.e.,  $\{u_E(n),n\in\Z\}$ is an eigenfunction of the finite-range operator $L^{2\cos}_{v,\alpha,0}$.

Similarly, if we take $v_E(n)=\widehat{b_E^{21}}(n)$, we have $\{v_E(n),n\in\Z\}$ is also an eigenfunction of the finite-range operator $L^{2\cos}_{v,\alpha,0}$.

If $k=2m+1$, we only need to replace $b_E^{11}(x)$ by $e^{-\pi
  ix}b_E^{11}(x)$, $b_E^{21}(x)$ by $e^{\pi ix}b_E^{21}(x)$. By Lemma
\ref{degree_psl}, $e^{-\pi ix}b_E^{11}(x), e^{\pi ix}b_E^{21}(x)\in
C^\omega(\T,\C)$. The rest of the proof is exactly the same.
\end{pf}

\begin{Theorem}\label{para}
Assume $\beta(\alpha)=0.$ For any $0<|\lambda|<1$ and  a real
1-periodic trigonometric polynomial $f$, there is
$\delta_0(\lambda,f)>0$ such that if $|\delta|\leq \delta_0$, then
for any $k\in\Z,$ $(\alpha,S_{E_k}^{2\lambda\cos+\delta f})$ is
reducible to a parabolic matrix.
\end{Theorem}
\begin{pf}
We only need to prove $c\neq 0$ in Theorem \ref{red}. Otherwise $c=0$
and, by Proposition \ref{dual}, we have that, depending on the parity
of $k$, either $L^{2\cos}_{2\lambda\cos+\delta f,\alpha,\alpha/2}$ or $L^{2\cos}_{2\lambda\cos+ \delta f,\alpha,0}$ has two linearly independent eigenfunctions $u_{E_k}$ and $v_{E_k}$. 
By Corollaries \ref{deigen_pamo} and \ref{Iph2}, operators
$L^{2\cos}_{2\lambda\cos+\delta f,\alpha,x}$ are $PH2$. Thus both cases
contradict Theorem \ref{tsimp}. 
\end{pf}
Theorem \ref{main14} follows from Theorem \ref{para} and a standard Moser-Poschel's argument.\qed

\section{Appendix}
Assume $(\Omega,T)$ is minimal, $f:\Omega\rightarrow \R$ is continuous
and $L_E^f$ is $PH2$. Then for any $E\in\Sigma_f$,  there exist continuous invariant decompositions
$$
\C^{2d}=E^s(\omega)\oplus E^c(\omega)\oplus E^u(\omega).
$$
Moreover, there are $C(E),\delta(E)>\delta'(E)>0$, such that for any $\omega\in\Omega$ and $n\geq 1$, we have
\begin{align}\label{eq10a}
\left\|(L_{E}^{f})_{-n}(\omega)v\right\|>C^{-1}e^{\delta n},\ \  \forall v\in E^s(\omega)\backslash\{0\},\ \ \|v\|=1,
\end{align}
\begin{align}\label{eq11a}
\left\|(L_{E}^{f})_{n}(\omega)u\right\|>C^{-1}e^{\delta n},\ \  \forall u\in E^u(\omega)\backslash\{0\},\ \ \|u\|=1.
\end{align}
\begin{align}\label{eq13a}
\left\|(L_{E}^{f})_{\pm n}(\omega)u\right\|<Ce^{\delta' n},\ \  \forall w\in E^c(\omega)\backslash\{0\},\ \ \|w\|=1.
\end{align}
\begin{align}\label{eq12a}
{\rm dim} E^c(\omega)=2.
\end{align}

For any $\begin{pmatrix}
u(d-1)\\
u(d-2)\\
\vdots\\
u(-d)
\end{pmatrix}\in \C^{2d}$, there exist
$$
\begin{pmatrix}
u^s(d-1)\\
u^s(d-2)\\
\vdots\\
u^s(-d)
\end{pmatrix}\in E^s(\omega), \ \ \begin{pmatrix}
u^c(d-1)\\
u^c(d-2)\\
\vdots\\
u^c(-d)
\end{pmatrix}\in E^c(\omega), \ \
\begin{pmatrix}
u^u(d-1)\\
u^u(d-2)\\
\vdots\\
u^u(-d)
\end{pmatrix}\in E^u(\omega)
$$
such that
$$
\begin{pmatrix}
u(d-1)\\
u(d-2)\\
\vdots\\
u(-d)
\end{pmatrix}=\begin{pmatrix}
u^s(d-1)\\
u^s(d-2)\\
\vdots\\
u^s(-d)
\end{pmatrix}+\begin{pmatrix}
u^c(d-1)\\
u^c(d-2)\\
\vdots\\
u^c(-d)
\end{pmatrix}+\begin{pmatrix}
u^u(d-1)\\
u^u(d-2)\\
\vdots\\
u^u(-d)
\end{pmatrix}.
$$

\begin{Lemma}\label{lee2}
If $\left\|\begin{pmatrix}
u(\pm n+d-1)\\
u(\pm n+d-2)\\
\vdots\\
u(\pm n-d)
\end{pmatrix}\right\|\leq Ce^{\delta'|n|}$ for some $n$ sufficiently large, then
$$
\left\|\begin{pmatrix}
u^s(d-1)\\
u^s(d-2)\\
\vdots\\
u^s(-d)
\end{pmatrix}\right\|,\ \ \left\|\begin{pmatrix}
u^u(d-1)\\
u^u(d-2)\\
\vdots\\
u^u(-d)
\end{pmatrix}\right\|\leq Ce^{-(\delta-\delta') |n|}.
$$
\end{Lemma}
\begin{pf}
By \eqref{eq10a}-\eqref{eq13a}, there exists $N$ such that if $n>N$, we have
\begin{equation}\label{eqq1}\small
\left\|\begin{pmatrix}
u^s(-n+d-1)\\
u^s(-n+d-2)\\
\vdots\\
u^s(-n-d)
\end{pmatrix}\right\|\leq \left\|\begin{pmatrix}
u(-n+d-1)\\
u(-n+d-2)\\
\vdots\\
u(-n-d)
\end{pmatrix}\right\|+\left\|\begin{pmatrix}
u^c(-n+d-1)\\
u^c(-n+d-2)\\
\vdots\\
u^c(-n-d)
\end{pmatrix}\right\|+\left\|\begin{pmatrix}
u^u(-n+d-1)\\
u^u(-n+d-2)\\
\vdots\\
u^u(-n-d)
\end{pmatrix}\right\|\leq Ce^{\delta' |n|},
\end{equation}
\begin{equation}\label{eqq2}\small
\left\|\begin{pmatrix}
u^u(n+d-1)\\
u^u(n+d-2)\\
\vdots\\
u^u(n-d)
\end{pmatrix}\right\|\leq \left\|\begin{pmatrix}
u(n+d-1)\\
u(n+d-2)\\
\vdots\\
u(n-d)
\end{pmatrix}\right\|+\left\|\begin{pmatrix}
u^c(n+d-1)\\
u^c(n+d-2)\\
\vdots\\
u^c(n-d)
\end{pmatrix}\right\|+\left\|\begin{pmatrix}
u^s(n+d-1)\\
u^s(n+d-2)\\
\vdots\\
u^s(n-d)
\end{pmatrix}\right\|\leq Ce^{\delta' |n|},
\end{equation}
By \eqref{eq10a}, \eqref{eq11a}, \eqref{eqq1} and \eqref{eqq2}, we obtain
$$
\left\|\begin{pmatrix}
u^s(d-1)\\
u^s(d-2)\\
\vdots\\
u^s(-d)
\end{pmatrix}\right\|,\ \ \left\|\begin{pmatrix}
u^u(d-1)\\
u^u(d-2)\\
\vdots\\
u^u(-d)
\end{pmatrix}\right\|\leq Ce^{-(\delta-\delta') |n|}.
$$
\end{pf}

\section*{Acknowledgements}
 J. You   was partially supported by NSFC grant (11871286) and
Nankai Zhide Foundation. SJ's work was supported by NSF DMS-2052899,
DMS-2155211, and Simons 681675. She is also grateful to School of
Mathematics at Georgia Institute of Technology where a part of this
work was done.

\end{document}